\documentclass[a4paper,12pt,reqno]{amsproc}

\usepackage{amsmath,amsfonts}

\DeclareMathOperator{\expect}{{\mathbb E}}

\makeatletter
\newcommand{\eqcolon}{\mathrel{\mathord{=}\raise.2\p@\hbox{:}}}
\newcommand{\coloneq}{\mathrel{\raise.2\p@\hbox{:}\mathord{=}}}
\makeatother
\newcommand{\der}{\delta\!}
\newcommand{\XX}{\mathbb{X}}
\newcommand{\RR}{\mathbb{R}}
\newcommand{\NN}{N}
\newcommand{\CC}{\mathcal{C}}
\newcommand{\DD}{\mathcal{D}}
\newcommand{\ZZ}{\mathcal{Z}}
\newcommand{\OCC}{\Omega \mathcal{C}}

\newcommand{\norm}[1]{\lVert #1\rVert}

\newcommand{\abs}[1]{\lvert #1\rvert}
\newcommand{\Abs}[1]{\left\lvert #1\right\rvert}
\newcommand{\wt}{\widetilde}

\newtheorem{definition}{Definition} 
\newtheorem{lemma}{Lemma} 
\newtheorem{proposition}{Proposition} 
\newtheorem{corollary}{Corollary} 
\newtheorem{theorem}{Theorem} 
\newtheorem{remark}{Remark} 
\newtheorem{problem}{Problem} 

\begin{document}
\title{Controlling Rough Paths}
\author{M. Gubinelli}
\address{Dipartimento di  Matematica Applicata ``U.Dini''
Via Bonanno Pisano, 25 bis - 56125 Pisa, ITALIA}
\email{m.gubinelli@dma.unipi.it}

\abstract
 We formulate indefinite integration with respect to an irregular
 function as an algebraic problem which has a unique solution under
 some analytic constraints.
This allows us to define a
 good notion of integral with respect to  irregular paths with
 H\"older exponent greater than $1/3$ (e.g.
 samples of Brownian motion) and study the problem of the existence,
 uniqueness and continuity of solution of
 differential equations driven by such paths.
 We recover 
 Young's theory of integration  and
 the main results of Lyons' theory of rough paths in H\"older topology.  
\endabstract
\maketitle
\section{Introduction}

This work has grown out from the attempt of the author to understand
the integration theory of T.~Lyons~\cite{Lyons,lyonsbook} which gives a meaning
and nice continuity properties to integrals of the form
\begin{equation}
\label{eq:line-integral}
  \int_s^t \langle \varphi(X_u), dX_u\rangle
\end{equation}
where $\varphi$ a differential 1-form on some vector space $V$
and $t \mapsto X_t$ is a path in $V$  not necessarily of bounded
variation. 
From the point of view of Stochastic Analysis Lyons' theory provide a
path-wise formulation of stochastic
integration and stochastic differential equations. The main feature of
this theory is
that a path in a vector space $V$ should not be considered determined 
by a function from an interval $I \subset \RR$ to $V$ but, if this
path is not regular enough,  some additional
information is needed which would play the r\^ole of the iterated integrals for
regular paths: e.g. quantities like the rank two tensor:
\begin{equation}
  \label{eq:area}
\XX^{2,\mu\nu}_{st} = \int_s^t \int_s^u dX^\mu_v dX_u^\nu
\end{equation}
and its generalizations (see the works of K.-T.~Chen~\cite{ChenAll} for 
applications of iterated integrals to Algebraic Geometry and Lie Group Theory).
For irregular paths the r.h.s. of eq.~(\ref{eq:area}) cannot in
general be
understood as a classical Lebesgue-Stieltjes integral however if we have
\emph{any} reasonable definition for this integral then (under some
mild regularity conditions) 
 all the integrals of the form given in
eq.~(\ref{eq:line-integral}) can be defined to depend continuously on
$X,\XX^2$ and $\varphi$ (for suitable topologies). A \emph{rough} path
is the original path together with its iterated integrals of low
degree. The theory can then
be extended to cover the case of more irregular paths (with H\"older exponents less than
$1/3$) by a straightforward but cumbersome generalization of the arguments (the more the
path is irregular the more
iterated integrals are needed to characterize a rough path).
 
\bigskip

With this work we would like provide an alternative formulation of
integration over rough paths which leads to the same results of that
of Lyons' but in some extent is simpler and more
straightforward. We will encounter an algebraic structure
which is interesting by itself and corresponds to a kind of
finite-difference calculus. 
In the original work of  Lyons~\cite{Lyons}  roughness is measured in
$p$-variation norm, instead here we prefer to work with H\"older-like
(semi)norms, in Sec.~\ref{sec:probability} we prove that Brownian
motion satisfy our requirements of regularity. In a recent work 
Friz~\cite{fritz} has established H\"older regularity of Brownian rough
paths (according to  Lyons' theory) and used this result to give an
alternative proof of the support theorem for diffusions. This work has
been extended later by Friz and Victoir~\cite{fritz2} by interpreting Brownian rough
paths as suitable processes on the free nilpotent group of step $2$:
regularity of Brownian rough paths can then be seen as a consequence
of standard H\"older regularity results for stochastic processes on groups. 

We will start by reformulating in Sec.~\ref{sec:algebraic_prelude} the classical integral as the unique
solution of an algebraic problem (adjoined with some analytic
condition to enforce uniqueness) and then generalizing this problem
and building an abstract  tool for its solution.
As a first application we rediscover in Sec.~\ref{sec:young} the integration theory of Young~\cite{young} which was the prelude to
the more deep theory of Lyons. Essentially, Young's theory define the integral
$$
\int_s^t f_u dg_u
$$
when $f$ is $\gamma-$H\"older continuous, $g$ is $\rho-$H\"older
continuous and $\gamma+\rho > 1$ (actually, the original argument was given in term of
$p$-variation norms).
This will be mainly an
exercise to familiarize with the approach before discussing the
integration theory for more irregular paths in
Sec.~\ref{sec:irregular}. We will define integration for a large class
of paths whose increments are controlled by a fixed reference rough
path. This is the main difference with the approach of Lyons.
Next, to illustrate an application of the theory, we discuss  the existence and
uniqueness of solution of ordinary differential equation
driven by irregular paths (Sec.~\ref{sec:ode}). In particular,
sufficient conditions will be given for the existence in the case of
$\gamma$-H\"older paths with $\gamma > 1/3$ which are weaker than
those required to get uniqueness. This point answer a question raised
in Lyons~\cite{Lyons}. In Sec.~\ref{sec:probability} we prove that
Brownian motion and the second iterated integral provided by It\^o or
Stratonovich integration are H\"older regular rough paths for which
the theory outlined above can be applied.
Finally we show how to prove the main results of Lyons' theory
(extension of multiplicative paths and the existence of a map from
almost-multiplicative to multiplicative paths) within this
approach. This last section is intended only for readers already
acquainted with Lyons' theory (extensive accounts are present in literature, see e.g.~\cite{Lyons,lyonsbook}). 

In
appendix~\ref{app:proofs} we collect some lengthy proofs.

\section{Algebraic prelude}
\label{sec:algebraic_prelude}
Consider the following observation. Let $f$ be a bounded continuous
function on $\RR$ and $x$ a function on $\RR$ with continuous first
derivative. Then there exists a unique couple $(a,r)$ with $a \in
C^1(\RR)$, $a_0 = 0$ and $r \in C(\RR^2)$ such that
\begin{equation}
\label{eq:problem}
  f_s(x_t-x_s) = a_t-a_s - r_{st}
\end{equation}
and
\begin{equation}
\label{eq:condition}
  \lim_{t\to s} \frac{|r_{st}|}{|t-s|} = 0.
\end{equation}
This unique couple $(a,r)$ is given by 
$$
a_t = \int_0^t f_u dx_u, \qquad r_{st} = \int_s^t (f_u-f_s) dx_u.
$$ 

The indefinite integral $\int f dx$ is the
unique solution $a$ of the algebraic problem~(\ref{eq:problem}) with the
additional requirement~(\ref{eq:condition}) on the remainder $r$. 
Since the eq.~(\ref{eq:problem}) make sense for arbitrary
functions $f,x$ it is natural to investigate the possible existence
and uniqueness of regular solutions. This will lead to the
generalization of the integral $\int f dx$ for functions $x$ not
necessarily of finite-variation. 

\medskip

\subsection{Framework}
Let $\CC$ the algebra of bounded continuous functions from $\RR$ to
$\RR$
and $\OCC_n$ ($n > 0$) the subset of bounded continuous functions from
$\RR^{n+1}$ to $\RR$ which are zero on the main diagonal where all the
arguments are equal, i.e. $R \in \OCC_n$ implies $R_{t_1\dots t_n}=0$
if $t_1 = t_2 = \cdots = t_n$ . In this paper we will call elements from
$\OCC_n$ (for any $n > 0$) \emph{processes} to distinguish them from
\emph{paths} which are elements of $\CC$.
The vector spaces $\OCC_n$
are  $\CC$-bimodules with left multiplication  $(AB )_{t_1\cdots
  t_{n+1}} \coloneq A_{t_1} B_{t_1\cdots t_{n+1}}$ and right
multiplication  $( B A)_{t_1\cdots t_{n+1}} \coloneq A_{t_{n+1}}
B_{t_1\cdots t_{n+1}}$ for all $(t_1,\dots,t_{n+1}) \in \RR^{n+1}$, $A
\in \CC$ and $B \in \OCC_n$. Moreover if $A \in \OCC_n$ and $B \in
\OCC_m$ it is defined their external product $AB \in \OCC_{m+n-1}$ as
$(AB)_{t_1\cdots t_{m+n-1} } = A_{t_1\cdots t_n} B_{t_{n}\cdots
  t_{n+m-1}}$.
In the following we will write $\OCC$ for
$\OCC_1$. 

The application $\der : \CC \to \OCC$ defined as
\begin{equation}
  \label{eq:derivation}
  (\der A)_{st} \coloneq A_t - A_s 
\end{equation}
is a derivation on $\CC$ since $\der (AB) = A \der B + 
\der A B = B \der A   +  \der B A$.

Let $\OCC^\gamma$ be the subspace of elements $X \in \OCC$ such that
\begin{equation*}
  \|X\|_\gamma \coloneq \sup_{t,s \in \RR^2} \frac{|X_{st}|}{|t-s|^\gamma}
< \infty
\end{equation*}
and let $\CC^\gamma$ be the subspace of the elements $A \in \CC$ such that
$\|\der A\|_\gamma < \infty$.

Define $\OCC_2^{\rho,\gamma}  $  as the
subspace of elements $X$ of $\OCC_2$ such that
\begin{equation*}
  \|X\|_{\rho,\gamma} \coloneq \sup_{(s,u,t) \in
  \RR^3}\frac{|X_{sut}|}{|u-s|^\rho |t-u|^{\gamma}} < \infty
\end{equation*}
Let $\OCC_2^z \coloneq \oplus_{\rho > 0} \OCC_2^{\rho,z-\rho}$: an
element $A \in \OCC_2^z$ is a finite linear combination of elements
$A_i \in \OCC_2^{\rho_i,z-\rho_i} $ for some $\rho_i \in (0,z)$. 

Define the linear operator $\NN : \OCC \to \OCC_2$ as 
\begin{equation*}
  (\NN R)_{sut} \coloneq R_{st}-R_{ut}-R_{su}.
\end{equation*}
and let $\ZZ_2 \coloneq \NN( \OCC )$ and  $\ZZ_2^z \coloneq \OCC^z_2 \cap \ZZ_2$.

We have that $\text{Ker} \NN = \text{Im}\, \der$. Indeed $\NN \der A =
0$ for all $A \in \CC$ and it is easy to see that for each $R \in
\OCC$ such that $\NN R = 0$ we can let $A_t = R_{t0}$ to obtain that
$\der A = R$.

If $F \in \CC$ and $R \in \OCC$ then a straightforward computation
shows that
\begin{equation}
  \label{eq:leibnitz_n}
  \begin{gathered}
N(F R)_{sut} = F_s N(R)_{sut} - \der F_{su}  R_{ut} = (F N(R) - \der F  R )_{sut};
\\
N(R F)_{sut} = F_{t} N(R)_{sut} +   R_{su} \der F_{ut}=  (N(R) F +   R \der F)_{sut}.     
  \end{gathered}
\end{equation}

These equations suggest that the operators $\delta$ and $\NN$ enjoys
remarkable algebraic properties. Indeed they are just the first two
members of a family of linear operators which acts as derivations on
the modules $\OCC_k$, $k=0,1,\dots$ and which can be characterized as
the coboundaries of a cochain complex
which we proceed to define.

\subsection{A cochain complex}
Consider the following chain complex: a \emph{simple} chain of degree $n$ is a
a string 
$[t_1 t_2 \cdots t_n]$ of real numbers and a chain of degree $n$ is
a formal linear combination of simple chains of the same degree with
coefficients in $\mathbb{Z}$. The boundary operator $\partial$ is
defined as
\begin{equation}
  \label{eq:boundary}
 \partial[t_1 \dots t_n] = \sum_{i=1}^n (-1)^{i} [t_1 \cdots \hat t_i \cdots t_n]  
\end{equation}
where $\hat t_i$ means that this element is removed from the
string. For example
\begin{equation*}
  \partial[st] = -[t]+[s], \qquad \partial [sut] = -[su]+[ts]-[ut];
\end{equation*}
It is easy to verify that $\partial \partial = 0$. To this chain
complex is adjoined in a standard way a complex of cochains (which are
linear functionals on chains). A cochain $A$ of degree $n$ is such
that, on simple chains  of degree $n$, act as
\begin{equation*}
 \langle [t_1\cdots t_n],A\rangle = A_{t_1\cdots t_n};
\end{equation*}
The coboundary $\partial^*$ acts on cochains of degree $n$ as
\begin{equation}
  \label{eq:coboundary}
  \begin{split}
(\partial^* A)_{t_1\cdots t_{n+1}} & = \langle[t_1\cdots
t_{n+1}],\partial^* A \rangle = \langle \partial [t_1\cdots
t_{n+1}], A \rangle      
\\ & = \sum_{i=1}^{n+1} (-1)^i
 \langle \partial [t_1\cdots \hat t_i \cdots
t_{n+1}], A \rangle =   \sum_{i=1}^{n+1} (-1)^i A_{t_1\cdots \hat t_i \cdots
t_{n+1}}
  \end{split}
\end{equation}
\emph{e.g.} for cochains $A,B$ of degree $1$ and $2$ respectively, we have
\begin{equation*}
  (\partial^* A)_{st} = A_s - A_t, \qquad (\partial^* B)_{sut} = B_{st}-B_{ut}-B_{su}
\end{equation*}
so that we have natural identifications of $\partial^*$ with $-\delta$ when acting on
$1$-cochains and with $\NN$ when acting on $2$-cochains. We recognize
also that elements of $\OCC_{n-1}$ ($\OCC_0 =
\CC$) are $n$-cochains  and that
we have the following complex of modules:
\begin{equation*}
 0 \rightarrow \RR \rightarrow \CC \stackrel{\partial^*}{\longrightarrow} \OCC \stackrel{\partial^*}{\longrightarrow} \OCC_2 \stackrel{\partial^*}{\longrightarrow} \OCC_3 \rightarrow \cdots 
\end{equation*}
As usual $\partial^* \partial^* = 0$ which means that the image of
$\partial^*|_{\OCC_n}$ is in the kernel of
$\partial^*|_{\OCC_{n+1}}$. Since $\text{Ker}\NN = \text{Im}\der\ $ the
above sequence is exact at $\OCC$. Actually, the sequence is exact at
every $\OCC_n$: let $A$ an $n+1$-cochain such that $\partial^* A =
0$. Let us show that there exists an $n$-cochain $B$ such that $A =
\partial^* B$. Take 
$$
B_{t_1\cdots t_n} = (-1)^{n+1} A_{t_1\cdots t_n s}
$$  
where $s$ is an arbitrary reference point. Then compute
\begin{equation*}
\begin{split}
(\partial^* B)_{t_1\cdots t_{n+1}} & = - B_{t_2\cdots t_{n+1}} +
B_{t_1 \hat t_2 \cdots t_{n+1}} + \cdots + (-1)^{n+1} B_{t_1 \cdots
  t_{n}}
\\ &
= (-1)^{n+1}[- A_{t_2\cdots t_{n+1} s} +
A_{t_1 \hat t_2 \cdots t_{n+1} s} + \cdots + (-1)^{n+1} A_{t_1 \cdots
  t_{n} s}]
\\ & 
= (-1)^{n+1}[(\partial^* A)_{t_1 t_2 \cdots t_{n+1} s} - (-1)^{n+2}
A_{t_1 \cdots t_{n+1}}]
= A_{t_1 \cdots t_{n+1}}.
\end{split}
\end{equation*}

As an immediate corollary we can introduce the operator $\NN_2 : \OCC_2
\to \OCC_3$ such that $N_2 \coloneq \partial^*|_{\OCC_2}$ to
characterize the image of $\NN$ as the kernel of $\NN_2$. Note that,
for example,  $\NN_2$ satisfy a Leibnitz rule: if $A,B \in \OCC_2$,
\begin{equation}
\label{eq:leibnitz_n2}
  \begin{split}
\NN_2 (AB)_{suvt} & = \partial^*(AB)_{suvt} =
-(AB)_{uvt}+(AB)_{svt}-(AB)_{sut}+(AB)_{suv}  
\\ & = 
-A_{uv} B_{vt}+A_{sv}B_{vt}-A_{su}B_{ut}+A_{su}B_{uv}        
\\ & = (\NN A)_{suv} B_{vt} - A_{su} (\NN B)_{uvt}
\\ & = (\NN A B - A \NN B)_{suvt}
  \end{split}
\end{equation}

To understand the relevance of this discussion to our problem let us
reformulate the observation at the beginning of this section as
follows:
\begin{problem}
\label{prob:0}
Given
two paths $F,X \in \CC$ is it possible to find 
a (possibly) unique  decomposition
\begin{equation}
\label{eq:problem00}
  F \der X = \der A - R 
\end{equation}
where 
$A \in \CC$ and $R \in
\OCC$?
\end{problem}

To have uniqueness of this decomposition we should require that $\der A$
should be (in some sense) orthogonal to $R$.  
So we are looking to a canonical decomposition of $\OCC \simeq \der\, \CC
\oplus \mathcal{B}$ where $\mathcal{B}$ is a linear subspace of $\OCC$
which should contain the remained $R$. This decomposition is equivalent to the
possibility of splitting the
short exact sequence
\begin{equation*}
 0 \rightarrow \CC/\RR \stackrel{\der\ }{\longrightarrow} \OCC
 \stackrel{\NN}{\longrightarrow} \ZZ_2 \rightarrow 0.
\end{equation*}
We cannot hope to achieve the splitting in full generality and we must
resort to consider an appropriate linear subspace $\mathcal{E}$ of
$\OCC$ which contains $\der\,\CC$ and for
which we can show that there exists a linear function $\Lambda_{\mathcal{E}} : \NN
\mathcal{E} \to \mathcal{E}$ such that 
$$
\NN \Lambda_{\mathcal{E}} = 1_{\NN \mathcal{E}}.
$$ 
Then $\Lambda_{\mathcal{E}}$ splits the short exact sequence
\begin{equation*}
 0 \rightarrow \mathcal{C}/\mathbb{R} \stackrel{\der\ }{\longrightarrow} \mathcal{E}
 \stackrel{\NN}{\longrightarrow} \NN{\mathcal{E}} \rightarrow 0
\end{equation*}
which implies $\mathcal{E} = \delta \mathcal{C} \oplus \NN \mathcal{E}$.

In this case, if $F \der X \in \mathcal{E}$ we can recover $\der A$ as
\begin{equation}
  \label{eq:split1}
\der A = F \der X - \Lambda_{\mathcal{E}} \NN (F\der X).  
\end{equation}

To identify a subspace $\mathcal{E}$ for which the splitting is possible we
note that $$\text{Im} \delta \cap \OCC^z = \{ 0 \}$$ for all
$z > 1$, indeed, if $X = \der A$ for some $A \in
\CC$ and $X \in \OCC^z$ then  $A \in \CC^z$ which implies $A = \text{const}$ if $z > 1$. 

Then we can reformulate the algebraic characterization of integration at
the beginning of this section as the following problem:
\begin{problem}
\label{prob:1}
Given
two paths $F,X \in \CC$ is it possible to find $A \in \CC$ and $R \in
\OCC^z$ for some $z > 1$ such that  the decomposition
\begin{equation}
\label{eq:problem0}
  F \der X = \der A - R
\end{equation}
holds?  
\end{problem}
Note that if such a decomposition exists then it is
automatically unique since if $F \der X = \der A' - R'$ is another we have that
$R - R' = \der(A-A')$ but since $R-R' \in \OCC^{z} \cap \ker \NN$ we get $R =
R'$ and thus $A = A'$ modulo a constant. 

That Problem~\ref{prob:1} cannot always be solved is clear from the following
consideration: let $F = X$ and apply $\NN$
to both sides of eq.~(\ref{eq:problem0}) to obtain
\begin{equation*}
  \der X_{su} \der X_{ut} =- \NN R_{sut}
\end{equation*}
for all $(s,u,t) \in \RR^3$. Then
\begin{equation*}
  \der X_{st} \der X_{st} =- \NN R_{tst} = R_{st}+R_{st}
\end{equation*}
for all $(t,s) \in \RR^2$. Now, if $R \in \OCC^z$ with $z > 1$ then 
\begin{equation}
\label{eq:counterex0}
  |\der X_{st}| |\der X_{st}| \le 2 \|R\|_z |t-s|^z.
\end{equation}
which implies that $X \in \CC^{z/2}$. So unless this last condition is
fulfilled we cannot solve Problem~(\ref{eq:problem0}) with the
required regularity on $R$.

A sufficient condition for a solution to Problem~\ref{prob:1} to
exists is given by the following result which states sufficient
conditions on $A \in \OCC_2$ for which the algebraic problem
$$
\NN R = A
$$
has a unique solution $R \in \OCC/\der\, \CC$.

\subsection{The main result}
For every $A \in \ZZ_2^z$ with $z > 1$ there exists a unique $R \in
\OCC^z$ such that $\NN R = A$:
\begin{proposition}
\label{prop:main}
If $z > 1$ there exists a unique linear  map $\Lambda:
\ZZ^{z}_2 \to \OCC^z$ such that $N \Lambda  = 1_{\ZZ_2}$ and
 such that for all
$A \in \ZZ^z_2$ we have 
\begin{equation*}
\|\Lambda A\|_{z} \le \frac{1}{2^z-2} \sum_{i=1}^n \|A_i\|_{\rho_i,z-\rho_i}  
\end{equation*}
if
$A = \sum_{i=1}^n A_i $ with
$n \ge 1$, 
$0 < \rho_i < z$ and $A_i \in \OCC_2^{\rho_i,z-\rho_i}$ for
$i=1,\dots,n$. 
\end{proposition}

\subsection{Localization}
If $I \subset J$
denote with $A|_I$ the restriction on $I$ of the function $A$ defined
on $J$.
 
The operator $\Lambda$ is local in the following sense:
\begin{proposition}
\label{prop:local}
If $I \subset \RR$ is an interval and $A,B \in \ZZ_2^z$ with $z > 1$
then
\begin{equation*}
  A |_{I^3} = B |_{I^3} \implies \Lambda A |_{I^2}= \Lambda B|_{I^2}  
\end{equation*}
\end{proposition}
\proof
This follows essentially from the same argument which gives the
uniqueness of $\Lambda$. Indeed if $Q = \Lambda A - \Lambda B$ we have that
$\NN Q = A-B $ which vanish when restricted to $I^2$. So for $(t,s) \in
I^2$, $t \le u \le s$ we have 
$$
Q_{ut}= Q_{st} - Q_{su} 
$$
but since $Q \in \OCC^z$ with $z>1$ we get $Q |_{I^2} = 0$.
\qed

Given an interval $I =[a,b]\subset \RR$ and defining in an obvious way
the corresponding spaces $\CC^\gamma(I)$, $\OCC_n^\gamma(I)$,
etc\dots we can introduce the operator $\Lambda_I : \ZZ_2^z(I) \to
\OCC^z(I)$ as $\Lambda_I A \coloneq \Lambda \tilde A |_{I^2}$ where
$\tilde A \in \ZZ_2^z$ is any arbitrary extension of the element $A
\in \ZZ_2^z(I)$. By the locality of $\Lambda$ any choice of the
extension $\tilde A$ will give the same result, moreover the specific
choice
$  \tilde A_{sut} := A_{\tau(t),\tau(u),\tau(s)} $
where $\tau(t) := (t \wedge b)\vee a$ has the virtue to satisfy the
following bound
\begin{equation*}
  \|\tilde A_i\|_{\rho_i,z-\rho_i} \le \| A_i\|_{\rho_i,z-\rho_i,I} 
\end{equation*}
where $\| \cdot \|_{\rho_i,z-\rho_i,I}$ is the norm on $\OCC_2^z(I)$
and $A = \sum_i A_i$ is a decomposition of $A$ in $\OCC_2^z(I)$ so
that we have
\begin{equation*}
  \|\Lambda_I A\|_{z,I} \le \frac{1}{2^z-2} \sum_i \| A_i\|_{\rho_i,z-\rho_i,I}
\end{equation*}
We will write $\Lambda$ instead of $\Lambda_I$ whenever the interval
$I$ can be deduced  from the  context.

\subsection{Notations}
In the following we will have to deal with tensor products of vector
spaces and we will use the ``physicist'' notation for
tensors.  We will use $V,V_1,V_2,\dots$ to denote vector spaces which
will be always finite dimensional\footnote{In many of the arguments this will
be not necessary, but to handle infinite-dimensional Banach spaces
some care should be excercized in the definition of norms on tensor
products. We prefer to skip this issue for the sake of clarity.}.
 Then, if $V$ is a vector space, $A \in V$ will be denoted
by  $A^\mu$, where $\mu$ is the corresponding vector index (in an
arbitrary but fixed basis), ranging from $1$ to the dimension of $V$,
elements in $V^*$ (the linear dual of $V$) are denoted by $A_\mu$ with
lower indexes,  
elements in $V\otimes V$ will be denoted
by $A^{\mu\nu}$, elements of $V^{\otimes 2}\otimes V^*$ as $A^{\mu\nu}_{\kappa}$, etc\dots 
Summation over repeated indexes is understood whenever not explicitly
stated otherwise: $A_\mu B^\mu$ is the scalar obtained by contracting
$A \in V^*$ with $B \in V$.

Symbols like $\bar \mu, \bar \nu,\dots$ (a bar
over a greek letter) will be vector multi-indexes, i.e. if $\bar \mu =
(\mu_1,\dots,\mu_n)$ then
$
A^{\bar \mu} 
$ is the element $A^{\mu_1,\dots,\mu_n}$ of $V^{\otimes n}$. 
Given two multi-indexes $\bar\mu$ and $\bar\nu$ we
can build another multi-index $\bar\mu\bar\nu$ which is composed of all
the indices of $\bar\mu$ and $\bar \nu$ in sequence. With $|\bar\mu|$
we denote the degree of the multi-index $\bar \mu$, i.e. if $\bar \mu =
(\mu_1,\dots,\mu_n)$ then $|\bar\mu| = n$. Then for example $|\bar \mu
\bar \nu| = |\bar \mu|+|\bar \nu|$. By convention we introduce also the
empty multi-index denoted by $\emptyset$ such that $\bar \mu \emptyset
= \emptyset \bar \mu = \bar \mu$ and $|\emptyset| = 0$.

Symbols like $\CC(V)$, $\OCC(V)$, $\CC(I,V)$, etc\dots (where $I$ is an
interval) will denote paths and
processes with values in the vector space $V$.

Moreover the symbol $K$ will denote arbitrary strictly positive constants, maybe
different from equation to equation and not depending on anything.

\section{Young's theory of integration}

\label{sec:young}
Proposition~\ref{prop:main} allows to solve Problem~\ref{prob:1} when
$F \in \CC^\rho$, $X \in \CC^\gamma$ with $\gamma+\rho > 1$: in this
case $$\NN (F \der X)_{sut} = - \der F_{su} \der X_{ut}$$ so that
$\NN(F \der X) \in \ZZ_2^{\gamma+\rho}$. Then since $\NN(F \der X -
\Lambda \NN (F \der X)) =0 $ there exists a unique $A \in \CC$ (modulo a
constant) such that
\begin{equation*}
  \der A = F \der X - \Lambda \NN (F \der X).
\end{equation*}

\begin{proposition}[Young]
\label{prop:young}
Fix an interval $I \subseteq \mathbb{R}$. If $F \in \CC^\rho(I)$ and $X \in \CC^\gamma(I)$ with $\gamma + \rho >
1$ define
\begin{equation}
  \label{eq:integral_young}
      \int_s^t F_u dX_u \coloneq \left[F \der X - \Lambda \NN(F \der
      X)\right]_{st}, \qquad s,t \in I. 
\end{equation}
Then we have
\begin{equation}
\label{eq:bound_young}
  \left|\int_s^t (F_u-F_s) dX_u \right|\le 
 \frac{1}{2^{\gamma+\rho}-2} |t-s|^{\gamma+\rho} \|F\|_{\rho,I} \|X\|_{\gamma,I}, \qquad s,t \in I.
\end{equation}
\end{proposition}
\proof Is immediate observing that by definition
\begin{equation*}
 \int_s^t (F_u-F_s) dX_s = - [\Lambda\NN (F \der X)]_{st} = 
 [\Lambda(\der F \der X)]_{st}   
\end{equation*}
and using the previous results.
\qed

\medskip
Another justification of this definition of the integral comes
from the the following convergence of discrete sums which also
establish the equivalence of this theory of integration with that of Young.
\begin{corollary}
\label{cor:sums_young}
In the hypothesis of the previous Proposition we have
\begin{equation*}
\int_s^t F_u dX_u = \lim_{|\Pi|\to 0} \sum_{\{t_i\}\in \Pi} F_{t_i}
(X_{t_{i+1}}-X_{t_{i}}) , \qquad s,t \in I   
\end{equation*}
where the limit is taken over  partitions $\Pi = \{t_0,t_1,\dots,t_n\}$ of the interval
$[s,t] \subseteq I$ such that $t_0 = s, t_n = t$, $t_{i+1}>t_i$, $|\Pi| = \sup_i |t_{i+1}-t_i|$.
\end{corollary}
\proof
For any partition $\Pi$ write
\begin{equation*}
  \begin{split}
 S_\Pi &=  \sum_{i=0}^{n-1} F_{t_i} (X_{t_{i+1}}-X_{t_{i}}) 
= \sum_{i=0}^{n-1} (F \der X)_{t_i t_{i+1} } 
 =     \sum_{i=0}^{n-1} (\der A + R)_{t_i t_{i+1} }
  \end{split}
\end{equation*}
with $R \in \OCC^{\gamma+\rho}(I)$ given by $R = \Lambda (\der F \der
X)$ and such that (cfr. Prop.~\ref{prop:young}):
$$
\norm{R}_{\gamma+\rho,I} \le \frac{1}{2^{\gamma+\rho}-2}
\norm{F}_{\rho,I} \norm{X}_{\gamma,I}.
$$
Then
\begin{equation}
  \begin{split}
 S_\Pi & =    A_t-A_s-  \sum_{i=0}^{n-1} R_{t_i t_{i+1}}
= \int_s^t F_u dX_u -  \sum_{i=0}^{n-1} R_{t_i t_{i+1}}.
  \end{split}
\end{equation}
But now, since $\gamma+\rho>1$,
\begin{equation*}
  \sum_{i=0}^{n-1} |R_{t_i t_{i+1}}| \le \norm{R}_{\gamma+\rho,I} 
\sum_{i=0}^{n-1} |t_{i+1}-t_i|^{\gamma+\rho} \le
\norm{R}_{\gamma+\rho,I} |\Pi|^{\gamma+\rho-1} |t-s| \to 0 
\end{equation*}
as $|\Pi|\to 0$.
\qed

\section{More irregular paths}
\label{sec:irregular}
In order to solve Problem~\ref{prob:0}
for a wider class of $F$ and $X$ we are led to dispense with the condition $R
\in \OCC^{z}$ with $z > 1$ and thus loose the uniqueness of the
decomposition: if the couple $(A,R)$ solve the problem, then also $(A+B,R+\der
B)$ solve the problem with a nontrivial $B \in \CC^z$. So our aim is
actually to find a distinguished couple $(A,R)$ which will be
characterized by some additional conditions. 

Up to now we have considered only paths with values in $\mathbb{R}$,
since the general case of vector-valued paths can be easily derived
however in the case of more irregular paths the vector features of
the paths will play a prominent r\"ole so from now on we will consider
paths with valued in (finite-dimensional) Banach spaces $V$,$V_1$,\dots

Let $X \in \CC^\gamma(V)$ a path with values in the Banach space $V$
for some $\gamma > 0$ and \emph{assume} that we are given a tensor
process $\mathbb{X}^2$ in
$\OCC^{2\gamma}(V^{\otimes 2})$ such that
\begin{equation}
\label{eq:two-process}
  \NN (\XX^{2,\mu\nu}) =  \der X^\mu \der X^\nu.
\end{equation}
If $\gamma \le 1/2$ we cannot obtain this process using
prop.~\ref{prop:main} but (as we will see in Sec.~\ref{sec:probability}) there are other natural ways
to build such a process for special paths $X$. We can think at the
arbitrary choice of $\mathbb{X}^2$ among all the possible solutions
(with given regularity $2\gamma$) of eq.~(\ref{eq:two-process}) as a
way to resolve the ambiguity of the decomposition  in
Problem~\ref{prob:0}, since in this case
$$
X^\mu \der X^\nu = \der I^{\mu\nu} - \mathbb{X}^{2,\mu\nu}
$$
and so we are able to integrate any component of $X$ with respect to
each other and we can write
$$
\int_s^t X_u^\mu dX^\nu_u = \der I^{\mu\nu}_{st}
$$
meaning that the integral on the l.h.s. is defined by the r.h.s.,
definition which depends on our choice of $\mathbb{X}^2$. Of course
in this case Corollary~\ref{cor:sums_young} does not hold anymore and discrete
sums of $X\der X$ are not guaranteed to converge to $\int X dX$.

Note that in the scalar case the equation
$$
X \der X = \der I - R
$$
with $X \in \CC^\gamma$
has always a solution given by 
$
I_{t} = X_t^2/2 + \text{const} 
$
for which  
$$
\der I_{st} = \frac{1}{2}X_t^2 - \frac{1}{2} X_s^2 =\frac{1}{2}
X_t(X_t-X_s) + \frac{1}{2} X_s(X_t-X_s) = X_s \der
X_{st} + \frac{1}{2} (\der X_{st})^2  
$$
giving the decomposition $\der I = X \der X + R$ with $R \in
\OCC^{2\gamma}$.
The same argument works for the symmetric part of the two-tensor
$\mathbb{X}^2$: If $X \in \CC^\gamma(V)$ there exists a two-tensor $S \in \OCC^{2\gamma}(V
\otimes V)$ given by
$$
S_{st}^{\mu\nu} = \frac{1}{2}\der X_{st}^\mu \der X_{st}^\nu
$$
for which
$$
\NN  S^{\mu\nu} = \frac{1}{2} (\der X^\mu \der X^\nu + \der X^\nu \der X^\mu) .
$$
of course $S$ is not unique as soon as $\gamma \le 1/2$.

\bigskip
Since one of the feature of the integral we wish to retain is
linearity we must agree that 
if $A$ is a linear application from $V$ to $V$ and $Y^\mu_t =
A^{\mu}_\nu X^\nu_t$ then the integral $\der I = \int Y d X$ must be
such that
\begin{equation*}
  Y^\mu \der X^\nu = A^{\mu}_{\kappa} X^\kappa \der X^\nu = 
\der I^{\mu\nu} - A^{\mu}_{\kappa} \mathbb{X}^{2,\kappa\nu}
\end{equation*}
so
\begin{equation*}
  \der I^{\mu\nu} =   Y^\mu \der X^\nu + A^{\mu}_{\kappa} \mathbb{X}^{2,\kappa\nu}
\end{equation*}
and we have fixed at once the values of all the integrals of linear
functions of the path $X$ w.r.t. $X$. Then consider a path $Y$ which
is only \emph{locally} a linear function of $X$, i.e. such that
\begin{equation}
  \label{eq:expansion}
\der Y^\mu = G^{\mu}_\nu \der X^\nu + Q^\mu
\end{equation}
where $Q$ is a ``remainder'' in $\OCC(V)$ and $G$ is a
path in $\CC(V \otimes V^*)$. 
In order to be able to show that $Y$  is integrable w.r.t. $X$ we must
find a solution $R$ of the equation
\begin{equation*}
   \NN R^{\mu\nu} = \der Y^\mu \der X^\nu.
\end{equation*}
but then, using the local expansion give in eq.~(\ref{eq:expansion}),
\begin{equation*}
  \begin{split}
   \NN R^{\mu\nu} & =  G^{\mu}_\kappa \der X^\kappa\der X^\nu + Q^\mu \der X^\nu
    \\ &
 =  G^{\mu}_\kappa \NN (\XX^{2,\kappa\nu}) + Q^\mu \der X^\nu 
\\ & = \NN(G^{\mu}_\kappa \XX^{2,\kappa\nu}) + \delta G^{\mu}_\kappa 
    \XX^{2,\kappa\nu} + Q^\mu  \der X^\nu  
  \end{split}
\end{equation*}
where we have used eq.~(\ref{eq:leibnitz_n}) (the Leibnitz rule for $\NN$).
To find a solution $R$ is then equivalent to let 
$$
\widetilde R^{\mu\nu} = R^{\mu\nu}   -
G^{\mu}_\kappa \XX^{2,\kappa\nu}
$$ 
and solve
\begin{equation}
\label{eq:rtildexx}
  \NN \widetilde R =  \delta G^{\mu}_\kappa 
    \XX^{2,\kappa\nu} + Q^\mu  \der X^\nu.  
\end{equation}
Sufficient conditions to apply Prop.~\ref{prop:main} to solve
eq.~(\ref{eq:rtildexx}) are that $G \in \CC^{\eta-\gamma}(V\otimes
V^*)$, $Q \in \OCC^\eta(V)$ with $\eta + \gamma = z > 1$. In this case
there exists a unique $\widetilde R \in \OCC^z$
solving~(\ref{eq:rtildexx}) and  we have obtained  the distinguished decomposition
\begin{equation}
\label{eq:decomposition_area}
  Y^\mu \der X^\nu = \der I^{\mu\nu} - G^{\mu}_{\kappa}
  \XX^{2,\kappa\nu} - \widetilde R^{\mu\nu}.
\end{equation}
Note that the path $Y$ lives
a-priori only in $\CC^\gamma$ and this implies that uniqueness of the
solution of Problem~\ref{prob:1} can be achieved only if $\gamma >
1/2$. On the other hand the request that $Y$ can be decomposed as in
eq.~(\ref{eq:expansion}) with prescribed regularity on $G$
and $Q$ has allowed us to show that the  ambiguity in the solution of
Problem~\ref{prob:0} can be reduced to the choice of a process $\XX^2$
satisfying eq.~(\ref{eq:two-process}). Of course if $\gamma > 1/2$
there is only one solution to~(\ref{eq:two-process}) with the
prescribed regularity and the
decomposition~(\ref{eq:decomposition_area}) (into a gradient and a
remainder) coincides with the unique solution of Problem~\ref{prob:1}.

\medskip

Another way to look at this result is to consider the ``non-exact''
differential
$$
F \der X + G \mathbb{X}^2
$$
where $F,G$ are arbitrary paths and ask in which case it admits a
unique decomposition
$$
F \der X + G \mathbb{X}^2 = \der A + R
$$
as a sum of an exact differential plus a remainder term. Of course
to have uniqueness is enough that $R \in \OCC^z$, $z>1$.
Compute
$$
\NN (F \der X + G \mathbb{X}^2) = - \der F \der X - \der G
\mathbb{X}^2 + G \der X \der X = (- \der F + G \der X) \der X - \der G
\mathbb{X}^2
$$
so in order to have $R \in \OCC^z$, $z>1$
condition~(\ref{eq:expansion}) and suitable regularity of $G$ and $Q$, are sufficient to apply
Prop.~\ref{prop:main}.

\subsection{Weakly-controlled paths.} 
The analysis laid out above leads to the following definition.
\begin{definition}
Fix an interval $I \subseteq \mathbb{R}$ and let $X \in \CC^\gamma(I,V)$.
 A path $Z \in \CC^\gamma(I,V)$ is said to be \emph{weakly-controlled by $X$ 
 in $I$ with a remainder of order $\eta$} if it exists a path $Z' \in
 \CC^{\eta-\gamma}(I,V\otimes V^*)$ and a
 process $R_Z \in \OCC^{\eta}(I,V)$ with $\eta > \gamma$ such that
$$
\der Z^\mu =  Z^{\prime\,\mu\nu} \der X^\nu + R^\mu_Z.
$$
If this is the case we will write $(Z,Z') \in \DD^{\gamma,\eta}_X(I,V)$
and we will consider on the linear space $\DD^{\gamma,\eta}_X(I,V)$ the
semi-norm
$$
\|Z\|_{D(X,\gamma,\eta),I} \coloneq
\|Z'\|_{\infty,I}+\|Z'\|_{\eta-\gamma,I}+\|R_Z\|_{\eta,I} + \norm{Z}_{\gamma,I}. 
$$
\end{definition}
(The last contribution is necessary to enforce $Z \in \CC^\gamma(I,V)$ when
$I$ is unbounded).

The decomposition $\der Z^\mu =  Z^{\prime\,\mu\nu} \der X^\nu + R^\mu$ is a-priori not
unique, so a path in $D_{\gamma,\eta}(I,X)$ must be understood as a
pair $(Z,Z')$ since then $R_Z$ is uniquely determined. However we will
often omit to specify $Z'$ when it will be clear from the context.

The term \emph{weakly-controlled} is inspired by the fact that paths
which are solution of differential equations controlled by $X$ (see
Sec.~\ref{sec:ode}) belongs to the class of weakly-controlled paths
(wrt. $X$). In general however, a weakly-controlled path $Z$ is 
uniquely determined knowing $X$ and the ``derivative'' $Z'$ only when
$\eta > 1$. 

Weakly-controlled paths enjoy  a transitivity property:
\begin{lemma}
  \label{lemma:transitivity}
If $Z \in
\DD_Y^{\gamma,\eta}(I,V)$ and $Y \in \DD_X^{\gamma,\sigma}(I,V)$ then $Z \in
\DD_X^{\gamma,\min(\sigma,\eta)}(I,V)$ and
\begin{equation*}
\|(Z,Z')\|_{D(X,\gamma,\delta),I} 
\le K \|Z\|_{D(Y,\gamma,\eta),I} (1+
   \|Y\|_{D(X,\gamma,\sigma),I})(1+ \|X\|_{\gamma,I})
\end{equation*}
where $K$ is some fixed constant.
\end{lemma}
\proof The Proof is in the Appendix,
Sec.~\ref{sec:proof-transitivity}. \qed

\bigskip
Another important property of the class of weakly-controlled paths 
is that it is
stable under smooth maps. Let $C^{n,\delta}(V,V_1)$ the space of $n$-times
differentiable maps from $V$ to the vector space $V_1$ with
$\delta$-H\"older $n$-th derivative and consider the norm
$$
\|\varphi\|_{0,\delta} = \|\varphi\|_\infty + \| \varphi\|_\delta
\qquad \|\varphi\|_{n,\delta} = \|\varphi\|_\infty + \sum_{k=1}^n
\|\partial^k \varphi \|_{\infty} + \|\partial^n \varphi\|_\delta
$$
where $\varphi \in C^{n,\delta}(V,V_1)$, $\partial^k \varphi$ is the
$k$-th derivative of $\varphi$ seen as a  function with values in
$V_1 \otimes V^{* \otimes k}$  and 
$$
\|\varphi \|_\infty = \sup_{x\in V} |\varphi(x)|,
$$
$$
\|\partial^n \varphi\|_\delta = \sup_{x,y \in V} \frac{|\partial^n
  \varphi(x) - \partial^n \varphi(y)|}{|x-y|^\delta}.
$$

\begin{proposition}
\label{prop:functionD}
Let $Y \in \DD^{\gamma,\eta}_X(I,V)$ and  $\varphi \in C^{1,\delta}(V,V_1)$, then the path
 $Z$ such that $Z_t^\mu = \varphi(Y_t)^\mu$ is in $\DD^{\gamma,\sigma}_X(I,V_1)$ with
 $\sigma = \min(\gamma (\delta+1),\eta)$. Its decomposition is
$$
\der Z^\mu = \partial_\nu \varphi(Y)^\mu Y^{\prime\,\nu}_{\kappa} \der X^\kappa + R^\mu_Z 
$$
with $R_Z \in \OCC^\sigma(I,V_1)$ and
 \begin{equation}
\label{eq:bound_functionD}
   \begin{split}
\|Z\|_{D(X,\gamma,\sigma),I} & \le 
K \|\varphi\|_{1,\delta}(
\norm{Y}_{D(X,\gamma,\eta),I} + \norm{Y}^{1+\delta}_{D(X,\gamma,\eta),I} + \norm{Y}^{\sigma/\gamma}_{D(X,\gamma,\eta),I})
   \end{split}
 \end{equation}
and if $\varphi \in C^{2,\delta}(V,V_1)$ we have also
\begin{equation}
\label{eq:bound_lipshitz_functionD}
 \norm{\varphi(Y)-\varphi(\widetilde Y)}_{D(X,\gamma,(1+\delta)\gamma),I} \le C \norm{Y-\widetilde Y}_{D(X,\gamma,(1+\delta)\gamma),I} 
\end{equation}
for $Y,\widetilde Y \in \DD^{\gamma,(1+\delta)\gamma}_X(I,V)$
with 
$$
C =  K
 \norm{\varphi}_{2,\delta} (1+\norm{X}_{\gamma,I})
 (1+\norm{Y}_{D(X,\gamma,(1+\delta)\gamma),I}+\norm{\widetilde Y}_{D(X,\gamma,(1+\delta)\gamma),I})^{1+\delta}.
$$
Moreover if $\widetilde Y \in \DD_{\widetilde
  X}^{\gamma,(1+\delta)\gamma}(I,V)$, $\widetilde Z = \varphi(\widetilde
Y)$ and
$$
\der Y^\mu = Y^{\prime\,\mu}_{\nu} \der X^\nu + R^\mu_Y, \qquad 
\der \widetilde Y^\mu = \widetilde Y^{\prime,\mu}_{\nu} \der \widetilde X^\nu + R^\mu_{\widetilde Y}, \qquad 
$$
$$
\der Z^\mu = Z^{\prime\,\mu}_{\nu} \der X^\nu + R^\mu_Z, \qquad 
\der \widetilde Z^\mu = \widetilde Z^{\prime\,\mu}_{\nu} \der \widetilde X^\nu + R^\mu_{\widetilde Z}, \qquad 
$$
with $Z^{\prime\,\mu}_{\nu,t} = \partial_\kappa \varphi(Y_t)^{\mu}
Y^{\prime\,\kappa}_{\nu,t}$, $\widetilde Z^{\prime\,\mu}_{\nu,t} = \partial_\kappa
\varphi(\widetilde Y_t)^{\mu} \widetilde Y^{\prime\,\kappa}_{\nu,t}$ 
then
\begin{equation}
  \label{eq:function_difference}
\|Z'-\widetilde Z'\|_\infty + \|Z'-\widetilde Z'\|_{\delta\gamma,I} +
\|R_Z-R_{\widetilde Z}\|_{(1+\delta)\gamma,I} +\|Z-\widetilde Z\|_{\gamma,I} \le C
(\|X-\widetilde X\|_{\gamma,I} + \epsilon_I)  
\end{equation}
with
$$
\epsilon_I = \|Y'-\widetilde Y'\|_{\infty,I} + \|Y'-\widetilde Y'\|_{\delta\gamma,I} +
\|R_Y-R_{\widetilde Y}\|_{(1+\delta)\gamma,I} +\|Y-\widetilde Y\|_{\gamma,I}.
$$
\end{proposition}
\proof The proof is given in the Appendix,
Sec.~\ref{sec:proof_of_function_D}. \qed

\subsection{Integration of weakly-controlled paths.} 

Let us given a reference path $X \in \CC^\gamma(I,V)$ and an associated
process $\mathbb{X}^2 \in \OCC^{2\gamma}(I,V\otimes V)$ satisfying the algebraic relationship
\begin{equation}
\label{eq:Hrelation}
\NN \mathbb{X}^{2,\mu\nu}_{sut} = \der X^\mu_{su} \der X^\nu_{ut} \qquad s,u,t \in I.
\end{equation}
Following Lyons we will call the couple $(X,\mathbb{X}^2)$ a \emph{rough
path} (of roughness $1/\gamma$).   

We are going to show that 
 weakly-controlled paths  can be integrated one against the other. 

Take two paths $Z,W$ in $V$ weakly-controlled by $X$ with remainder of order $\eta$.
By an argument similar to that at the beginning of this section we can
obtain a unique decomposition of $Z \der W$ as
\begin{equation*}
  Z^\mu \der W^\nu = \der A^{\mu\nu} - F^{\mu\mu'} G^{\nu\nu'}  \mathbb{X}^{2,\mu'\nu'} + \Lambda \NN(Z^\mu \der W^\nu +
  F^{\mu\mu'} G^{\nu\nu'} \mathbb{X}^{2,\mu'\nu'})
\end{equation*}
and we can state the following Theorem:

\begin{theorem}
 \label{th:rough}
For every $(Z,Z') \in D_X^{\gamma,\eta}(I,V)$ and $(W,W') \in
D_X^{\gamma,\eta}(I,V)$ with $\eta+\gamma = \delta >1$ define
\begin{equation}
  \label{eq:rough-definition}
  \begin{split}
  \int_s^t Z^\mu_{u} dW^\nu_u  := Z^\mu_{s} \der W^\nu_{st} + Z^{\prime\,\mu}_{\mu',s}
  W^{\prime\,\nu}_{\nu',s} \mathbb{X}^{2,\mu'\nu'}_{st} - [\Lambda \NN (Z^\mu \der W^\nu +
   Z^{\prime\,\mu}_{\mu'} W^{\prime\,\nu}_{\nu'}
  \mathbb{X}^{2,\mu'\nu'} ) ]_{st}, \qquad s,t \in I     
  \end{split}
\end{equation}
then this integral extends that defined in prop.~\ref{prop:young} and
the following bound holds:
\begin{equation}
  \label{eq:rough-bound}
  \left| \int_s^t (Z^\mu_u-Z^\mu_s) dW^\nu_u - Z^{\prime\,\mu}_{\mu',s} W^{\prime\,\nu}_{\nu',s} \mathbb{X}^{2,\mu'\nu'}_{st} \right| \le \frac{1}{2^{\delta}-2}
  |t-s|^{\delta} \|(Z,Z')\|_{D(X,\gamma,\eta)} \|(W,W')\|_{D(X,\gamma,\eta)},
\end{equation}
which implies the continuity of the bilinear application
$$
((Z,Z'),(W,W')) \mapsto \left(\int_0^\cdot Z dW,Z W'\right) 
$$
from $\DD_X^{\gamma,\eta}(V)\times \DD_X^{\gamma,\eta}(V)$ to
$\DD_X^{\gamma,\min(2\gamma,\eta)}(V\otimes V)$.
\end{theorem}
\proof
Compute
\begin{equation*}
  \begin{split}
  Q^{\mu\nu}_{sut} & = \NN (Z^\mu \der W^\nu +
  Z^{\prime\,\mu}_{\mu'} W^{\prime\,\nu}_{\nu'} \mathbb{X}^{2,\mu'\nu'})_{sut} 
\\ & = - \der Z^\mu_{su}\der W^\nu_{ut} + (Z^{\prime\,\mu}_{\mu'}
  W^{\prime\,\nu}_{\nu'})_s N \mathbb{X}^{2,\mu'\nu'}_{sut} -
   \der(Z^{\prime\,\mu}_{\mu'} W^{\prime\,\nu}_{\nu'})_{su} \mathbb{X}^{2,\mu'\nu'}_{ut}
\\ &= - Z^{\prime\,\mu}_{\mu',s} \der X^{\mu'}_{su}  
  W^{\prime\,\nu}_{\nu',u} \der X^{\nu'}_{ut} - R^\mu_{Z,su} \der
  W^\nu_{ut} - Z^{\prime\,\mu}_{\mu',s} \der X^{\mu'}_{su}
  R^\nu_{W,ut} \\ & \qquad  - 
   \der(Z^{\prime\,\mu}_{\mu'} W^{\prime\,\nu}_{\nu'})_{su}
  \mathbb{X}^{2,\mu'}_{\nu',ut}
 +  (Z^{\prime\,\mu}_{\mu'} W^{\prime\,\nu}_{\nu'})_s  N \mathbb{X}^{2,\mu'\nu'}_{sut} 
\\ &=  - R^\mu_{Z,su} \der W^\nu_{ut} - Z^{\prime\,\mu}_{\mu',s}\der X^{\mu'}_{su}
  R^\nu_{W,ut}\\ & \qquad  - 
  \der(Z^{\prime\,\mu}_{\mu'} W^{\prime\,\nu}_{\nu'})_{su} \mathbb{X}^{2,\mu'\nu'}_{ut} 
- Z^{\prime\,\mu}_{\mu',s} \der X^{\mu'}_{su}
  \der W^{\prime\,\nu}_{\nu',su} \der X^{\nu'}_{ut}
  \end{split}
\end{equation*}
and observe that all the terms are in $\OCC_2^\delta(I,V^{\otimes 2})$  so that $Q \in
\ZZ_2^\delta(I,V^{\otimes 2})$ is in the domain of $\Lambda$, then
\begin{equation*}
  \begin{split}
 \norm{\Lambda Q}_{\delta,I} & \le \frac{1}{2^\delta-2} \left[ \norm{R_Z}_{\eta,I} \norm{W}_{\gamma,I} +
   \norm{Z'}_{\infty,I} \norm{X}_{\gamma,I} \norm{R_W}_{\eta,I}  \right. 
\\ & \qquad \left. + \norm{\mathbb{X}^2}_{2\gamma,I}
   (\norm{Z'}_{\infty,I} \norm{W'}_{\eta-\gamma,I}
+\norm{W'}_{\infty,I}\norm{Z'}_{\eta-\gamma,I})
   + \norm{Z'}_{\infty,I} \norm{W'}_{\eta-\gamma,I} \norm{X}_{\gamma,I}^2
\right]     
\\ & \le \frac{1}{2^\delta-2} (1+\norm{X}_{\gamma,I}^2  + \norm{\mathbb{X}^2}_{2\gamma,I})
   \norm{(Z,Z')}_{D(X,\gamma,\eta),I} \norm{(W,W')}_{D(X,\gamma,\eta),I}
  \end{split}
\end{equation*}
and the
bound~(\ref{eq:rough-bound}) together with the stated continuity
easily follows. 

To prove that this new integral extends the previous definition note
that when $2\gamma > 1$ eq.~(\ref{eq:Hrelation}) has a
unique solution and since $Z,W \in \CC^\gamma(I,V)$  let
$\tilde A_{st} = \int_s^t Z dW$ where the integral is understood in
the sense of prop.~\ref{prop:young}. Then we have
$$
Z^\mu \der W^\nu = \der \tilde A^{\mu\nu} - \tilde R^{\mu\nu} 
$$
with $\tilde R \in \OCC^{2\gamma}(I,V\otimes V)$, at the same time
$$
Z^\mu \der W^\nu = \der A^{\mu\nu} -  Z^{\prime\,\mu}_{\mu'} W^{\prime\,\nu}_{\nu'}  \mathbb{X}^{2,\mu'\nu'} - R^{\mu\nu}
$$
with $R \in \OCC^{\delta}(I,V^{\otimes 2})$. Comparing
these two expressions and taking into account that $2\gamma > 1$
we get $\der A = \der \tilde A$ and $\tilde R^{\mu\nu} = 
Z^{\prime\,\mu}_{\mu'} W^{\prime\,\nu}_{\nu'}  \mathbb{X}^{2,\mu'\nu'} - R^{\mu\nu}$
proving the equivalence of the two integrals. 
\qed

\bigskip
Note that, in the hypothesis of Th.~\ref{th:rough}, we have
\begin{equation*}
  \mathbb{X}^{2,\mu\nu}_{st} = \int_s^t (X^{\mu}_u-X^{\mu}_s) dX^{\nu}_u.  
\end{equation*}

Even if the notation does not make it explicit it is important to
remark that the integral depends on the rough path $(X,\mathbb{X}^2)$,
however if there is another rough path  $(Y,\mathbb{Y}^2)$ and $X \in \DD_Y^{\gamma,\eta}(I,V)$  we have shown that
$\DD_X^{\gamma,\eta}(I,V)\subseteq \DD_Y^{\gamma,\eta}(I,V)$ (see
Lemma~\ref{lemma:transitivity}) and the integral
defined according to  $(X,\mathbb{X}^2)$ is equal to that defined
according to $(Y,\mathbb{Y}^2)$ if and only if we have
\begin{equation*}
  \mathbb{X}^{2,\mu\nu} = \int_s^t \der X^\mu_{su} dX^\nu_u
\end{equation*}
where this last integral is understood based on $(Y,\mathbb{Y}^2)$.
Necessity is obvious, let us prove sufficiency.
Let the decomposition of $X$ according to $Y$ be
$$
\der X^\mu = A^{\mu}_{\nu} \der Y^\nu + R^\mu_X
$$
and write
$$
\der Z^\mu = Z^{\prime\,\mu}_{\nu} \der X^\nu + R^\mu_Z, 
\qquad \der W^\mu = W^{\prime\,\mu}_{\nu} \der W^\nu + R^\mu_W
$$
then if 
$$
\der I^{\mu\nu}_{st} = \int_s^t Z^\mu d_{(X,\mathbb{X}^2)}W^\nu
$$
is the integral based on $(X,\mathbb{X}^2)$,
$$
\der \widetilde I^{\mu\nu}_{st} = \int_s^t Z^\mu d_{(Y,\mathbb{Y}^2)}W^\nu
$$
the one based on $(Y,\mathbb{Y}^2)$; we have by definition of integral
$$
\der I^{\mu\nu} = Z^{\mu} \der W^\nu + Z^{\prime\,\mu}_{\kappa}
W^{\prime,\nu}_{ \rho} \mathbb{X}^{2,\kappa\rho} + R^{\mu\nu}_I
$$
$$
\der \widetilde I^{\mu\nu} = Z^{\mu} \der W^\nu + Z^{\prime,\mu}_{\kappa}
A^{\kappa}_{\kappa'} W^{\prime\,\nu}_{\rho} A^{\rho}_{\rho'} \mathbb{Y}^{2,\kappa'\rho'} + R^{\mu\nu}_{\widetilde I}
$$
and
$$
\mathbb{X}^{2,\kappa\rho} = A^{\kappa}_{\kappa'} A^{\rho}_{\rho'} \mathbb{Y}^{2,\kappa'\rho'} + R^{\kappa\rho}_{\mathbb{X}^2}
$$
where $R_I, R_{\widetilde I}, R_{\mathbb{X}^2} \in \OCC^{\gamma+\eta}(V^{\otimes 2})$.
Then
\begin{equation*}
  \begin{split}
\der (I^{\mu\nu}-\widetilde I^{\mu\nu}) & = Z^{\prime\,\mu}_{\kappa} W^{\prime,\nu}_{\rho} (\mathbb{X}^{2,\kappa\rho}-A^{\kappa}_{\kappa'} A^{\rho}_{\rho'} \mathbb{Y}^{2,\kappa'\rho'}) + R^{\mu\nu}_I -
R^{\mu\nu}_{\widetilde I}    
\\ & = Z^{\prime\,\mu}_{\kappa} W^{\prime\,\nu}_{\rho}
R_{\mathbb{X}^2}^{\kappa\rho} + R^{\mu\nu}_I - R^{\mu\nu}_{\widetilde I}    
  \end{split}
\end{equation*}
but then $\der(I -\widetilde I) \in \OCC^{\gamma+\eta}(I,V^{\otimes 2})$ with
$\gamma+\eta > 1$ so it must be $\der I = \der \widetilde I$.\qed

Given another rough path $(\wt X, \mathbb{\wt X}^2)$ and paths $\wt W,
\wt Z \in \DD_{\wt X}^{\gamma,\eta}(I,V)$ then
it takes not so much effort to show that
the difference
$$
\Delta_{st} := \int_s^t Z dW - \int_s^t \widetilde Z d\widetilde W
$$
(where the first integral is understood with respect to $(X,\mathbb{X}^2)$
and the second w.r.t. $(\widetilde X, \mathbb{\widetilde X}^2)$)
can be bounded as 
\begin{equation}
\label{eq:continuity}
\|\Delta -Z \der W + \widetilde Z \der \widetilde W+ \widetilde W' \widetilde Z' \mathbb{\widetilde X}^2 - W'  Z'
\mathbb{X}^2\|_{\delta,I} \le \frac{1}{2^z-2} (D_1+D_2+D_3)
\end{equation}
where
$$
D_1 = (1+\norm{X}_{\gamma,I}^2  + \norm{\mathbb{X}^2}_{2\gamma,I})
(\norm{(Z,Z')}_{D(X,\gamma,\eta),I}+\norm{(\widetilde Z,\widetilde
  Z')}_{D(\widetilde X,\gamma,\eta),I}) \epsilon_W
$$
$$
D_2 = (1+\norm{X}_{\gamma,I}^2  + \norm{\mathbb{X}^2}_{2\gamma,I})
(\norm{(W,W')}_{D(X,\gamma,\eta),I}+\norm{(\widetilde W,\widetilde
  W')}_{D(\widetilde X,\gamma,\eta),I}) \epsilon_Z
$$
\begin{equation*}
  \begin{split}
D_3  & = (\norm{(W,W')}_{D(X,\gamma,\eta),I}+\norm{(\widetilde W,\widetilde
  W')}_{D(\widetilde X,\gamma,\eta),I}) 
\\ & \qquad \cdot 
(\norm{(Z,Z')}_{D(X,\gamma,\eta),I}+\norm{(\widetilde Z,\widetilde
  Z')}_{D(\widetilde X,\gamma,\eta),I}) (\norm{X-\widetilde X}_{\gamma,I}  +
\norm{\mathbb{X}^2-\mathbb{\widetilde X}^2}_{2\gamma,I})    
  \end{split}
\end{equation*}
and
\begin{equation*}
  \begin{split}
\epsilon_Z & = \|Z'-\widetilde Z'\|_{\infty,I} + \|Z'-\widetilde
Z'\|_{\eta-\gamma,I} + \|R_Z-\widetilde R_Z\|_{\eta,I}+ \|Z-\widetilde
Z\|_{\gamma,I} 
  \end{split}
\end{equation*}
\begin{equation*}
  \begin{split}
\epsilon_W & = 
\|W'-\widetilde W'\|_{\infty,I} + \|W'-\widetilde
W'\|_{\eta-\gamma,I} + \|R_W-\widetilde R_W\|_{\eta,I}+ \|W-\widetilde
W\|_{\gamma,I}        
  \end{split}
\end{equation*}
so that the integral possess reasonable continuity properties also
with respect to the reference rough path $(X, \mathbb{X}^2)$.

\begin{remark}
It is trivial but cumbersome to generalize the statement of Theorem~\ref{th:rough} in
the case of inhomogeneous degrees of smoothness, i.e. when we have $Z \in \DD_X^{\gamma,\eta}(V)$, $W \in \DD_Y^{\rho,\eta'}(V)$
with $X \in \CC^\gamma(V)$, $Y\in \CC^\rho(V)$ and there is a process $H \in
\OCC^{\gamma+\rho}(V^{\otimes 2})$ which satisfy
$$
N H^{\mu\nu} = \der X^\mu \der Y^\nu.
$$   
In this case the condition to be satisfied in order to be able to
define the integral is $\min(\gamma+\eta',\rho+\eta) = \delta > 1$.  
\end{remark}

As in Sec.~\ref{sec:young} we can give an approximation result of the
integral defined in Theorem~\ref{th:rough} as a limit of sums of increments:
\begin{corollary}
\label{cor:sums_rough}
 In the hypothesis of the previous Proposition we have
\begin{equation*}
\int_s^t Z^\mu_u dW^\nu_u = \lim_{|\Pi|\to 0} \sum_{i=0}^{n-1} \left(Z^\mu_{t_i}
\der W^\nu_{t_{i},t_{i+1}} + Z^{\prime\,\mu}_{\mu',t_i} W^{\prime\,\nu}_{\nu',t_i} \mathbb{X}^{2,\mu'\nu'}_{t_{i},t_{i+1}}\right)    
\end{equation*}
where the limit is taken over  partitions $\Pi = \{t_0,t_1,\dots,t_n\}$ of the interval
$[s,t]$ such that $t_0 = s, t_n = t$, $t_{i+1}>t_i$, $|\Pi| = \sup_i
|t_{i+1}-t_i|$.
\end{corollary}
\proof The proof is analogous to that of Corollary~\ref{cor:sums_young}.\qed

\bigskip

Better bounds can be stated in the case where we are integrating a
path controlled by $X$ against $X$ itself 
\begin{corollary}
\label{cor:betterbounds}
When $W \in \DD_X^{\gamma,\eta}(I,V_1 \otimes V^*)$ the integral
$$
\der A^\mu_{st} = \int_s^t W^\mu_{\nu,u} dX^\nu_u
$$  
belongs to $\DD_{X}^{\gamma,2\gamma}(I,V_1)$ and satisfy
\begin{equation}
\label{eq:boundAsimple}
\norm{\der A - W_\nu \der X^\nu - W^{\prime}_{\nu\kappa}
  \mathbb{X}^{2,\nu\kappa}}_{D(X,\gamma,\eta+\gamma),I} \le
\frac{1}{2^{\eta+\gamma}-2}
(\norm{X}_{\gamma,I}+\norm{\mathbb{X}^2}_{2\gamma,I}) \norm{W}_{D(X,\gamma,\eta),I}  
\end{equation}
Moreover if $(\wt X, \mathbb{\wt X}^2)$ is another rough path and $\wt
W \in \DD_{\wt X}^{\gamma,\eta}(I,V_1\otimes V^*)$ then
$$
\der B^\mu_{st} =  \int_s^t W^\mu_{\nu,u} dX^\nu_u -  \int_s^t \wt W^\mu_{\nu,u} d\wt X^\nu_u
$$
and
$$
\der B^\mu = W^\mu_\nu \der X^\nu - \wt W^\mu_\nu \der \wt X^\nu - W^{\prime\,\mu}_{\nu\kappa} \mathbb{ X}^{2,\nu\kappa} - \wt W^{\prime\,\mu}_{\nu\kappa}
\mathbb{ \wt X}^{2,\nu\kappa} + R^\mu_{B}
$$
with  $R_B$
satisfying the bound
\begin{equation}
\label{eq:betterbound-difference}
\norm{R_B}_{\eta+\gamma,I} \le \frac{1}{2^{\eta+\gamma}-2}\left[
C_{X,I} \epsilon_{W,I} +
(\norm{W}_{D(X,\gamma,\eta),I}+\norm{\wt
  W}_{D(\wt X,\gamma,\eta),I}) \rho_I \right]
\end{equation}
with
$$
\epsilon_{W,I} = \norm{ R_W - R_{\wt W}}_{\eta,I} + \norm{W'-\wt W'}_{\eta-\gamma,I}
$$
and
$$
\rho_I = \norm{X-\wt X}_{\gamma} + \norm{\mathbb{X}^2-\mathbb{\wt X}^2}_{2\gamma,I} 
$$
$$
C_{X,I} =
 \norm{X}_{\gamma,I}+\norm{\mathbb{X}^2}_{2\gamma,I}+
 \norm{\wt X}_{\gamma,I}+\norm{\mathbb{\wt X}^2}_{2\gamma,I} 
$$
\end{corollary}
\proof
The integral path $\der A$ has the following decomposition
$$
\der A^\mu = W^\mu_\nu \der X^\nu + W^{\prime\,\mu}_{\nu\kappa} \mathbb{X}^{2,\nu\kappa} + R^\mu_A
$$
with $R_A$ satisfying
$$
\NN R_A^\mu = \der W^{\prime\,\mu}_{\nu\kappa} \mathbb{X}^{2,\nu\kappa} + R^\mu_{W,\nu} \der X^\nu
$$
then eq.~(\ref{eq:boundAsimple}) follows immediately from the
properties of $\Lambda$. Next, let $\der \wt A = \int \wt W d \wt X$
and
$$
\der \wt A^\mu = \wt W^\mu_\nu \der \wt X^\nu + \wt
W^{\prime\,\mu}_{\nu\kappa} \mathbb{\wt X}^{2,\nu\kappa} + R^\mu_{\wt A}
$$
then
$$
\NN R^\mu_B = \der W^{\prime\,\mu}_{\nu\kappa} \mathbb{X}^{2,\nu\kappa} +
R^\mu_{W,\nu} \der X^\nu - \der \wt W^{\prime\,\mu}_{\nu\kappa}
\mathbb{\wt X}^{2,\nu\kappa} + R^\mu_{\wt W,\nu} \der \wt X^\nu
$$
and
\begin{equation*}
  \begin{split}
\norm{R_B}_{\eta+\gamma,I}& \le \frac{1}{2^{\eta+\gamma}-2}\left[ 
\norm{W'-\wt W'}_{\eta-\gamma,I}
\norm{\mathbb{X}^2}_{2\gamma,I}+\norm{\wt W'}_{\eta-\gamma,I}
\norm{\mathbb{X}^2-\mathbb{\wt X}^2}_{2\gamma,I}
\right. \\ &\qquad  
\left. +
\norm{X-\wt X}_{\gamma,I} \norm{R_W}_{\eta,I}+\norm{\wt X}_{\gamma,I}
\norm{R_W-R_{\wt W}}_{\eta,I}
\right]    
\\ & \le  \frac{1}{2^{\eta+\gamma}-2}\left[ C_{X,I} \epsilon_{W,I} +
  (\norm{W}_{D(X,\gamma,\eta),I} +\norm{\wt
    W}_{D(X,\gamma,\eta),I})\rho_I\right]
  \end{split}
\end{equation*}

\qed

\section{Differential equations driven by paths in $\CC^\gamma(V)$}
\label{sec:ode}
The continuity of the integral defined in eq.~(\ref{eq:integral_young}) allows
to prove existence and uniqueness of solutions of differential
equations driven by paths in $\CC^\gamma(V)$ for $\gamma$ not too small.

Fix an interval $J \subseteq \mathbb{R}$ and
let us given $X \in \CC^\gamma(J,V)$ and a function $\varphi \in
C(V,V \otimes V^*)$. 
A solution $Y$ of the differential
equation
\begin{equation}
\label{eq:diff_eq_young}
dY^\mu_t = \varphi(Y_t)^\mu_\nu dX^\nu_t, \qquad Y_{t_0} = y, \quad
t_0 \in J  
\end{equation}
in $J$ will be a continuous path $Y \in \CC^\gamma(V,J)$ such that
\begin{equation}
\label{eq:integralODE}
Y_t^\mu = y + \int_{t_0}^t \varphi(Y_u)^\mu_\nu dX^\nu_u.  
\end{equation}
for every $t \in J$. If $\gamma > 1/2$ sufficient conditions must be
imposed on $\varphi$ such that the integral
in~(\ref{eq:integralODE}) can be understood in the sense of
prop.~\ref{prop:young}. If $1/3 < \gamma \le 1/2$ the integral must be
understood in the sense of Theorem~\ref{th:rough}. Then in this case we want to
show that, given a driving rough path $(X,\mathbb{X}^2)$ it is possible to find a path $Y \in
\DD_X^{\gamma,2\gamma}(V,J)$ that satisfy eq.~(\ref{eq:integralODE}).

The strategy of the proof will consist in introducing a map $Y \mapsto
G(Y)$ on suitable  paths $Y \in \CC(J,V)$
depending implicitly on $X$ (and eventually on $\mathbb{X}^2$) such that
\begin{equation}
  \label{eq:mapG}
  G(Y)_t = Y_{t_0} + \int_{t_0}^t \varphi(Y_u)^\mu_\nu dX^\nu_u. 
\end{equation}
Existence of solutions will follow from a fixed-point theorem applied
to $G$ acting on a suitable compact and convex  subset of the Banach
space of H\"older continuous functions on $J$
(this require $V$ to be finite dimensional). 
To show uniqueness we will prove that under stronger conditions on
$\varphi$ the map $G$ is locally a strict contraction. 
Next we show also that the It\^o map (in the terminology of
Lyons~\cite{Lyons})  $Y=F(y,\varphi,X) $ (or
$Y=F(y,\varphi,X,\mathbb{X}^2)$) which sends the data of the
differential equation to the corresponding  solution  $Y = G(Y)$, is a
Lipschitz continuous map (in compact intervals $J$) in each of its argument, where
on $X$ and $\mathbb{X}^2$ we are considering the norms of
$\CC^\gamma(J,V)$ and $\OCC^{2\gamma}(J,V^{\otimes 2})$ respectively.

Note that, in analogy with the classical setting, the solution of the differential equation
is ``smooth'' in the sense that it will be of
the form
\begin{equation}
 \label{eq:diffeq-1} 
\der Y = \varphi(Y) \der X + R_Y
\end{equation}
with $R_Y \in \OCC^{z}(V,J)$ with $z>1$ in the case of $\gamma > 1/2$
and of the form
\begin{equation}
  \label{eq:diffeq-2}
\der Y = \varphi(Y) \der X + \partial \varphi(Y) \varphi(Y) \mathbb{X}^2 + Q_Y  \end{equation}
with $R_Y \in \OCC^{z}(V,J)$ with $z>1$ in the case of $1/3 < \gamma
\le 1/2$. 

Natural conditions for existence of solutions will be $\varphi \in
C^{\delta}(V,V\otimes V^*)$ if $\gamma > 1/2$ and $(1+\delta)\gamma >
1$, while $\varphi \in
C^{1,\delta}(V,V\otimes V^*)$ if $1/3 < \gamma \le 1/2$ where $\delta \in
(0,1)$ such that $(2+\delta)\gamma > 1$ while uniqueness will
hold if $\varphi \in C^{1,\delta}(V,V\otimes V^*)$ or $\varphi \in
C^{2,\delta}(V,V\otimes V^*)$ 
respectively with analogous conditions on $\delta$.

\begin{remark}
Another equivalent approach to the definition of a differential equation
in the non-smooth setting is to say that $Y$ solves a differential
equation driven by $X$ if eq.~(\ref{eq:diffeq-1}) or
eq.~(\ref{eq:diffeq-2}) is satisfied with remainders $R_Y$ or $Q_Y$
in $\OCC^z(V)$ for some $z$. This would have the natural meaning of
describing the local dynamical behaviour of $Y_t$ as the parameter $t$
is changed in terms of the control $X$. This point of view has been
explored previously in an unpublished work by A.~M.~Davie~\cite{Davie}
which also gives some examples showing that the conditions on the
vector field $\varphi$ cannot be substantially relaxed.  
\end{remark}

\begin{remark}
In a recent work~\cite{lilyons} Li and Lyons show that, under natural
hypotesis on $\varphi$, the It\^o map $F$ can be differentiated
with respect to the control path $X$ (when extended to a  rough path).
\end{remark}

\subsection{Some preliminary results}
In the proofs of the Propositions below it will be useful the following
comparison of norms which holds for locally H\"older continuous paths:
\begin{lemma}
\label{lemma:improved_bound}
 Let $\eta > \gamma$, $b>a$ then $\OCC^\eta([a,b]) \subseteq
 \OCC^\gamma([a,b])$ and 
 \begin{equation*}
    \norm{X}_{\gamma,[a,b]} \le |b-a|^{\eta-\gamma} \norm{X}_{\eta,[a,b]}
 \end{equation*}
\end{lemma}
for any $X \in \OCC^\eta([a,b])$.
\proof
Easy:
\begin{equation*}
  \norm{X}_{\gamma,[a,b]}  = \sup_{t,s \in [a,b]}
  \frac{\abs{X_{st}}}{|t-s|^\gamma} = \sup_{t,s \in [a,b]}
  \frac{\abs{X_{st}}}{|t-s|^\eta} |t-s|^{\eta-\gamma} \le
  |b-a|^{\eta-\gamma}  \sup_{t,s \in [a,b]}
  \frac{\abs{X_{st}}}{|t-s|^\eta}.
\end{equation*}
\qed

Moreover we will need to patch together local H\"older bounds for
different intervals:
\begin{lemma}
  \label{eq:holder-patching}
Let $I,J$ be two adjacent intervals on $\mathbb{R}$ (i.e. $I \cap J
\neq 0$) then if $X \in \OCC^{\gamma}(I,V)$, $X \in \OCC^{\gamma}(J,V)
$ and $\NN X \in \OCC^{\gamma_1,\gamma_2}(I \cup J,V)$ with $\gamma =
\gamma_1+\gamma_2$, then we
have $X \in \CC^{\gamma}(I \cup J,V)$ with
\begin{equation}
  \label{eq:patchbound}
  \norm{X}_{\gamma,I \cup J} \le      2
  (\norm{X}_{\gamma,I}+\norm{X}_{\gamma,J}) + \norm{\NN X}_{\gamma_1,\gamma_2,I
  \cup J}.
\end{equation}
\end{lemma}
\proof  See the Appendix, Sec.~\ref{sec:proof-holder-patching}. \qed


\subsection{Existence and uniqueness when $\gamma > 1/2$}

 First we will
formulate the results for the case $\gamma > 1/2$ since they are  simpler and
require weaker conditions.

\begin{proposition}[Existence $\gamma > 1/2$]
\label{eq:existence_young} If $\gamma > 1/2$ and $\varphi \in
C^\delta(V,V\otimes V^*)$ with $\delta \in (0,1)$ and $(1+\delta)\gamma > 1$ 
there exists a path $Y \in
\CC^\gamma(V)$ which solves eq.~(\ref{eq:diff_eq_young}) (where the
integral is the one defined in Sec.~\ref{sec:young}).    
\end{proposition}
\proof
Consider an interval $I=[t_0,t_0+T] \subseteq J$, $T>0$ and
note that $W = \varphi(Y)$ is  in $\CC^{\delta\gamma}(I,V\otimes
V^*)$ with
$$
\norm{W}_{\delta\gamma,I}= \norm{\varphi(Y)}_{\delta\gamma,I} \le \norm{\varphi}_{\delta}\norm{Y}_{\gamma,I}^\delta
$$
so that if $(1+\delta)\gamma > 1 $ it is meaningful, according to
Prop.~\ref{prop:young} to consider the application
$\CC^{\gamma}(I,V) \to \CC^{\gamma}(I,V)$ defined as in
eq.~(\ref{eq:mapG}). Moreover the path $Z = G(Y) \in \CC^\gamma(I,V)$
satisfy
\begin{equation*}
\der Z^\mu = \varphi(Y)^{\mu}_{\nu} \der X^\nu + Q^\mu_Z
\end{equation*}
with
$$
\norm{Q_Z}_{(1+\delta)\gamma,I} \le \frac{1}{2^{(1+\delta)\gamma}-2} \norm{X}_{\gamma,I}
\norm{\varphi(Y)}_{\delta \gamma,I} \le  \frac{1}{2^{(1+\delta)\gamma}-2} \norm{\varphi}_{\delta}  \norm{X}_{\gamma,I}\norm{Y}_{\gamma,I}^\delta
$$
then, using Lemma~\ref{lemma:improved_bound},
\begin{equation*}
  \begin{split}
\norm{Z}_{\gamma,I} & \le \norm{\varphi(Y) \der X}_{\gamma,I} + \norm{Q_Z}_{\gamma,I}
\\ & \le \norm{\varphi}_{0,\delta} \norm{X}_{\gamma,I} +
T^{\gamma\delta} \norm{Q_Z}_{(1+\delta)\gamma,I}
\\ & \le K C_{X,I} \norm{\varphi}_{0,\delta} (1+T^{\delta
  \gamma}\norm{Y}^\delta_{\gamma,I}) 
\\ & \le K C_{X,J} \norm{\varphi}_{0,\delta} (1+T^{\delta
  \gamma}\norm{Y}^\delta_{\gamma,I}) 
  \end{split}
\end{equation*}
with
$$
C_{X,I} = \norm{X}_{\gamma,I}
$$
For any $T$ let $A_T > 0$ be the solution to
\begin{equation}
  \label{eq:eqA}
 A_T = K C_{X,J} \norm{\varphi}_{0,\delta} (1+T^{\delta
  \gamma}A_T^\delta).
\end{equation}
Then $\norm{G(Y)}_{\gamma,I} \le A_T$ whenever
$\norm{Y}_{\gamma,I} \le A_T$ and moreover $G(Y)_{t_0} = Y_{t_0}$. 
Then for any $y \in V$, the application $G$ maps the compact and convex
set 
\begin{equation}
  \label{eq:setQ}
Q_{y,[t_0,t_0+T]} = \{ Y \in \CC^\gamma([t_0,t_0+T],V) : Y_{t_0} = y, \norm{Y}_{\gamma,[t_0,t_0+T]}
\le A_T \}   
\end{equation}
into itself. Let us show that $G$ on $Q_{y,[t_0,t_0+T]}$ is at least
H\"older continuous with respect to the norm $\|\cdot\|_\gamma$. This
will allow us to conclude (by the Leray-Schauder-Tychonoff theorem)  the existence of a fixed-point in
$Q_{y,[t_0,t_0+T]}$.
To prove continuity take $Y,\wt Y \in Q_{y,I}$ and denote $\wt Z =
G(\wt Y)$ so that
$$
\der \wt Z^\mu = \varphi(\wt Y)^{\mu}_{\nu} \der X^\nu + \wt Q^\mu_Z
$$
as for $Z = G(Y)$.
Then
\begin{equation}
\label{eq:continuitybound0}
  \begin{split}
\norm{Z-\wt Z}_{\gamma,I} \le \norm{\varphi(Y)-\varphi(\wt
  Y)}_{\infty,I}
\norm{X}_{\gamma,I} + \norm{Q_Z-Q_{\wt Z}}_{\gamma,I}    
  \end{split}
\end{equation}
but now taking $0 < \alpha < 1$ such that $(1+\alpha\delta)\gamma > 1$
\begin{equation*}
\norm{Q_Z-Q_{\wt Z}}_{(1+\alpha\delta)\gamma,I} \le 
\frac{1}{2^{(1+\alpha\delta)\gamma}-2} \norm{X}_{\gamma,I}
\norm{\varphi(Y)-\varphi(\wt Y)}_{\alpha \delta \gamma,I}  
\end{equation*}
To bound $\norm{\varphi(Y)-\varphi(\wt Y)}_{\alpha \delta \gamma,I}$
we interpolate between the following two bounds:
\begin{equation*}
\norm{\varphi(Y)-\varphi(\wt Y)}_{0,I} \le 2
\norm{\varphi(Y)-\varphi(\wt Y)}_{\infty,I} \le 2
\norm{\varphi}_{\delta} \norm{\wt Y - Y}_{\infty,I}^\delta  
\end{equation*}
and
\begin{equation*}
  \begin{split}
\norm{\varphi(Y)-\varphi(\wt Y)}_{\delta \gamma,I} \le
\norm{\varphi(Y)}_{\delta\gamma,I} + \norm{\varphi(\wt
  Y)}_{\delta\gamma,I}
\le \norm{\varphi}_{\delta} (\norm{Y}_{\gamma,I}^\delta+\norm{\wt
  Y}_{\gamma,I}^\delta) \le \norm{\varphi}_\delta 2 A_T^\delta    
  \end{split}
\end{equation*}
obtaining
\begin{equation*}
  \norm{\varphi(Y)-\varphi(\wt Y)}_{\alpha \delta \gamma,I} \le 2
  \norm{\varphi}_\delta \norm{\wt Y - Y}_{\infty,I}^{(1-\alpha)\delta}
  A_T^{\alpha \delta} 
\end{equation*}
Eq.~(\ref{eq:continuitybound0}) becomes
\begin{equation*}
  \begin{split}
\norm{Z-\wt Z}_{\gamma,I} & \le \norm{\varphi(Y)-\varphi(\wt
  Y)}_{\infty,I}
\norm{X}_{\gamma,I} + T^{\alpha\delta \gamma}\norm{Q_Z-Q_{\wt Z}}_{(1+\alpha\delta)\gamma,I}    
\\ &  \le K \norm{\varphi}_\delta  \norm{X}_{\gamma,I} \left[ \norm{Y-\wt
  Y}_{\infty,I}^\delta
 + 
 \norm{\wt Y - Y}_{\infty,I}^{(1-\alpha)\delta}
  A_T^{\alpha \delta} \right]
  \end{split}
\end{equation*}
Since $\norm{Y-\wt Y}_{\infty,I} \le \norm{Y-\wt Y}_{\gamma,I}$
(recall that $T < 1$) we have that $G$ is continuous on $Q_{y,I}$ for
the topology induced by the norm $\norm{\cdot}_{\gamma,I}$ (the paths
all have  a common starting point).
 
Since all these arguments does not
depend on the location of the interval $I$  we can patch together
local solutions to get the existence of a global solution on all $J$. 
\qed

\begin{proposition}[Uniqueness $\gamma > 1/2$]
\label{eq:uniqueness_young}
Assume $\varphi \in C^{1,\delta}(V,V\otimes V^*)$ with $(1+\delta)\gamma > 1$, then 
there exists a unique solution of
eq.~(\ref{eq:diff_eq_young}). The It\^o map $F(y,\varphi,X)$ is
Lipschitz in the sense that satisfy the following bound
$$
\norm{F(y,\varphi,X) - F(\wt y,\wt \varphi,\wt X)}_{\gamma,J} \le M
(\norm{X-\wt X}_{\gamma,J}+\norm{\varphi-\wt\varphi}_{1,\delta}+|y-\wt
y|)
$$
for some constant $M$ depending only on $\norm{X}_{\gamma,J}$,
$\norm{\wt X}_{\gamma,J}$, $\norm{\varphi}_{1,\delta}$, $\norm{\wt
  \varphi}_{1,\delta}$ and $J$.
\end{proposition}
\proof
Let us continue to use the notations of the previous proposition.
Let $Y,\widetilde Y$ be two paths in $\CC^\gamma(J,V)$, and $X,\wt X
\in \CC^\gamma(J,V)$.
Let $W= \varphi(Y)$, $\wt W = \varphi(\wt Y)$, $Z = G(Y)$, $\wt Z =
\wt G(\wt Y)$ where $\wt G$
is the map corresponding to the driving path $\wt X$:
$$
\wt Y \mapsto \wt G(\wt Y)^\mu \coloneq  \wt
Y^\mu_{t_0}+\int_{t_0}^\cdot \varphi(\wt Y_u)^{\mu}_{\nu} d\wt X^\nu_u.
$$
Then
\begin{equation*}
\der \wt Z^\mu = \varphi(\wt Y_s)^{\mu} \der \wt X^\nu
+ Q^\mu_{\wt Z}
\end{equation*}
Introduce the following shorthands:
$$
\epsilon_{Z,I} = \norm{Z-\wt Z}_{\gamma,I}, \quad
\epsilon_{W,I}^* = \norm{W-\wt W}_{\delta\gamma,I}, \quad
\epsilon_{Y,I} = \norm{Y-\wt Y}_{\gamma,I}, \quad
\epsilon_{Y,I}^* = \norm{Y-\wt Y}_{\delta\gamma,I};
$$
$$
\rho_I = \norm{X-\wt X}_{\gamma,I} + |Y_0 - \wt Y_0| + \norm{\varphi-\wt\varphi}_{1,\delta}
$$
$$
C_{X,I} = \norm{X}_{\gamma,I} + \norm{\wt X}_{\gamma,I}
\qquad
C_{Y,I} = \norm{Y}_{\gamma,I} + \norm{\wt Y}_{\gamma,I}.
$$
With these notations, Lemma~\ref{lemma:zetabound-young} states that, when $T <1$ :
\begin{equation}
\label{eq:epsilonZbound-young}
  \begin{split}
  \epsilon_{Z,I} &  \le K  C_{X,I}
  C_{Y,I}^\delta [(1+\norm{\varphi}_{1,\delta})\rho_I +
  \norm{\varphi}_{1,\delta} T^{\gamma\delta} \epsilon_{Y,I}]
  \end{split}
\end{equation}

As we showed before in Prop.~\ref{eq:existence_young} there
exists a constant $A_T$ such that the set $Q_{y,I}\coloneq \{Y \in
C^\gamma(I,V): Y_{t_0} = y, \|Y\|_{\gamma,I} \le
A_T\}$ is invariant under $G$. Take $Y,\widetilde Y
\in Q_{y,I}$ and $X = \wt X$. Then we have $\rho_I = 0$, $C_{Y,I} \le 2 A_T$
and
\begin{equation*}
  \epsilon_{Z,I}   \le K \norm{\varphi}_{1,\delta} C_{X,J} A_T^\delta
   T^{\gamma\delta} \epsilon_{Y,I}.  
\end{equation*}
Choosing $T$ small enough such that $K \norm{\varphi}_{1,\delta}
C_{X,J} A_T^\delta T^{\gamma\delta} = \alpha < 1$ implies  
$$
\norm{G(Y) - G(\wt Y)}_{\gamma,I} = \epsilon_{Z,I} \le \alpha
\norm{Y-\wt Y}_{\gamma,I}.
$$
The map $G$ is then a strict contraction on $Q_{y,I}$ and has a unique
fixed-point. Again, since the estimate does not depend on the
location of $I \subset J$ we can extend the unique solution to all
$J$.
\qed

\subsection{Existence and uniqueness for $\gamma > 1/3$}

\begin{proposition}[Existence $\gamma > 1/3$]
\label{eq:existence_rough} 
If $\gamma > 1/3$ and $\varphi \in C^{1,\delta}(V,V)$ with $(2+\delta)\gamma>1$
there exists a path $Y \in
\DD^{\gamma,2\gamma}_X(V)$ which solves eq.~(\ref{eq:diff_eq_young})
where the integral is understood in the sense of
Theorem~\ref{th:rough} based on the couple $(X,\mathbb{X}^2)$.    
\end{proposition}
\proof
By Prop.~\ref{prop:functionD} for any  $Y \in \DD_{X}^{\gamma,2\gamma}(J,V)$,  
 the path $W = \varphi(Y)$ is in
 $\DD_{X}^{\gamma,(1+\delta)\gamma}(J,V)$ with
\begin{equation}
\label{eq:boundWxx}
\begin{split}
\norm{W}_{D(X,\gamma,(1+\delta)\gamma),I} & =  \norm{\varphi(Y)}_{D(X,\gamma,(1+\delta)\gamma),I} \le  
K \|\varphi\|_{1,\delta}(
\norm{Y}_{*,I} + \norm{Y}^{1+\delta}_{*,I} + \norm{Y}^{2}_{*,I})
\\ & 
\le 3 K \|\varphi\|_{1,\delta} (1+\norm{Y}_{*,I})^2  
\end{split}
\end{equation}
where we introduced the notation  $\|\cdot\|_{*,I} = \|\cdot \|_{D(X,\gamma,2\gamma),I}$.

Then we can integrate $W$
 against $X$ as soon as $(2+\delta)\gamma > 1$ and  define the map
 $G$  as   $G :
\DD^{\gamma,2\gamma}_X(I,V) \to
\DD^{\gamma,2\gamma}_X(I,V)$ with the formula~(\ref{eq:mapG}). Let $Y$
 be a path such that $Y'_{t_0} = \varphi(Y_{t_0})$.

The decomposition of $Z$ (as above  $Z = G(Y)$) reads
\begin{equation*}
\der Z^\mu = Z^{\prime\,\mu}_\nu \der X^\nu + R_Z^\mu = \varphi(Y)^{\mu}_{\nu} \der X^\nu + 
\partial^\kappa \varphi(Y)^{\mu}_{\nu} Y^{\prime\,\kappa}_{\rho} \XX^{2,\nu\rho} + Q^\mu_Z
\end{equation*}
with (use eq.~(\ref{eq:boundAsimple}))
\begin{equation}
\label{eq:Qbound3}
\norm{Q_Z}_{(2+\delta)\gamma,I} \le K C_{X,I}
\norm{\varphi(Y)}_{D(X,\gamma,(1+\delta)\gamma),I} 
\end{equation}
where
$$
C_{X,I} = 1+\norm{X}_{\gamma,I}+\norm{\mathbb{X}^2}_{2\gamma,I}.
$$

Our aim is to bound $Z$ in $\DD_X^{\gamma,2\gamma}(I,V)$. To achieve
this we already have the good bound~(\ref{eq:Qbound3}) for $Q_Z$ so we
need  bounds  for $\norm{\partial_\kappa \varphi(Y)^{\cdot}_{\nu}
 Y^{\prime\,\kappa}_{\rho} \XX^{2,\nu\rho}}_{2\gamma,I}$,
$\norm{\varphi(Y)}_{\gamma,I}$ and $\norm{Z}_{\gamma,I}$. 
To simplify the arguments assume that
$T < 1$ since at the end we will need to take $T$ small anyway.

Let us start with $\norm{\partial_\kappa \varphi(Y)^{\cdot}_{\nu}
 Y^{\prime\,\kappa}_{\rho} \XX^{2,\nu\rho}}_{2\gamma,I}$: 
\begin{equation}
\label{eq:zbound-part1}
  \begin{split}
  \norm{\partial_\kappa \varphi(Y)^{\cdot}_{\nu}
 Y^{\prime\,\kappa}_{\rho} \XX^{2,\nu\rho}}_{2\gamma,I} & \le
 \norm{\partial_\kappa \varphi(Y)^{\cdot}_{\nu}}_{\infty,I}
 \norm{Y^{\prime\,\kappa}_{\rho}}_{\infty,I}
 \norm{\XX^{2,\nu\rho}}_{2\gamma,I} 
\\ & \le \norm{\partial
 \varphi}_{\infty}(|Y'_{t_0}|+T^{\gamma}\norm{Y'}_{\gamma,I})
 \norm{\XX^{2,\nu\rho}}_{2\gamma,I} 
\\ & \le \norm{
 \varphi}_{1,\delta}(\norm{\varphi}_{1,\delta}+T^{\gamma}\norm{Y'}_{\gamma,I})
 \norm{\XX^{2,\nu\rho}}_{2\gamma,I} 
  \end{split}
\end{equation}

Next, using the fact that 
\begin{equation*}
  \begin{split}
\norm{\partial \varphi(Y)}_{\infty,I}& \le |\partial \varphi(Y_{t_0})|+
\norm{\partial \varphi(Y)}_{0,I} 
\\ & \le \norm{\varphi}_{1,\delta} +
T^{\delta\gamma} \norm{\partial \varphi(Y)}_{\delta\gamma,I}
\\ & \le \norm{\varphi}_{1,\delta} +
T^{\delta\gamma} \norm{ \varphi(Y)}_{D(X,\gamma,(1+\delta)\gamma),I}    
  \end{split}
\end{equation*}
obtain 
\begin{equation}
\label{eq:zbound-part2}
  \begin{split}
 \norm{\varphi(Y)}_{\gamma,I} & \le 
\|X\|_{\gamma,I} \norm{\partial \varphi(Y)}_{\infty,I} +
\norm{R_{\varphi(Y)}}_{\gamma,I}  
\\ & \le \norm{\varphi}_{1,\delta}\|X\|_{\gamma,I} +T^{\delta\gamma} (\|X\|_{\gamma,I} \norm{\partial \varphi(Y)}_{D(X,\gamma,(1+\delta)\gamma),I} +
 \norm{R_{\varphi(Y)}}_{(1+\delta)\gamma,I}  ) 
\\ & \le C_{X,I} (\norm{\varphi}_{1,\delta} + T^{\delta\gamma} \norm{\varphi(Y)}_{D(X,\gamma,(1+\delta)\gamma),I}  )
    \end{split}
\end{equation}

To finish consider
\begin{equation}
\label{eq:zbound-part3}
  \begin{split}
 \norm{Z}_{\gamma,I} & \le \norm{Z'\der X}_{\gamma,I} +
 \norm{R_Z}_{\gamma,I} 
   \\ & \le \norm{ \varphi(Y)}_{\infty,I} \norm{X}_{\gamma,I}
 + \norm{\partial \varphi(Y) Y' \mathbb{X}^2}_{2\gamma,I} + \norm{Q_Z}_{2\gamma,I}
  \end{split}
\end{equation}

Putting together the bounds given in
eqs.~(\ref{eq:Qbound3}), (\ref{eq:zbound-part1}), (\ref{eq:zbound-part2})
and eq.~(\ref{eq:zbound-part3}) we get
\begin{equation}
\label{eq:boundstoghether0}
  \begin{split}
\norm{Z}_{*,I} & = \norm{\varphi(Y)}_\infty +
\norm{\varphi(Y)}_{\gamma,I} + \norm{\partial_\kappa \varphi(Y)^{\cdot}_{\nu}
 Y^{\prime\,\kappa}_{\rho} \XX^{2,\nu\rho}}_{2\gamma,I} +
\norm{Q_Z}_{2\gamma,I} + \norm{Z}_{\gamma,I}
\\ & \le
2 (1+\norm{X}_{\gamma,I}) \norm{\varphi(Y)}_\infty +
\norm{\varphi(Y)}_{\gamma,I} 
 + 2 \norm{\partial_\kappa \varphi(Y)^{\cdot}_{\nu}
 Y^{\prime\,\kappa}_{\rho} \XX^{2,\nu}_{\rho}}_{2\gamma,I} +
 2T^{\delta\gamma} \norm{Q_Z}_{(2+\delta)\gamma,I}
\\ & \le K C_{X,I} (\norm{\varphi}_{1,\delta} + \norm{\varphi}_{1,\delta}^2 
+ T^{\delta\gamma} \norm{\varphi}_{1,\delta}\norm{Y}_{*,I}
+ T^{\delta\gamma}\norm{\varphi(Y)}_{D(X,\gamma,(1+\delta)\gamma),I})
  \end{split}
\end{equation}
Eq.~(\ref{eq:boundWxx}) is used to conclude that
\begin{equation}
\label{eq:bound_on_gy}
  \begin{split}
\norm{G(Y)}_{*,I} & \le K \|\varphi\|_{1,\delta}
C_{X,I}(1+\|\varphi\|_{1,\delta}+T^{\delta\gamma} (1+\norm{Y}_{*,I}))^2 
\\ & \le K \|\varphi\|_{1,\delta} C_{X,J}(1+\|\varphi\|_{1,\delta}+T^{\delta\gamma} (1+\norm{Y}_{*,I}))^2
  \end{split}
\end{equation}
There exists $T_*$ such that for any $T < T_*$
the equation
$$
A_T  = K \|\varphi\|_{1,\delta} C_{X,J}(1+\|\varphi\|_{1,\delta}+T^{\delta\gamma} (1+A_T))^2
$$
has at least a solution $A_T > 0$. Then
we get that $\norm{G(Y)}_{*,I} \le A_T$ whenever
$\norm{Y}_{*,I} \le A_T$. Let us now prove that in the set
$$
Q'_{y,I} = \{ Y \in \DD_X^{\gamma,2\gamma}(I,V) : Y_{t_0}=y, Y'_{t_0}
= \varphi(y),
\norm{Y}_{*,I} \le A_T \} 
$$
the map $G$ is continuous (in the topology induced by the
$\norm{\cdot}_{*,I}$ norm).
Take $Y,\wt Y \in Q'_{y,I}$ with $Z = G(Y)$, $\wt Z = G(\wt Y)$ and
\begin{equation*}
\der \wt Z^\mu = \wt Z^{\prime\,\mu}_\nu \der X^\nu + R_{\wt Z}^\mu =
\varphi(\wt Y)^{\mu}_{\nu} \der X^\nu + 
\partial^\kappa \varphi(\wt Y)^{\mu}_{\nu} \wt
Y^{\prime\,\kappa}_{\rho} \XX^{2,\nu\rho} + Q^\mu_{\wt Z}
\end{equation*}
Take $0 < \alpha < 1$ and $(2+\alpha \delta)\gamma > 1$:
a bound similar to Eq.~(\ref{eq:boundstoghether0}) exists for
$\norm{Z-\wt Z}_{*,I}$:
\begin{equation*}
  \begin{split}
\norm{Z-\wt Z}_{*,I} & 
 \le
2 (1+\norm{X}_{\gamma,I}) \norm{\varphi(Y)-\varphi(\wt Y)}_\infty +
\norm{\varphi(Y)-\varphi(\wt Y)}_{\gamma,I} 
\\ & \qquad 
 + 2 \norm{(\partial_\kappa \varphi(Y)^{\cdot}_{\nu}
 Y^{\prime\,\kappa}_{\rho}-\partial_\kappa \varphi(\wt Y)^{\cdot}_{\nu}
 \wt Y^{\prime\,\kappa}_{\rho}) \XX^{2,\nu}_{\rho}}_{2\gamma,I} +
 2 \norm{Q_Z-Q_{\wt Z}}_{(2+\alpha \delta)\gamma,I}
\\ & \le K C_{X,I}\left[
\norm{\varphi(Y)-\varphi(\wt Y)}_{\gamma,I} 
+\norm{\partial \varphi(Y) + \partial \varphi(\wt Y)}_{\infty,I} A_{T} +
 \norm{Y'-\wt Y'}_{\infty,I} \norm{\varphi}_{\infty}
\right]
\\ & \qquad +
 2 \norm{Q_Z-Q_{\wt Z}}_{(2+\alpha \delta)\gamma,I}
  \end{split}
\end{equation*}
when $\norm{Y-\wt Y}_{*,I} \le \varepsilon < 1$ we have
$$
\norm{\varphi(Y)-\varphi(\wt Y)}_{\gamma,I} 
+\norm{\partial \varphi(Y) + \partial \varphi(\wt Y)}_{\infty,I} A_{T} +
 \norm{Y'-\wt Y'}_{\infty,I} \norm{\varphi}_{\infty} \le K
 \norm{\varphi}_{1,\delta} (1+A_T) \varepsilon^\delta 
$$
moreover we can bound $\norm{Q_Z-Q_{\wt Z}}_{(2+\alpha
  \delta)\gamma,I}$ as
$$
\norm{Q_Z-Q_{\wt Z}}_{(2+\alpha
  \delta)\gamma,I} \le \frac{1}{2^{(2+\alpha\delta)\gamma}-2}
C_{X,I} \left[\norm{R_{W}-R_{\wt W}}_{(1+\alpha\delta)\gamma,I} + \norm{\partial
    \varphi(Y) - \partial \varphi(\wt Y)}_{\alpha\delta\gamma,I} \right]
$$
with $W = \varphi(Y)$, $\wt W = \varphi(\wt Y)$.
Both of the terms in the r.h.s. will be bounded by interpolation: the
first between
$$
\norm{R_{W}-R_{\wt W}}_{(1+
  \delta)\gamma,I} \le
\norm{\varphi(Y)}_{D(X,\gamma,(1+\delta)\gamma)}
+ \norm{\varphi(\wt Y)}_{D(X,\gamma,(1+\delta)\gamma)}
$$
and
\begin{equation*}
  \begin{split}
\norm{R_{W}-R_{\wt W}}_{\gamma,I} & = 
\norm{(\der \varphi(Y)-\der \varphi(\wt Y)) - (\partial \varphi(Y) -
  \partial \varphi(\wt Y)) \der X}_{\gamma,I}
\\ &\le \norm{\varphi(Y)- \varphi(\wt Y)}_{\gamma,I} + C_{X,I}
\norm{\partial \varphi(Y) - \partial \varphi(\wt Y)}_{\infty,I}    
\\ & \le \norm{\varphi}_{1,\delta} \varepsilon 
+ C_{X,I} \norm{\varphi}_{1,\delta} \varepsilon^\delta
  \end{split}
\end{equation*}
while the second between
$$
\norm{\partial
    \varphi(Y) - \partial \varphi(\wt Y)}_{\delta\gamma,I} \le 
\norm{\partial
    \varphi(Y) }_{\delta\gamma,I} \le 
+\norm{\partial
    \varphi(\wt Y) }_{\delta\gamma,I} 
$$
and
$$
\norm{\partial
    \varphi(Y) - \partial \varphi(\wt Y)}_{0,I} \le 2
\norm{\partial
    \varphi(Y) - \partial \varphi(\wt Y)}_{\infty,I}
\le \norm{\varphi}_{1,\delta} \norm{Y-\wt Y}^\delta_{\infty,I} \le
\norm{\varphi}_{1,\delta} \varepsilon^{\delta}.
$$
These estimates are enough to conclued that $\norm{Z-\wt Z}_{*,I}$ goes to zero whenever
$\norm{\wt Y - Y}_{*,I}$ does.

Reasoning as in Prop.~\ref{eq:existence_young}
we can prove that a solution exists in $\DD_X^{\gamma,2\gamma}(I,V)$
for any $I \subseteq J$ such that $|I|$ is sufficiently small.
Cover $J$ by a sequence $I_1,\dots,I_n$ of intervals of size $T <
T_*$.
 Patching
together local solutions we have a continuous solution $\overline Y$ defined on
all $J$ with
$$
\der \overline Y = \overline Y' \der X + R_{\overline Y}
$$
where $R_{\overline Y} \in \cup_i \OCC^{2\gamma}(I_i,V)$ and $\overline Y' \in \cup_i \OCC^{\gamma}(I_i,V)$. 
 It remains to prove that $\overline Y \in
\DD_X^{\gamma,2\gamma}(J,V)$. 
 Since the
restriction of $\overline Y$ on $I_i$ is in $Q_{y,I_i}$ for some $y
\in V$  we have that (with abuse of notation)
$\norm{\overline Y}_{*,I_i} \le A_T$ for any $i$.

 Using
Lemma~\ref{eq:holder-patching} iteratively we can obtain that
$$
\norm{\overline Y}_{\gamma,J} \le 2^{n+1} \sup_i \norm{\overline
  Y}_{\gamma,I_i} \le 2^{n+1} A_T
$$
and by the same token
$$
\norm{\overline Y'}_{\gamma,J} \le 2^{n+1} A_T
$$
Next consider $R_{\overline Y}$: write $J_k = \cup_{i=1}^k I_i$ and by
the very same lemma  get ($J_{i+1} = J_i \cup I_{i+1}$)
\begin{equation*}
  \begin{split}
\norm{R_{\overline Y}}_{2\gamma,J_{i+1}} & \le 2
\norm{R_{\overline Y}}_{2\gamma,J_i} + 2\norm{R_{\overline
    Y}}_{2\gamma,I_{i+1}}    + \norm{\der \overline Y'\der
  X}_{\gamma,\gamma,J_{i+1}}    
\\ & \le  2
\norm{R_{\overline Y}}_{2\gamma,J_i} + 2\norm{R_{\overline
    Y}}_{2\gamma,I_{i+1}}    + \norm{\overline Y'}_{\gamma,J} \norm{X}_{\gamma,J}   \end{split}
\end{equation*}
since $\NN R_{\overline Y} = - \der \overline Y'\der X$. By induction
over $i$ we end up with
\begin{equation*}
  \norm{R_{\overline Y}}_{2\gamma,J} \le 2^{n+1} \sup_i \norm{R_{\overline
    Y}}_{2\gamma,I_{i}} + n \norm{\overline Y'}_{\gamma,J}
    \norm{X}_{\gamma,J} 
\le (2^{n+1}+ 2^{2n+2} n) A_T
\end{equation*}
and this is enough to conclude that $\overline Y \in \DD_X^{\gamma,2\gamma}(J,V)$.
\qed

\begin{proposition}[Uniqueness $\gamma > 1/3$]
   If $\gamma > 1/3$ and $\varphi \in C^{2,\delta}(V,V)$ with $(2+\delta)\gamma>1$
there exists a unique path $Y \in
\DD^{\gamma,2\gamma}_X(J,V)$ which solves eq.~(\ref{eq:diff_eq_young})
based on the couple $(X,\mathbb{X}^2)$. Moreover the It\^o map
$F(y,\varphi,X,\mathbb{X}^2)$ is Lipschitz continuous in the following
sense. Let $Y = F(y,\varphi,X,\mathbb{X}^2)$ and $\wt Y = F(\wt y, \wt
\varphi, \wt X, \mathbb{\wt X}^2)$ where  $(X,\mathbb{X}^2)$ and $(\wt
X, \mathbb{\wt X}^2)$ are two rough paths, then defining 
\begin{equation*}
  \epsilon_{Y,I} = \norm{Y'-\wt Y'}_{\infty,I} +  \norm{Y'-\wt
  Y'}_{\gamma,I} + \norm{R_Y-R_{\wt Y}}_{2\gamma,I} + \norm{\varphi-\wt\varphi}_{2,\delta}
\end{equation*}
\begin{equation*}
  \rho_I = |Y_{t_0} - \wt Y_{t_0}| + \norm{X-\wt X}_{\gamma,I} +
  \norm{\mathbb{X}^2-\mathbb{\wt X}^2}_{2\gamma,I}
\end{equation*}
and
$$
C_{X,I} = (1+\norm{X}_{\gamma,I}+\norm{\widetilde
  X}_{\gamma,I}+\norm{\mathbb{X}^2}_{2\gamma,I}+\norm{\mathbb{\wt X}^2}_{2\gamma,I})
$$
$$
C_{Y,I} = (1+\norm{Y}_{*,I}+\norm{\widetilde Y}_{*,I}).
$$
we have that there exists a constant $M$ depending only on $C_{X,J}$,
$C_{Y,J}$, $\norm{\varphi}_{2,\delta}$ and $\norm{\wt \varphi}_{2,\delta}$ such that
$$
\epsilon_{Y,J} \le M \rho_{J}.
$$
\end{proposition}

\proof
The strategy will be the same as in the proof of
Prop.~\ref{eq:uniqueness_young}.
Take two paths $Y,\widetilde Y \in \DD_X^{\gamma,2\gamma}(J,V)$ and
let as above 
$Z = G(Y)$, $\wt Z = \wt G (\wt Y)$. Write
the decomposition for each of the paths $Y, \wt Y, Z,\wt Z$ as
$$
\der Y^\mu = Y^{\prime\,\mu}_{\nu} \der X^\nu + R^\mu_Y, \qquad 
\der \widetilde Y^\mu = \widetilde Y^{\prime\,\mu}_{\nu} \der \widetilde X^\nu + R^\mu_{\widetilde Y}, \qquad 
$$
and
$$
\der Z = Z' \der X + R_Z = \varphi(Y) \der X + \partial \varphi(Y) \mathbb{X}^2 + Q_Z
$$
$$
\der \wt Z = \wt Z' \der \wt X + R_{\wt Z} = \wt \varphi(\wt Y) \der
\wt X + \partial \wt \varphi(\wt Y)
\mathbb{\wt X}^2 + Q_{\wt Z}
$$

The key point is to bound $\epsilon_{Z,I}$ defined as
\begin{equation*}
  \begin{split}
\epsilon_{Z,I} & = \norm{\varphi(Y)-\wt \varphi(\wt Y)}_{\infty,I} +
\norm{\varphi(Y)-\wt \varphi(\wt Y)}_{\gamma,I} + \norm{R_Z-R_{\wt
    Z}}_{2\gamma,I}    
  \end{split}
\end{equation*}
and the result of Lemma~\ref{lemma:zetabound-rough} (in the Appendix) tells us that, when
$T < 1$, $\epsilon_{Z,I}$ can be bounded by
\begin{equation}
\label{eq:zetabound}
   \epsilon_{Z,I} \le K  [(1+\norm{\varphi}_{2,\delta})C_{X,I}^2 C_{Y,I}^3
   \rho_I + \norm{\varphi}_{2,\delta} T^{\delta\gamma} C_{X,I}^3 C_{Y,I}^2 \epsilon_{Y,I}].
 \end{equation}

Taking $Y_0 = \wt Y_0$, $\wt X = X$, $\mathbb{\wt X}^2 =
\mathbb{X}^2$ and $\varphi = \wt \varphi$ we have $\rho_I = \rho_J
= 0$. As shown in the proof of Prop.~\ref{eq:existence_rough} if $T < T_*$  for any $y \in V$ there exists a set
$Q_{y,I} \subset \DD_{X}^{\gamma,2\gamma}(I,V)$ invariant under $G$.
Moreover if  $Y,\wt Y \in Q_{y,I}$ for some $y$ then
$\norm{Y}_{*,I} \le A_T$, $\norm{\wt Y}_{*,I} \le A_T$ and
letting
$$
\bar C_{Y,T} = 1+2 A_T
$$
we can rewrite eq.~(\ref{eq:zetabound}) as
$$
  \epsilon_{Z,I} \le K  \norm{\varphi}_{2,\delta} T^{\delta\gamma}
  C_{X,J}^3 \bar C_{Y,T}^2 \epsilon_{Y,I}
$$
so choosing $ T$ small enough such that 
\begin{equation}
  \label{eq:tcritical}
  T^{\delta\gamma}
C_{X,J}^3 \bar C_{Y,T}^2 = \alpha < 1 
\end{equation}
we have 
$$
\|G(Y)-G(\wt Y)\|_{*,I} = \epsilon_{Z,I} \le \alpha \epsilon_{Y,I} =
\alpha  \|Y-\wt Y\|_{*,I}.
$$
Then $G$ is a strict contraction in $\DD_X^{\gamma,2\gamma}(I,V)$ and
thus has a unique fixed-point. Again, patching together local solutions we
get a global one defined on all $J$ and belonging to
$\DD_X^{\gamma,2\gamma}(J,V)$.

Now let us discuss the continuity of the It\^o map
$F(y,\varphi,X,\mathbb{X}^2)$.  Let $Y, \wt Y$ be the solutions based
on $(X,\mathbb{X}^2)$ and $(\wt X,\mathbb{\wt X}^2)$ respectively. We have $Y
= G(Y) = Z$, $\wt Y
= \wt G(\wt Y) = \wt Z$ so that $\epsilon_{Z,I} = \epsilon_{Y,I}$ for any
interval $I \subset J$ and we   can use
eq.~(\ref{eq:zetabound}) to write
$$
\epsilon_{Y,I} = \epsilon_{Z,I} \le
K  [(1+\norm{\varphi}_{2,\delta})C_{X,I}^2 C_{Y,I}^3 \rho_I + \norm{\varphi}_{2,\delta}T^{\delta\gamma}
C_{X,I}^3 C_{Y,I}^2 \epsilon_{Y,I}].
$$
Fix $T$ small enough for~(\ref{eq:tcritical}) to hold so that
$$
\epsilon_{Y,I} \le (1-\alpha)^{-1} K  (1+\norm{\varphi}_{2,\delta})
C_{X,J}^2 C_{Y,J}^3 \rho_I = M_1 \rho_I
$$
Cover $J$ with intervals $I_1,\dots,I_n$ of width $T$ and let $J_k =
\cup_{i=1}^k I_k$ with $J_n = J$.

To patch together the bounds for different $I_i$ into a global bound
for $\epsilon_{Y,J}$ we use again
 Lemma~\ref{eq:holder-patching}  to estimate
\begin{equation*}
  \begin{split}
 \norm{R_Y-R_{\wt Y}}_{2\gamma,J_{i+1}} & \le \norm{R_Y-R_{\wt
 Y}}_{2\gamma,J_{i}}+\norm{R_Y-R_{\wt Y}}_{2\gamma,I_{i+1}}+\norm{\der
 Y' \der X - \der \wt Y' \der \wt X}_{\gamma,\gamma,J_{i+1}}   
\\ & \le  \norm{R_Y-R_{\wt
 Y}}_{2\gamma,J_{i}}+\norm{R_Y-R_{\wt Y}}_{2\gamma,I_{i+1}}
\\ & \qquad + \norm{
 Y'-\wt Y'}_{2\gamma,J_{i+1}} \norm{  X}_{\gamma,J} + \norm{ \wt
 Y'}_{\gamma,J} \norm{ X- \wt X}_{\gamma,J_{i+1}} 
  \end{split}
\end{equation*}
then we obtain easily that
\begin{equation*}
  \epsilon_{Y,J_{i+1}} \le C_{X,J}
  (\epsilon_{Y,J_{i}}+\epsilon_{Y,I_{i+1}} ) + C_{Y,J} \rho_J.
\end{equation*}
Proceeding by induction we get
\begin{equation*}
  \begin{split}
  \epsilon_{Y,J_{n}} & \le (C_{X,J} n + \sum_{k=1}^n C_{X,J}^k) \sup_i \epsilon_{Y,I_{i}}
  + n C_{Y,J} \rho_J \\ & \le [ 2 \sum_{k=1}^n C_{X,J}^k
   M_1 + n
  C_{Y,J}] \rho_J    
  \end{split}
\end{equation*}
which implies that  there exists a constant $M$ depending only on $C_{X,J}$,
$C_{Y,J}$, $\norm{\varphi}_{2,\delta}$ such that
$$
\epsilon_{Y,J} \le M \rho_{J}.
$$

\qed

\section{Some probability}
\label{sec:probability}
So far we have developed our arguments using only analytic and
algebraic properties of paths. In this section we  show how
probability theory provides concrete examples of non-smooth paths for which the
theory outlined above applies.  
 
Let $(\Omega,\mathcal{F},\mathbb{P})$ a probability space where is
defined a standard Brownian motion $X$ with values in $V=\mathbb{R}^n$
(endowed with the Euclidean scalar product).
 It is well known that $X$
is almost surely locally H\"older continuous for any exponent $\gamma < 1/2$,
so that we can fix $\gamma < 1/2$ and choose a version of $X$ living
in  $ \CC^\gamma(I,V)$ on any bounded interval $I$.   
In this case solutions $\mathbb{X}^2$ of eq.~(\ref{eq:two-process})
can be obtained by stochastic integration: let
$$
W^{\mu\nu}_{\text{It\^o},st} \coloneq \int_s^t (X^\mu_u-X^\mu_s) \hat
d X^\nu_u
$$ 
where the hat indicates that the integral is understood in It\^o's
sense with respect to the 
forward filtration $\mathcal{F}_t = \sigma(X_s; s\le t)$. Then it is
easy to show that, for any $s,u,t \in \mathbb{R}$
\begin{equation}
  \label{eq:ito1}
W^{\mu\nu}_{\text{It\^o},st} - W^{\mu\nu}_{\text{It\^o},su} -
W^{\mu\nu}_{\text{It\^o},ut} = (X^\mu_u - X^\mu_s)(X^\nu_t - X^\nu_u)  
\end{equation}
which means that
$$
\NN W^{\mu\nu}_{\text{It\^o}} = \der X^\mu \der X^\nu
$$
then we can choose a continuous version $\mathbb{X}^2_{\text{It\^o}}$
of $(t,s) \mapsto
W_{\text{It\^o,st}}$ for which  eq.~(\ref{eq:ito1}) holds
a.s. for all $t,u,s \in \mathbb{R}$.
It remains to show that $\mathbb{X}^2_{\text{It\^o}} \in
\OCC^{2\gamma}(I,V^{\otimes 2})$ (for any $\gamma <1/2$ and bounded interval $I$).

To prove this result we will develop a small variation on a well known argument
first introduced by Garsia, Rodemich and Rumsey (cfr.~\cite{Rosinski,stroock}) to control H\"older-like seminorms of
continuous stochastic processes
with a corresponding integral norm.

Fix an interval $T  \subset \mathbb{R}$.  A Young function
$\psi$ on $\mathbb{R}^+$ is an increasing, convex function such
that $\psi(0)=0$.

\begin{lemma}
\label{lemma:besov}
For any process $R \in \OCC(T)$ 
let
\begin{equation*}
  U = \int_{T \times T} \psi\left(\frac{|R_{st}|}{p(|t-s|/4)} 
  \right) dt\,ds
\end{equation*}
where $p: \mathbb{R}^+ \to \mathbb{R^+}$ is an increasing function
with $p(0)=0$ and $\psi$ is a Young function. Assume there exists a constant $C$ such that
\begin{equation}
\label{eq:ext-bound-n}
 \sup_{(u,v,r) \in [s,t]^3} |\NN R_{u v r}| \le \psi^{-1}\left(
 \frac{C}{|t-s|^2}\right) p(|t-s|/4), 
\end{equation}
for any couple $s<t$ such that $[s,t] \subset T$. Then
\begin{equation}
\label{eq:control-besov}
  |R_{st}| \le 16 \int_0^{|t-s|}
   \left[\psi^{-1}\left(\frac{U}{ r^2}\right)+\psi^{-1}\left(\frac{C}{ r^2}\right) \right]dp(r)
\end{equation}
for any $s,t \in T$.
\end{lemma}
\proof See the Appendix, Sec.~\ref{sec:proof_besov}.\qed

\begin{remark}
 Lemma~\ref{lemma:besov} reduces to well known results in the case
 $\NN R = 0$ since we can take $C=0$. Condition~(\ref{eq:ext-bound-n})
 is not very satisfying and we conjecture that an integral control
 over $\NN R$ would suffice to obtain~(\ref{eq:control-besov}).
However in its
 current formulation it is enough to prove the following useful corollary.     
\end{remark}

\begin{corollary}
For any $\gamma > 0$ and $p \ge 1$  there exists a constant
$C$ such that for any $R \in \OCC$ 
\begin{equation}
\label{eq:generalboundxx}
\|R\|_{\gamma,T} \le C (U_{\gamma+2/p,p}(R,T)+\|\NN R\|_{\gamma,T}).
\end{equation}
where
\begin{equation*}
 U_{\gamma,p}(R,T) = \left[ \int_{T \times T}
 \left(\frac{|R_{s t}|}{|t-s|^\gamma}\right)^p dt ds \right]^{1/p}.
\end{equation*}
\end{corollary}
\proof
In the previous proposition take $\psi(x) = x^p$, $p(x) =
x^{\gamma+2/p}$; 
the conclusion easily follows.
\qed

\medskip

In the case of $\mathbb{X}^2$ we have, fixed $T = [t_0,t_1] \in
\mathbb{R}$, $t_0 < t_1$, and using the scaling properties of Brownian motion, 
\begin{equation*}
  \begin{split}
  \mathbb{E}\left[
  U_{\gamma+2/p,p}(\mathbb{X}^2_{\text{It\^o}},T)^p\right] & =
  \mathbb{E} \int_{[t_0,t_1]^2}
  \frac{|\mathbb{X}^2_{\text{It\^o},uv}|^p}{|u-v|^{p\gamma+2}} du dv 
\\ & =  \mathbb{E}|\mathbb{X}^2_{\text{It\^o},0\,1}|^p \int_{[t_0,t_1]^2}
  |u-v|^{p(1-\gamma-2/p)} du dv 
< \infty    
  \end{split}
\end{equation*}
for any $\gamma < 1$ and $p > 1/(1-\gamma)$ so that,
a.s. $U_{\gamma+2/p,p}(\mathbb{X}^2_{\text{It\^o}},T)$ is finite
for any $\gamma < 1$ and $p$ sufficiently large.
Since
$$
\sup_{(u,v,w): s \le u \le v \le w\le t} |(\NN \mathbb{X}^2_{\text{It\^o}})_{uvw}| \le
\sup_{(u,v,w): s \le u \le v \le w\le t} |\der X_{uv}||\der X_{vw}| \le \|X\|_{\gamma,T}^2 |t-s|^{2\gamma}
$$
for any $t_0\le s \le t \le t_1$,
we have from~(\ref{eq:generalboundxx}) that  for any $\gamma < 1/2$, a.s.
$$
\|\mathbb{X}^2_{\text{It\^o},st}(\omega)\| \le C_{\gamma,T}(\omega) |t-s|^{2\gamma}
$$
for any $t,s \in I$, where $C_{\gamma,T}$ is a suitable random
constant.
Then  for any $\gamma < 1/2$ and bounded interval $I \subset \mathbb{R}$  we can choose a version such that $\mathbb{X}^2_{\text{It\^o}} \in
\OCC^{2\gamma}(I,V^{\otimes 2})$.

\bigskip
We can  introduce
$$
\mathbb{X}^{2,\mu\nu}_{\text{Strat.},st} \coloneq \int_s^t (X^\mu_u-X^\mu_s)
\circ \hat dX^\nu_u
$$ 
where the integral is understood in Stratonovich sense, then by well
known results in stochastic integration, we have
$$
\mathbb{X}^{2,\mu\nu}_{\text{Strat.},st} = \mathbb{X}^{2,\mu\nu}_{\text{It\^o},st} +\frac{g^{\mu\nu}}{2}
 (t-s) 
$$
where $g^{\mu\nu}=1$ if $\mu=\nu$ and $g^{\mu\nu} =0 $ otherwise. It
is clear that, also in this case, we can select a continuous version
of $\mathbb{X}^2_{\text{Strat.},st}$ which lives in $\OCC^{2\gamma}$
and such that $\NN\mathbb{X}^2_{\text{Strat.}} = \der X \der X$.

The connection between  stochastic integrals and the integral
we defined in Sec.~\ref{sec:irregular} starting from a couple
$(X,\mathbb{X}^2)$ is clarified in the next corollary:

\begin{corollary}
Let $\varphi \in C^{1,\delta}(V,V\otimes V^*)$ with $(1+\delta)\gamma > 1$, then the It\^o stochastic integral
$$
\der I_{\text{It\^o},st}^\mu = \int_s^t \varphi(X_u)^\mu_\nu \hat dX^\nu_u
$$
has a continuous version which is a.s. equal to
$$
\der I_{\text{rough},st}^\mu = \int_s^t \varphi(X_u)^\mu_\nu dX^\nu_u
$$
where the integral is understood in the sense of
Theorem~\ref{th:rough} based on the rough path
$(X,\mathbb{X}^2_{\text{It\^o}})$ moreover the Stratonovich integral
$$
\der I_{\text{Strat.},st}^\mu = \int_s^t \varphi(X_u)^\mu_\nu \circ
\hat dX^\nu_u
$$
is a.s. equal to the integral
$$
\der J_{st}^\mu = \int_s^t \varphi(X_u)^\mu_\nu dX^\nu_u
$$
defined based on the couple $(X,\mathbb{X}^2_{\text{Strat.}})$ and
the following relation holds
$$
\der J_{st}^\mu = \der I_{\text{rough},st}^\mu + \frac{g^{\nu\kappa}}{2}  \int_s^t \partial_\kappa \varphi(X_u)^\mu_\nu du
$$   
\end{corollary}
\proof
Recall that the It\^o integral $\der I_{\text{It\^o}}$ is the
limit in probability of the discrete sums
$$
S^\mu_\Pi = \sum_i \varphi(X_{t_i})^\mu_\nu (X^\nu_{t_{i+1}}-X^\nu_{t_i})
$$
while the integral $\der I_{\text{rough}}$ is the classical limit
as $|\Pi| \to 0$ of
$$
S^{\prime\,\mu}_\Pi = \sum_i \left[\varphi(X_{t_i})^\mu_\nu (X^\nu_{t_{i+1}}-X^\nu_{t_i})
+ \partial_\kappa \varphi(X_{t_i})^\mu_\nu \mathbb{X}^{2,\kappa\nu}_{\text{It\^o},t_i t_{i+1}}\right]
$$
(cfr. Corollary~\ref{cor:sums_rough}). Then it will suffice to show
that the limit in probability of
$$
R^\mu_\Pi = \sum_i \partial_\kappa \varphi(X_{t_i})^\mu_\nu \mathbb{X}^{2,\kappa\nu}_{\text{It\^o},t_i t_{i+1}}
$$
is zero. Since we assume $\partial \varphi$ bounded it will be enough to
show that $R_\Pi \to 0$ in $L^2(\Omega)$. By a standard argument, using the fact that $R_\Pi$
is a discrete martingale, we have
\begin{equation*}
  \begin{split}
\expect|R_\Pi|^2 & = \sum_i \expect |\partial_\kappa \varphi(X_{t_i})_\nu
\mathbb{X}^{2,\kappa\nu}_{\text{It\^o},t_i t_{i+1}}|^2      
 \le \|\varphi\|_{1,\delta} \sum_i \expect
|\mathbb{X}^{2}_{\text{It\^o},t_i t_{i+1}}|^2
\\& = \|\varphi\|_{1,\delta} \expect
|\mathbb{X}^{2}_{\text{It\^o},0 1}|^2 \sum_i |t_{i+1}-t_i|^2 
 \le \|\varphi\|_{1,\delta} \expect
|\mathbb{X}^{2}_{\text{It\^o},0 1}|^2 |\Pi| |t-s|
  \end{split}
\end{equation*}
which implies that $\expect|R_\Pi|^2 \to 0$ as $|\Pi| \to 0$.

As far as the integral $\der J$ is concerned, we have that it is the
classical limit of
\begin{equation*}
  \begin{split}
S^{\prime\prime\,\mu}_\Pi & = \sum_i \left[\varphi(X_{t_i})^\mu_\nu (X^\nu_{t_{i+1}}-X^\nu_{t_i})
+ \partial_\kappa \varphi(X_{t_i})^\mu_\nu \mathbb{X}^{2,\kappa\nu}_{\text{Strat.},t_i
  t_{i+1}}\right]
\\ & = 
\sum_i \left[\varphi(X_{t_i})^\mu_\nu (X^\nu_{t_{i+1}}-X^\nu_{t_i})
+ \partial_\kappa \varphi(X_{t_i})^\mu_\nu \mathbb{X}^{2,\kappa\nu}_{\text{It\^o},t_i
  t_{i+1}} + \frac{g^{\kappa\nu}}{2} \partial_\kappa  \varphi(X_{t_i})^\mu_\nu (t_{i+1}-t_i)
\right] 
\\ & = S^{\prime\,\mu}_\Pi + 
 \frac{g^{\kappa\nu}}{2}\sum_i  \partial_\kappa  \varphi(X_{t_i})^\mu_\nu (t_{i+1}-t_i)
  \end{split}
\end{equation*}
so that $$
\der I^\mu_{\text{rough},st} = \der J^\mu_{st} - \frac{g^{\kappa\nu}}{2}\int_s^t \partial_\kappa  \varphi(X_u)^\mu_\nu du
 $$
as claimed
and then, by the relationship between It\^o  and Stratonovich
integration:
$$
\der I_{\text{It\^o},st}^\mu = \der I_{\text{Strat.},st}^\mu - \frac{g^{\kappa\nu}}{2}\int_s^t \partial_\kappa  \varphi(X_u)^\mu_\nu du
$$
 we get $\der J = \der I_{\text{Strat.}}$. \qed

\section{Relationship with Lyons' theory of rough paths}
\label{sec:lyons}

The general abstract result given in Prop.~\ref{prop:main} can also be used to provide alternative proofs
of the main results in Lyons' theory of rough paths~\cite{Lyons},
i.e. the extension of multiplicative paths to any degree and the
construction of a multiplicative path from an almost-multiplicative
one. The main restriction is that we only consider control functions
$\omega(t,s)$ (cfr. Lyons~\cite{Lyons} for details and definitions)
which are given by
$$
\omega(t,s) = K |t-s|
$$
for some constant $K$. 

Given an integer $n$, $T^{(n)}(V)$ denote the truncated tensor algebra
up to degree $n$:   $T^{(n)}(V) := \oplus_{k=0}^n V^{\otimes k}$,
$V^{\otimes 0 } = \RR$. A tensor-valued path $Z : I^2
\to T^{(n)}(V)$ is of \emph{finite $p$-variation} if 
\begin{equation}
  \label{eq:finite_p_variation_path}
 \|Z^{\bar{\mu}}\|_{|\bar{\mu}|/p} \le K^{|\bar\mu|}, \qquad \forall
 \bar{\mu} : | \bar{\mu}| \le n  
\end{equation}
where $\bar{\mu}$ is a tensor multi-index. A path $Z$ of degree $n$
and finite $p$-variation is \emph{almost multiplicative} (of roughness $p$) if
$Z^\emptyset \equiv 1$, $n \ge \lfloor p \rfloor$ and
\begin{equation}
  \label{eq:almost_multiplicative}
  \NN Z^{\bar\mu} =  
\sum_{\bar\nu\bar\kappa =  \bar\mu}
  Z^{\bar\nu} Z^{\bar\kappa} + R^{\bar\mu}
\end{equation}
with $R^{\bar\mu} \in \OCC_2^z(I,T^{(n)}(V))$ for some $z>1$ uniformly
for all $\bar\mu$. By convention  the summation
$\sum_{\bar\nu\bar\kappa =  \bar\mu}$ does not include the terms where
either $\bar\mu =\emptyset$ or $\bar\kappa = \emptyset$. 

A path $Z$ is \emph{multiplicative} if
$Z^\emptyset \equiv 1$ and
\begin{equation}
  \label{eq:multiplicative}
  \NN Z^{\bar\mu} =  
\sum_{\bar\nu\bar\kappa =  \bar\mu}  Z^{\bar\nu} Z^{\bar\kappa} 
\end{equation}

Then the key result is contained in the following Proposition:
\begin{proposition}
 If $Z$ is an almost-multiplicative path of degree $n$ and finite
 $p$-variation,  $n \ge \lfloor p \rfloor$, then there exists a unique
 multiplicative path $\widetilde Z$ in $T^{(\lfloor p \rfloor)}(V)$ with finite
 $p$-variation such that
 \begin{equation}
   \label{eq:from_almost_to_multiplicative}
 \|Z^{\bar\mu}-\widetilde Z^{\bar\mu}\|_z \le K  
 \end{equation}
for some $z > 1$ and all multi-index $\bar\mu$ such that $|\bar\mu| \le \lfloor p \rfloor$.
\end{proposition}
\proof
Let us prove that there exists a multiplicative path $\widetilde Z$
such that
\begin{equation}
\label{eq:almost_remainder}
  Z = \widetilde Z + Q
\end{equation}
with $Q \in \OCC^z$, $z > 1$. We proceed by induction: if
$|\bar\mu|=1$:
\begin{equation*}
  \NN Z^{\bar\mu}_{sut} = R^{\bar\mu}_{sut}
\end{equation*}
which, given that $R^{\bar\mu} \in \OCC_2^z$, $z >1$ implies that
exists a unique $\widetilde Z^{\bar\mu}$ such that $\NN \widetilde Z^{\bar\mu}
= 0$ and
\begin{equation*}
  Z^{\bar\mu} =  \widetilde Z^{\bar\mu} + \Lambda R^{\bar\mu} = \widetilde Z^{\bar\mu} + Q^{\bar\mu} 
\end{equation*}
with $Q^{\bar\mu} \in \OCC^z$. Then assume that
eq.~(\ref{eq:almost_remainder}) is true up to degree $j-1$ and let us
show that it is true also for a multi-index $\bar\mu$ of degree $j$:
\begin{equation*}
\begin{split}
\NN Z^{\bar\mu} & = \sum_{\bar\nu \bar\kappa = \bar\mu}
Z^{\bar\nu} Z^{\bar\kappa} + R^{\bar\mu}
\\
& = \sum_{\bar\nu \bar\kappa = \bar\mu}
(\widetilde Z^{\bar\nu}+Q^{\bar\nu})(\widetilde
Z^{\bar\kappa}+Q^{\bar\kappa}) + R^{\bar\mu}
\\ & = \sum_{\bar\nu \bar\kappa = \bar\mu}
\widetilde Z^{\bar\nu} \widetilde
Z^{\bar\kappa} + \sum_{\bar\nu \bar\kappa = \bar\mu} [
Q^{\bar\nu} \widetilde
Z^{\bar\kappa} + \widetilde Z^{\bar\nu} Q^{\bar\kappa}
+ Q^{\bar\nu} Q^{\bar\kappa}]
 + R^{\bar\mu}
\\ 
& = \sum_{\bar\nu \bar\kappa = \bar\mu}
\widetilde Z^{\bar\nu} \widetilde
Z^{\bar\kappa} + \widetilde R^{\bar\mu}
    \end{split}
\end{equation*}
If we can prove that $\widetilde R^{\bar\mu}$ is in the image of
$\NN$, then writing 
$$
\widetilde Z^{\bar\mu}  = Z^{\bar\mu} - \Lambda \widetilde R^{\bar\mu}
= Z^{\bar\mu} + Q^{\bar\mu}
$$
we obtain the multiplicative property for $\widetilde Z^{\bar\mu}$
$$
\NN \widetilde Z^{\bar\mu}  = \sum_{\bar\nu \bar\kappa = \bar\mu}
\widetilde Z^{\bar\nu}_{ut} \widetilde
Z^{\bar\kappa}_{su}
$$
with $|\bar\mu| = j$, and we are done since uniqueness is obvious.
To prove $\widetilde R^{\bar\mu} \in \text{Im}\NN$ we must show that 
$\NN_2 \widetilde R^{\bar\mu} = 0$: 
\begin{equation*}
  \begin{split}
\NN_2 \widetilde R^{\bar\mu} & = 
\NN_2 \left[ \NN Z^{\bar \mu}-\sum_{\bar\nu \bar\kappa = \bar\mu}
\widetilde Z^{\bar\nu} \widetilde
Z^{\bar\kappa}\right]
 = 
N_2    \left[\sum_{\bar\nu \bar\kappa = \bar\mu}
\widetilde Z^{\bar\nu} \widetilde
Z^{\bar\kappa}\right]
\\ & = \sum_{\bar\nu \bar\kappa =
\bar\mu} N  \widetilde Z^{\bar\nu} \widetilde
Z^{\bar\kappa} - \sum_{\bar\nu \bar\kappa =
\bar\mu}  \widetilde Z^{\bar\nu} N \widetilde
Z^{\bar\kappa}
\\ & =     
\sum_{\bar\nu \bar\kappa =
\bar\mu} \sum_{\bar\sigma\bar\tau = \bar\nu} \widetilde Z^{\bar\sigma} \widetilde Z^{\bar\tau}\widetilde
Z^{\bar\kappa} - \sum_{\bar\nu \bar\kappa =
\bar\mu} \sum_{\bar\sigma\bar\tau = \bar\kappa} \widetilde
Z^{\bar\nu} \widetilde Z^{\bar\sigma} \widetilde
Z^{\bar\tau}   = 0
  \end{split}
\end{equation*}
where we used the  Leibnitz rule for $N_2$
(see eq.~(\ref{eq:leibnitz_n2})).

To finish we can take for the constant $K$ in
 eq.~(\ref{eq:from_almost_to_multiplicative}) the maximum of
 $\|Q^{\bar\mu}\|_z$ for all $|\bar\mu| \le \lfloor p \rfloor$. 
 \qed

\begin{proposition}
Let $Z$ be a multiplicative path of degree $n$ and finite $p$-variation such that
\begin{equation}
  \label{eq:multiplicative_bound}
 \sum_{\bar\mu : |\bar\mu|=k}\norm{Z^{\bar\mu}}_{k/p} \le C \frac{\alpha^{k}}{k!} 
\end{equation}
for all $k \le n$ and with $\alpha,C >0$; then if $(n+1) > p$ and $C$
is small enough (see eq.~(\ref{eq:smallness_of_C})) there exists a unique
multiplicative extension of $Z$ to any degree and
eq.~(\ref{eq:multiplicative_bound}) holds for every $k$.
\end{proposition}
\proof
By induction we can assume that $Z$ is a multiplicative path of degree
$k$ for which eq.~(\ref{eq:multiplicative_bound}) holds up to degree
$k$ and prove that it can be extended to degree $k+1$ with the same
bound. Note that $k \ge n$ and then $(k+1) > p$.
For $|\bar\mu|=k+1$ we should have
\begin{equation}
\label{eq:extension_decomposition}
  \NN Z^{\bar\mu} = \sum_{\bar\nu\bar \kappa = \bar\mu} Z^{\bar\nu}
  Z^{\bar\kappa} \in \ZZ_2^{(k+1)/p}
\end{equation}
Since $(k+1) > p$, this equation has a unique solution $Z^{\bar\mu}
\in \OCC^{(k+1)/p}(T^{k+1}(V))$.
Then observe
that, from eq.~(\ref{eq:extension_decomposition})
$$
 Z^{\bar\mu}_{st} = Z^{\bar\mu}_{ut} + Z^{\bar\mu}_{su} + \sum_{\bar\nu\bar \kappa = \bar\mu} Z^{\bar\nu}_{su}
  Z^{\bar\kappa}_{ut} 
$$
and taking as $u$ the mid-point between $t$ and $s$ we can bound
$Z^{\bar \mu}$ as
follows:
\begin{equation*}
  \sum_{|\bar\mu| = k+1} \norm{Z^{\bar\mu}_{st}}_{(k+1)/p} \le
  \frac{2}{2^{(k+1)/p}} \sum_{|\bar\mu| = k+1}
  \norm{Z^{\bar\mu}_{st}}_{(k+1)/p} +C^2 \alpha^{k+1}
\sum_{i=1}^k 
\frac{2^{-i/p}}{i!} \frac{2^{-(k+1-i)/p}}{(k+1-i)!} 
\end{equation*}
Now,
\begin{equation*}
  \begin{split}
\sum_{i=0}^{k+1} &
\frac{2^{-i/p}}{i!} \frac{2^{-(k+1-i)/p}}{(k+1-i)!} \le
    \sum_{i=0}^{k+1} 
\frac{2^{-i}}{i!} \frac{2^{-(k+1-i)}}{(k+1-i)!} +2 \sum_{i=0}^{\lfloor p \rfloor}
\frac{ (2^{-(k+1-i)/p} 2^{-i/p}-2^{-(k+1-i)} 2^{-i})}{i!(k+1-i)!}
\\ & = \frac{1}{(k+1)!}\left[1+2
\sum_{i=0}^{\lfloor p \rfloor}
\frac{(k+1)!}{i!(k+1-i)!} (2^{-(k+1)/p} -2^{-(k+1)})
\right]
\\ & \le \frac{1+D_p k^{\lfloor p \rfloor} 2^{-(k+1)/p}}{(k+1)!}
  \end{split}
\end{equation*}
which gives
\begin{equation*}
  \sum_{|\bar\mu| = k+1} \norm{Z^{\bar\mu}_{st}}_{(k+1)/p} \le
  C^2 \frac{(2^{(k+1)/p}-2)}{2^{(k+1)/p}}  
\frac{(1+D_p k^{\lfloor p \rfloor} 2^{-(k+1)/p})\alpha^{k+1}}{(k+1)!}
\le
C \frac{\alpha^{k+1}}{(k+1)!}
\end{equation*}
whenever $C$ is such that
\begin{equation}
\label{eq:smallness_of_C}
 0 < C \le  \min_{k \ge n} \frac{2^{(k+1)/p}}{(2^{(k+1)/p}-2)(1+D_p k^{\lfloor p \rfloor}
 2^{-(k+1)/p})}.
\end{equation}
This concludes the proof of the induction step.
\qed

\section*{Acknowledgments}
The author has been introduced to this problem by F.~Flandoli which
has sustained his work with useful discussions and  constant
encouragement. 
A due thank to M.~Franciosi for some advices on
homological algebra and to Y.~Ouknine for pointing out a mistake in an
earlier version of the paper.


\appendix
\section{Some proofs}
\label{app:proofs}
\subsection{Proof of prop.~\ref{prop:main}}
\label{sec:proof_main}
The basic technique to prove the existence of the map $\Lambda$ is
borrowed form~\cite{FGGT}.
Let $\eta(x)$ be a smooth function on $\mathbb{R}$ with compact support and
$\eta_{\alpha}(x) \coloneq \alpha^{-1} \eta(x/\alpha)$.

Define
\begin{equation*}
  (\Lambda_\beta A)_{st} \coloneq  - \int_s^t dx \iint d\tau d\sigma
  \mathcal{F}_\beta(x,s;\tau,\sigma)  A_{\tau x \sigma} 
\end{equation*}
where
\begin{equation*}
  \mathcal{F}_\beta(x,s;\tau,\sigma) \coloneq
[ \eta_{\beta}(x-\tau)- \eta_{\beta}(s-\tau)] \partial_x \eta_{\beta}(x-\sigma)
\end{equation*}
and the integrals in $\tau$ and $\sigma$ are extended over all $\RR$.

Given that $A \in \ZZ_2$ there exists  $R \in \OCC$ such that
$\NN R = A$ and  
\begin{equation*}
  \begin{split}
  (\Lambda_\beta A)_{st} & = - \int_s^t dx \iint d\tau d\sigma
  \mathcal{F}_\beta(x,s;\tau,\sigma)  (R_{\tau  \sigma}-R_{\tau x
  }-R_{ x \sigma})
\\ & =  - \int_s^t dx \iint d\tau d\sigma
  \mathcal{F}_\beta(x,s;\tau,\sigma)  R_{\tau  \sigma}
  \end{split}
\end{equation*}
since the other terms vanish after the integrations in $\tau$ or
$\sigma$.
Then the following decomposition holds:
\begin{equation}
\label{eq:decomposition}
  \Lambda_\beta A = \tilde R_\beta + \der \Phi_\beta(R)
\end{equation}
where
\begin{equation*}
  (\tilde R_\beta)_{st} \coloneq  \iint d\tau d\sigma \eta_\beta(s-\tau)
  [\eta_\beta(t-\sigma)-\eta_\beta(s-\sigma)] R_{\tau \sigma}
\end{equation*}
and
\begin{equation*}
  \der \Phi_\beta(R)_{st} \coloneq  -
  \int_s^t dx \iint d\sigma d\tau \eta_\beta(x-\tau) \partial_x
  \eta_\beta(x-\sigma) R_{\tau \sigma} 
\end{equation*}
 
In eq.~(\ref{eq:decomposition}) the
l.h.s. depends only on $A = \NN R$ while each of the terms in the r.h.s
depends explicitly on $R$. We have $\NN \Lambda_\beta A = \NN
\tilde R_\beta$ and since $\lim_{\beta \to 0} \tilde R_\beta = R$
pointwise we have that $\lim_{\beta \to 0} \NN \Lambda_\beta A = \NN
R = A$. So every accumulation point $X$ of $\Lambda_\beta A$ will solve
the equation $\NN X = A$. Moreover if it
exists $X \in \OCC^{z}$ with $z > 1$ and $\NN X = A$ then it is 
unique and 
$\lim_{\beta \to 0} \Lambda_\beta R = X$ in $\OCC^{1}$ since in this case
\begin{equation*}
  \Lambda_\beta A = \tilde R_\beta + \der \Phi_\beta(R) = \tilde X_\beta + \der \Phi_\beta(X)   
\end{equation*}
and it is easy to prove that $\Phi_\beta(X) \to 0$ in $\CC^1$.    

Now we will prove that $\lim_{\beta\to 0 }\Lambda_\beta A$ exists when
$A \in \ZZ_2^z$ with $z > 1$.

Define $f_\tau : \mathbb{R}^2 \times \mathbb{R}_+ \to V$ as $f_\tau (x,y,\alpha)
:= \eta_{\alpha}(x-\tau)$ and $g_\sigma : \mathbb{R}^2 \times \mathbb{R}_+ \to V$ as $g_\sigma (x,y,\alpha) := \eta_{\alpha}(y-\sigma)$.    
Apply Stokes Theorem to the exact
differential 2-form $\omega \coloneq df_\tau \wedge d g_\sigma =
d(f_\tau d g_\sigma)$  on $D \coloneq \Delta_{t,s} \times
[\beta,\beta']$ where $\Delta_{t,s} = \{(x,y) \in \mathbb{R}^2 : s < x < y <
t\}$. Then
\begin{equation*}
  \int_{\partial D} \omega =  \int_D d \omega = 0  
\end{equation*}
where the boundary $\partial D = - c_1 + c_2 + c_3 $ is composed of
$c_1 = \Delta_{t,s} \times \{ \beta \}$, $c_2 = \Delta_{t,s} \times \{
\beta' \}$, $c_3 = \partial \Delta_{t,s} \times [\beta,\beta']$.
So
\begin{equation*}
  \int_{\Delta_{t,s}} \omega |_{\alpha = \beta} =  \int_{\Delta_{t,s}}
  \omega |_{\alpha = \beta'} + \int_{ \partial \Delta_{t,s} \times [\beta,\beta']}
  \omega 
\end{equation*}
giving
\begin{equation*}
  \int_s^t \mathcal{F}_\beta(x,s;\tau,\sigma) dx
  =  \int_s^t \mathcal{F}_{\beta'}(x,s;\tau,\sigma) dx
+ \int_\beta^{\beta'} d\alpha \int_s^t
 \mathcal{K}(\alpha,x,t,s;\tau,\sigma) dx 
\end{equation*}
with
\begin{multline*}
  \mathcal{K}
  (\alpha,x,t,s;\tau,\sigma) = \partial_\alpha
  [\eta_\alpha(x-\sigma)-\eta_\alpha(s-\sigma)] \partial_x \eta_\alpha(x-\tau)
\\+\partial_\alpha
  [\eta_\alpha(t-\tau)-\eta_\alpha(x-\tau)] \partial_x \eta_\alpha(x-\sigma)
\end{multline*}

Then 
\begin{equation}
\label{eq:telescope}
  \Lambda_\beta A_{st} = \Lambda_{\beta'} A_{st}
- \int_\beta^{\beta'} d\alpha \int_s^t dx
 \iint d\tau d\sigma \,\mathcal{K}(\alpha,x,t,s;\tau,\sigma) R_{\tau \sigma}   
\end{equation}
Assume we can write $A = \sum_{i=1}^n A_i$ where $A_i \in
\OCC_2^{\rho_i,z-\rho_i}$ for a choice of $n$ and $\rho_i > 0$,
$i=1,\dots,n$. Write $\rho_i' = z-\rho_i$.

Then consider
\begin{equation*}
  \begin{split}
 I&(\alpha)  =
 - \int_s^t dx \iint d\tau d\sigma
 \mathcal{K}(\alpha,x,t,s;\tau,\sigma) R_{\tau \sigma} 
\\ & =
  \iint d\tau d\sigma \left\{\partial_\alpha \eta_\alpha(s-\sigma)
 [\eta_\alpha(t-\tau)-\eta_\alpha(s-\tau)] \right.\\ & \qquad \left. -
 \partial_\alpha \eta_\alpha(t-\tau)
 [\eta_\alpha(t-\sigma)-\eta_\alpha(s-\sigma)]\right\} R_{\tau \sigma} 
\\ & 
 + \int_s^t dx \iint d\tau d\sigma
[\partial_\alpha \eta_\alpha(x-\tau) \partial_x \eta(x-\sigma) -
 \partial_\alpha \eta_\alpha(x-\sigma) \partial_x \eta(x-\tau)]
 R_{\tau \sigma}
\\ & = 
  \iint d\tau d\sigma \partial_\alpha \eta_\alpha(\sigma)
 \eta_\alpha(\tau) [R_{t+\tau,s+\sigma }-R_{s+\tau,s+\sigma }-R_{t+\tau,t+\sigma}+R_{t+\tau,s+\sigma} ]
\\ & 
 + \int_s^t dx \iint d\tau d\sigma
\partial_\alpha \eta_\alpha(\tau) \partial_\sigma \eta(\sigma) 
 [R_{x+\sigma,x+\tau}-R_{x+\tau,x+\sigma}]
\\ & = 
  \iint d\tau d\sigma \partial_\alpha \eta_\alpha(\sigma)
 \eta_\alpha(\tau) [\NN R_{t+\tau,s+\tau,s+\sigma }+\NN R_{t+\tau,t+\sigma,s+\sigma} ]
\\ & 
 + \int_s^t dx \iint d\tau d\sigma
\partial_\alpha \eta_\alpha(\tau) \partial_\sigma \eta(\sigma) 
 [\NN R_{x+\sigma,x,x+\tau}-\NN R_{x+\tau,x,x+\sigma}]
  \end{split}
\end{equation*}
so that we can bound
\begin{equation*}
  \begin{split}
  |I(\alpha)| & \le     
  \iint d\tau d\sigma |\partial_\alpha \eta_\alpha(\sigma)|
 |\eta_\alpha(\tau)| \left[|\NN R_{t+\tau,s+\tau,s+\sigma }|+|\NN R_{t+\tau,t+\sigma,s+\sigma}| \right]
\\ & 
 \qquad + \int_s^t dx \iint d\tau d\sigma
|\partial_\alpha \eta_\alpha(\tau)|| \partial_\sigma \eta(\sigma)| 
 \left[|\NN R_{x+\sigma,x,x+\tau}|+|\NN R_{x+\tau,x,x+\sigma}|\right]
\\ & \le     \sum_{i=1}^n \|A_i\|_{\rho_i,\rho_i'}
  \iint d\tau d\sigma |\partial_\alpha \eta_\alpha(\sigma)|
 |\eta_\alpha(\tau)| \left[|t-s|^{\rho_i} |\tau-\sigma|^{\rho_i'}
+|\tau-\sigma|^{\rho_i} |t-s|^{\rho_i'}  \right]
\\ & 
 \qquad +    \sum_{i=1}^n \|A_i\|_{\rho_i,\rho_i'} \int_s^t dx \iint d\tau d\sigma
|\partial_\alpha \eta_\alpha(\tau)|| \partial_\sigma \eta(\sigma)| 
 \left[|\sigma|^{\rho_i} |\tau|^{\rho_i'}+|\tau|^{\rho_i} |\sigma|^{\rho_i'}\right]
  \end{split}
\end{equation*}
where each term can be bounded as follows:
\begin{equation*}
  \iint d\tau d\sigma |\partial_\alpha \eta_\alpha(\sigma)|
 |\eta_\alpha(\tau)| |\tau-\sigma|^{a} = \alpha^{a-1} \iint d\tau d\sigma
 |\eta(\sigma)- \sigma \eta'(\sigma)|
 |\eta(\tau)|\, |\tau-\sigma|^{a}  \le K \alpha^{a-1},
\end{equation*}
\begin{equation*}
   \iint d\tau 
|\partial_\alpha \eta_\alpha(\tau)| |\tau|^{a} = \alpha^{a-1}
    \iint d\tau 
|\eta(\tau)-\tau\eta'(\tau)| |\tau|^{a} \le K^{1/2} \alpha^{a-1}
\end{equation*}
for a suitable constant $K>0$
and obtain
\begin{equation*}
  \begin{split}
  |I(\alpha)| & \le K \sum_{i=1}^n (\alpha^{\rho_i-1} |t-s|^{\rho_i'}
   +\alpha^{\rho_i'-1} |t-s|^{\rho_i}) \|A_i\|_{\rho_i,\rho_i'}
\\   & + K |t-s| \sum \alpha^{z-2} \|A_i\|_{\rho_i,\rho_i'}
  \end{split}
\end{equation*}
Upon integration in $\alpha$  we get:
\begin{equation*}
 \int_0^1 |I(\alpha)| d\alpha \le K \sum_{i=1}^n \|A_i\|_{\rho_i,\rho_i'}  
\end{equation*}
if $|t-s| \le 1$. 
By dominated convergence of the integral in eq.~(\ref{eq:telescope}), 
$$\lim_{\beta \to 0} \Lambda_\beta A \eqcolon \Lambda A$$
exists (in $\OCC$ uniformly in bounded intervals).
If we also observe that
\begin{equation*}
  |(\Lambda_{\beta'} A)_{st}| \le K (\beta^{\prime})^{-1}|t-s| \sum_{i=1}^n \|A_i\|_{\rho_i,\rho_i'}
\end{equation*}
we get that
\begin{equation*}
  |(\Lambda A)_{t,s}| \le K \sum_{i=1}^n \|A_i\|_{\rho_i,\rho_i'}
\end{equation*}
for $|t-s| \le 1$.

Finally, let $J_{t,s}(x) \coloneq s+(t-s)(0 \vee (x\wedge 1))$
and $(J_{t,s}^* X)_{u,v,w} \coloneq
X_{J_{t,s}(u),J_{t,s}(v),J_{t,s}(w)}$ for all $X \in \OCC_2$.  
Then
\begin{equation*}
  \|J_{t,s}^* X\|_{\gamma,\gamma'} \le |t-s|^{\gamma+\gamma'} \|X\|_{\gamma,\gamma'}.
\end{equation*}

Since $\Lambda_\beta A_{t,s} = (J_{t,s}^* \Lambda_{|t-s|\beta} A)_{0,1} = \Lambda_{|t-s|\beta}
(J_{t,s}^* A)_{0,1}$ and 
\begin{equation*}
  |(\Lambda (J^*_{t,s} R))_{1,0}| \le 
K \sum_{i=1}^n \|J^*_{t,s} A_i\|_{\rho_i,\rho_i'}
\end{equation*}
this is enough to obtain the desired regularity:
\begin{equation*}
  |(\Lambda A)_{t,s}| \le 
K  |t-s|^{z} \sum_{i=1}^n \|A_i\|_{\rho_i,\rho_i'}.
\end{equation*}

The constant $K$ can be chosen to be equal to $1/(2^z-2)$. 
Let
 $\Phi = \sum_{i=1}^n
\|A_i\|_{\rho_i,\rho_i'}$.
and $R=\Lambda A$ and
since $N R = A$ write
\begin{equation*}
  R_{st} = R_{ut} + R_{su} + \sum_i A_{i,sut}
\end{equation*}
with $t>u>s$ and $u = s+|t-s|/2$. Then estimate
\begin{equation*}
  \begin{split}
  |R_{st}| & \le |R_{ut}| + |R_{su}| + \sum_i |A_{i,sut}|
    \\ & \le  \norm{R}_z (|t-u|^z+|u-s|^z) + \sum_i
    \|A_i\|_{\rho_i,\rho_i'} |u-s|^{\rho_i} |t-u|^{\rho'_i}  
\\ & =  \frac{2 \norm{R}_z + \Phi}{2^z}  |t-s|^z 
  \end{split}
\end{equation*}
so that
\begin{equation*}
  \norm{R}_z \le \frac{1}{2^z-2} \Phi.
\end{equation*}
\qed

\subsection{Some Proofs for Sec.~\ref{sec:irregular}}
\subsubsection{Proof of Lemma~\ref{lemma:transitivity}}
\label{sec:proof-transitivity}
\proof
Write down the decomposition for $Z$ and $Y$:
\begin{equation*}
  \begin{gathered}
  \der Z^\mu = F^{\mu}_{\nu} \der Y^\nu + R^\mu_{ZY},\\
  \der Y^\mu = G^{\mu}_{\nu} \der X^\nu + R^\mu_Y\\
  \end{gathered}
\end{equation*}
where $F \in \CC^{\eta-\gamma}(I,V \otimes V^*)$, $G \in \CC^{\sigma-\gamma}(I,V)$, $R_{ZY} \in
\OCC^\eta(I,V)$ and $R_Y \in \OCC^\sigma(I,V)$, then
\begin{equation*}
  \der Z^\mu = F^{\mu}_{\nu} G^{\nu\kappa} \der X^\kappa + R^\mu_{ZY} + F^{\mu}_{\nu} R^\nu_Y = Z^{\prime\,\mu}_{\kappa} \der X^\kappa + R^\mu_{ZX}
\end{equation*}
with $Z^{\prime\,\mu}_{\kappa} = F^{\mu}_{\nu} G^{\nu}_{\kappa}$ and $R^\mu_{ZX} = R^\mu_{ZY} + F^{\mu}_{\nu} R^\nu_{Y}$.
Let $\delta = \min(\sigma,\eta)$ and
note that for $R_{ZY}$ we have
\begin{equation*}
  \begin{gathered}
\|R_{ZY}\|_{\eta,I} \le \|Z\|_{D(Y,\gamma,\eta),I}\\
\|R_{ZY}\|_{\gamma,I} \le \|Z\|_{\gamma,I} + \|F\|_{\infty,I} \|Y\|_{\gamma,I} \le
\|Z\|_{D(Y,\gamma,\eta),I} (1+\|Y\|_{\gamma,I})    
  \end{gathered}
\end{equation*}
and by interpolation we obtain ($a=(\eta-\delta)/(\eta-\gamma) \le 1$)
$$
\|R_{ZY}\|_{\delta,I} \le \|R_{ZY}\|_{\eta,I}^{1-a} \|R_{ZY}\|_{\gamma,I}^{a}  \le \|Z\|_{D(Y,\gamma,\eta),I} (1+\|Y\|_{\gamma,I})^a 
\le \|Z\|_{D(Y,\gamma,\eta),I} (1+\|Y\|_{\gamma,I}) 
$$
and similarly
$$
\|R_{Y}\|_{\delta,I} \le \|Y\|_{D(X,\gamma,\sigma),I} (1+\|X\|_{\gamma,I})
$$
moreover
$$
\|F\|_{0,I} = \sup_{t,s \in I}|F_t-F_s| \le \sup_{t,s \in I} (|F_t|+ |F_s|) = 2
\|F\|_{\infty,I}  \le 2 \|Z\|_{D(Y,\gamma,\eta),I} 
$$
so, again by interpolation, we find
$$
\|F\|_{\delta-\gamma,I} \le \|Z\|_{D(Y,\gamma,\eta),I}
2^{1-(\delta-\gamma)/(\sigma-\gamma)} \le 2  \|Z\|_{D(Y,\gamma,\eta),I}
$$
and
$$
\|G\|_{\delta-\gamma,I} \le 2  \|Y\|_{D(X,\gamma,\sigma),I}
$$

To finish bound the norm of
$Z,Z'$ as 
\begin{equation*}
  \begin{split}
\|(Z,Z')\|_{D(X,\gamma,\delta),I} & = \|Z'\|_{\infty,I} +
\|Z'\|_{\delta-\gamma,I} + \|R_{ZX}\|_{\delta,I} + \|Z\|_{\gamma,I}   
\\ & \le
 \|F\|_{\infty,I} \|G\|_{\infty,I} +
\|F\|_{\delta-\gamma,I}\|G\|_{\infty,I} 
\\ & \qquad + \|F\|_{\infty,I} \|G\|_{\delta-\gamma,I} 
+ \|R_{ZY}\|_{\delta,I} + \|F\|_{\infty,I} \|R_Y\|_{\delta,I} + \|Z\|_{\gamma,I}   
\\ & \le K \|Z\|_{D(Y,\gamma,\eta),I} (1+
   \|Y\|_{D(X,\gamma,\sigma),I})(1+ \|X\|_{\gamma,I})
  \end{split}
\end{equation*}
\qed

\subsubsection{Proof of Prop.~\ref{prop:functionD}}
\label{sec:proof_of_function_D}
Let
$y(r) = (Y_t-Y_s)r+Y_s$
so that
\begin{equation*}
  \begin{split}
Z^\mu_t-Z^\mu_s & = \varphi(y(1))^\mu-\varphi(y(0))^\mu = \int_0^1
 \partial_\nu \varphi(y(r))^\mu y'(r)^\nu dr 
\\ & =  (Y^\nu_t-Y^\nu_s) \int_0^1
\partial_\nu \varphi(y(r))^\mu dr 
\\ & =  \partial_\nu\varphi(Y_s)^\mu (Y^\nu_t-Y^\nu_s) + (Y^\nu_t-Y^\nu_s) \int_0^1 \left[\partial_\nu \varphi(y(r))^\mu-\partial_\nu \varphi(Y_s)^\mu\right] dr 
  \end{split}
\end{equation*}
Then if $\der Y^\mu =  Y^{\prime\,\mu}_{\nu}  \der X^\nu + R^\mu$ we have
\begin{equation}
\label{eq:decomposition_of_Z}
  \begin{split}
Z^\mu_t-Z^\mu_s &  =  \partial_\nu \varphi(Y_s)^\mu Y^{\prime\,\nu}_{\kappa,s} (X^\kappa_t-X^\kappa_s) 
+  \partial_\nu\varphi(Y_s)^\mu R^\nu_{st} 
+ (Y^\nu_t-Y^\nu_s) \int_0^1 \left[\partial_\nu \varphi(y(r))^\mu-\partial_\nu
  \varphi(Y_s)^\mu\right] dr 
\\ & =  Z^{\prime\,\mu}_{\kappa,s} (X^\kappa_t-X^\kappa_s) + R^{\mu}_{Z,st}
  \end{split}
\end{equation}
with $Z^{\prime\,\mu}_{\kappa,s} = \partial_\nu \varphi(Y_s)^\mu
Y^{\prime\,\nu}_{\kappa,s}$, 
\begin{equation*}
  \begin{split}
\|Z'\|_{\sigma-\gamma} &\le \|\partial \varphi(Y_\cdot)\|_{\sigma-\gamma} \|Y'\|_\infty+
\|\partial \varphi(Y_\cdot)\|_{\infty} \|Y'\|_{\sigma-\gamma}
\\ &\le (\|\partial \varphi(Y_\cdot)\|_{\delta \gamma}+\|\partial \varphi(Y_\cdot)\|_{0}) \|Y'\|_\infty+
\|\partial \varphi(Y_\cdot)\|_{\infty} (\|Y'\|_{\eta-\gamma}+\|Y'\|_{0})
\\ &\le \|\varphi\|_{1,\delta}(\|Y\|^\delta_{\gamma}+2) \|Y'\|_\infty+
2 \|\varphi\|_{1,\delta} (\|Y'\|_{\eta-\gamma}+2\|Y'\|_{\infty})
\\ & \le K \|\varphi\|_{1,\delta}
(\norm{Y}_{D(X,\gamma,\eta)}+\norm{Y}_{D(X,\gamma,\eta)}^{1+\delta})
  \end{split}
\end{equation*}
as far as $R_Z$ is concerned we have
\begin{equation*}
  \begin{split}
 |R_{Z,st}| & = |Y_t-Y_s|    
   \left[\int_0^1 \left|\partial \varphi(y(r))-\partial
  \varphi(Y_s)\right| dr\right] 
\\ & \le \|\varphi\|_{1,\delta} \left|\int_0^1
  r^\delta  dr\right| |Y_t-Y_s|^{1+\delta} \le K \|\varphi\|_{1,\delta}
  \|Y\|_{\gamma}^{1+\delta} |t-s|^{\gamma (1+\delta)};
  \end{split}
\end{equation*}
and
\begin{equation*}
  \begin{split}
 |R_{Z,st}| & = |Y_t-Y_s|    
   \left[\int_0^1 \left|\partial \varphi(y(r))-\partial
  \varphi(Y_s)\right| dr\right] \le K \|\varphi\|_{1,\delta}
   \|Y\|_\gamma |t-s|^\gamma.
\end{split}
\end{equation*}
Interpolating these two inequalities we get
\begin{equation*}
  \norm{R_Z}_{\sigma} \le K \|\varphi\|_{1,\delta}
  \|Y\|^{\sigma/\gamma}_\gamma
\le K \|\varphi\|_{1,\delta} \|Y\|^{\sigma/\gamma}_{D(X,\gamma,\sigma)}
\end{equation*}
which together with the obvious bound
\begin{equation*}
  \norm{Z}_\gamma \le \|\varphi\|_{1,\delta}\norm{Y}_\gamma
\end{equation*}
implies
 \begin{equation*}
   \begin{split}
\|Z\|_{D(X,\gamma,\sigma)} & \le 
K \|\varphi\|_{1,\delta}(
\norm{Y}_{D(X,\gamma,\eta)} + \norm{Y}^{1+\delta}_{D(X,\gamma,\eta)} + \norm{Y}^{\sigma/\gamma}_{D(X,\gamma,\eta)})
   \end{split}
 \end{equation*}

If $\der\widetilde Y^\mu =  \widetilde Y^{\prime\,\mu}_{\nu}  \der X^\nu
+ \widetilde  R^\mu$ is another path, $\widetilde Z_t =
\varphi(\widetilde Y_t)$ and $H = Z -\widetilde Z$ we have (see
eq.~(\ref{eq:decomposition_of_Z})):
\begin{equation}
\label{eq:difference_H}
  \begin{split}
\der H^\mu &  = H^{\prime\,\mu}_{\kappa} \der X^\kappa + A^\mu + B^\mu
  \end{split}
\end{equation}
with 
$$
H^{\prime\,\mu}_{\kappa} = \partial_\nu \varphi(Y)^\mu
Y^{\prime\,\nu}_{\kappa}-\partial_\nu \varphi(\widetilde Y)^\mu \widetilde Y^{\prime\,\nu}_{\kappa}
$$
$$
A^{\mu}_{st} = \partial_\nu\varphi(Y_s)^\mu R^\nu_{st} -
\partial_\nu\varphi(\widetilde  Y_s)^\mu \widetilde  R^\nu_{st}
$$
and
\begin{equation*}
  \begin{split}
B^{\mu}_{st} & =
 \der Y^\nu_{st} \int_0^1 \left[\partial_\nu
  \varphi(y(r))^\mu-\partial^\nu \varphi(y(0))^\mu\right] dr -
\der \widetilde Y^\nu_{st} \int_0^1 \left[\partial_\nu
  \varphi(\widetilde y(r))^\mu-\partial_\nu \varphi(\widetilde
  y(0))^\mu\right] dr 
\\ & = \der (Y-\widetilde Y)^\nu_{st} \int_0^1 \left[\partial_\nu
  \varphi(y(r))^\mu-\partial_\nu \varphi(y(0))^\mu\right] dr 
\\ & \qquad +
\der \widetilde Y^\nu_{st} \int_0^1 \left[
\partial_\nu\varphi(y(r))^\mu-\partial_\nu
  \varphi(\widetilde y(r))^\mu-\partial_\nu \varphi(y(0))^\mu+\partial_\nu \varphi(\widetilde
  y(0))^\mu\right] dr 
  \end{split}
\end{equation*}
Let $y(r,r') = (y(r) -  \widetilde  y(r)) r' + \widetilde y(r)$ and bound the
second integral as
\begin{equation*}
  \begin{split}
\Big|\int_0^1 dr & \left[
\partial_\nu\varphi(y(r))^\mu-\partial_\nu
  \varphi(\widetilde y(r))^\mu-\partial_\nu
  \varphi(y(0))^\mu+\partial_\nu \varphi(\widetilde y(0))^\mu\right] \Big|
\\ & =
\Abs{\int_0^1 dr \int_0^1 dr'  \left[
\partial_{\kappa} \partial_\nu\varphi(y(r,r'))^\mu-\partial_{\kappa} \partial_\nu \varphi(y(0,r'))^\mu\right] (y(r) -  \widetilde  y(r))^{\kappa} }
\\ & \le \|\varphi\|_{2,\delta}
{\int_0^1 dr \int_0^1 dr'  |y(r,r')-y(0,r')|^\delta } \abs{y(r) -  \widetilde  y(r) }
\\ & \le K \|\varphi\|_{2,\delta} (\norm{Y}_{\gamma}+\norm{\widetilde
  Y}_{\gamma})^{\delta} \norm{Y-\widetilde Y}_\infty |t-s|^{\gamma \delta}
  \end{split}
 \end{equation*}
then
\begin{equation*}
  \begin{split}
\norm{B}_{(1+\delta)\gamma} \le  \norm{Y-\widetilde Y}_\gamma
\norm{\varphi}_{2,\delta} \norm{Y}_\gamma^{\delta}  + K\norm{\widetilde
Y}_\gamma  \|\varphi\|_{2,\delta} (\norm{Y}_{\gamma}+\norm{\widetilde
  Y}_{\gamma})^{\delta} \norm{Y-\widetilde Y}_\infty     
  \end{split}
\end{equation*}
and in the same way it is possible to obtain 
\begin{equation*}
\norm{B}_{\gamma} \le \|\varphi\|_{2,\delta}( \norm{Y-\widetilde Y}_\gamma
  + \norm{\widetilde
Y}_\gamma   \norm{Y-\widetilde Y}_\infty) .
\end{equation*}

Moreover
\begin{equation*}
  \norm{H'}_{\infty} \le \|\varphi\|_{2,\delta}
  \norm{Y'-\widetilde Y'}_{\infty} + \norm{Y'}_{\infty}
  \|\varphi\|_{2,\delta} \norm{Y-\widetilde Y}_{\infty}
\end{equation*}
\begin{equation*}
  \norm{H'}_{\gamma\delta} \le \|\varphi\|_{2,\delta}
  \norm{Y'-\widetilde Y'}_{\gamma\delta} + \norm{Y'}_{\gamma\delta}
  \|\varphi\|_{2,\delta} \norm{Y-\widetilde Y}_{\infty}
\end{equation*}
and
\begin{equation*}
  \begin{split}
  \norm{A}_{\gamma}& \le
 \|\varphi\|_{2,\delta}
  \norm{R-\widetilde R}_{\gamma} + \norm{R}_{\gamma}
  \|\varphi\|_{2,\delta} \norm{Y-\widetilde Y}_{\infty}    
\\ & \le
 \|\varphi\|_{2,\delta}
 ( \norm{Y'-\widetilde Y'}_{\infty}\norm{X}_\gamma+\norm{Y-\widetilde
 Y}_{\gamma}) + (\norm{Y'}_\infty \norm{X}_{\gamma}+\norm{Y}_{\gamma})
  \|\varphi\|_{2,\delta} \norm{Y-\widetilde Y}_{\infty}    
  \end{split}
\end{equation*}
\begin{equation*}
  \norm{A}_{(1+\delta)\gamma}\le
 \|\varphi\|_{2,\delta}
  \norm{R-\widetilde R}_{(1+\delta)\gamma} + \norm{R}_{(1+\delta)\gamma}
  \|\varphi\|_{2,\delta} \norm{Y-\widetilde Y}_{\infty}
\end{equation*}

And collecting all together these results we end up with
\begin{equation*}
 \norm{Z-\widetilde Z}_{D(X,\gamma,(1+\delta)\gamma)} \le C \norm{Y-\widetilde Y}_{D(X,\gamma,(1+\delta)\gamma)} 
\end{equation*}
with 
$$
C =  K
 \norm{\varphi}_{2,\delta} (1+\norm{X}_\gamma)
 (1+\norm{Y}_{D(X,\gamma,(1+\delta)\gamma)}+\norm{\widetilde Y}_{D(X,\gamma,(1+\delta)\gamma)})^{1+\delta}.
$$

To finish consider the case in which
 $\der\widetilde Y^\mu =  \widetilde Y^{\prime\,\mu}_{\nu}  \der
 \widetilde X^\nu
+ \widetilde  R^\mu_{\widetilde Y}$ is a path controlled by $\widetilde X$. If we let again
$\widetilde Z_t =
\varphi(\widetilde Y_t)$ and $H = Z -\widetilde Z$ we have 
\begin{equation*}
  \begin{split}
\der H^\mu &  =  \partial_\nu \varphi(\widetilde Y_\cdot)^\mu
\widetilde Y^{\prime\,\nu}_\kappa \der (X^\kappa - \widetilde X^\kappa )   +H^{\prime\,\mu}_{\kappa} \der X^\kappa + A^\mu + B^\mu
  \end{split}
\end{equation*}
where the only difference with the expression in
eq.~(\ref{eq:difference_H}) is in the first term in the r.h.s.
then
\begin{equation*}
  \|Z-\widetilde Z\|_\gamma +  \|Z'-\widetilde Z'\|_{\delta\gamma}
+ \|R_Z-R_{\widetilde Z}\|_{(1+\delta)\gamma} + \|Z'-\widetilde Z'\|_{\infty} \le C
(\epsilon+\|X-\widetilde X\|_\gamma)
\end{equation*}
with 
$$
\epsilon =   \|Y-\widetilde Y\|_\gamma +  \|Y'-\widetilde Y'\|_{\delta\gamma}
+ \|R_Y-R_{\widetilde Y}\|_{(1+\delta)\gamma} + \|Y'-\widetilde Y'\|_{\infty}
$$
and this concludes the proof of prop.~\ref{prop:functionD}.
\qed

\subsection{Some Proofs and Lemmata used in Sec.~\ref{sec:ode}}

\subsubsection{Proof of Lemma~\ref{eq:holder-patching}}
\label{sec:proof-holder-patching}
\proof
Take $ u \in I \cap J$:
\begin{equation*}
  \begin{split}
\sup_{t \in I\backslash J ,s \in J\backslash I}
  \frac{|X_{st}|}{|t-s|^\gamma} & \le
\sup_{t \in I\backslash J ,s \in J\backslash I}
  \frac{|X_{ut}|+|X_{su}|+|(\NN X)_{sut}|}{|t-s|^\gamma}  
\\ &  \le
\sup_{t \in I\backslash J ,s \in J\backslash I}
  \frac{|X_{ut}|}{|t-s|^\gamma}  
+\sup_{t \in I\backslash J ,s \in J\backslash I}
  \frac{|X_{su}|}{|t-s|^\gamma}      
+\sup_{t \in I\backslash J ,s \in J\backslash I}
  \frac{|(\NN X)_{sut}|}{|t-s|^\gamma}      
\\ &  \le
\sup_{t \in I\backslash J ,s \in J\backslash I}
  \frac{|X_{ut}|}{|t-u|^\gamma}  
+\sup_{t \in I\backslash J ,s \in J\backslash I}
  \frac{|X_{su}|}{|u-s|^\gamma}      
+\sup_{t \in I\backslash J ,s \in J\backslash I}
  \frac{|(\NN X)_{sut}|}{|t-u|^{\gamma_2} |s-u|^{\gamma_2}}      
\\ & \le
  \norm{X}_{\gamma,I}+\norm{X}_{\gamma,J}+\norm{X}_{\gamma_1,\gamma_2,I
  \cup J}
  \end{split}
\end{equation*}
then
\begin{equation*}
  \begin{split}
  \norm{X}_{\gamma,I \cup J} &= \sup_{t,s \in I \cup J}
  \frac{|X_{t}-X_s|}{|t-s|^\gamma} \le \sup_{t ,s \in I }
  \frac{|X_{t}-X_s|}{|t-s|^\gamma} + \sup_{t,s \in  J}
  \frac{|X_{t}-X_s|}{|t-s|^\gamma} + \sup_{t \in I\backslash J ,s \in
  J\backslash I}
  \frac{|X_{t}-X_s|}{|t-s|^\gamma} 
\\ & \le     2 (\norm{X}_{\gamma,I}+\norm{X}_{\gamma,J}) +\norm{X}_{\gamma_1,\gamma_2,I
  \cup J}
  \end{split}
\end{equation*}
as claimed.
\qed

\subsubsection{Lemmata for some bounds on the map $G$}
With the notation in the proof of Prop.~\ref{eq:existence_young} we have
\begin{lemma}
\label{lemma:zetabound-young}
For any interval $I=[t_0,t_0+T] \subseteq J$ such that $T<1$
the following bound holds
\begin{equation}
\label{eq:epsilonZbound-young-lemma}
  \begin{split}
  \epsilon_{Z,I} &  \le K  C_{X,I}
  C_{Y,I}^\delta [(1+\norm{\varphi}_{1,\delta})\rho_I + T^{\gamma\delta} \epsilon_{Y,I}]
  \end{split}
\end{equation}
\end{lemma}

\proof
Consider first the case when $\varphi =\wt \varphi$.
Eq.~(\ref{eq:continuity-young-x}) is a statement of continuity of the
The integral defined in Prop.~\ref{prop:young} is  a bounded bilinear
application $(A,B) \mapsto \int A dB$ then it is also continuous in
both arguments and it is easy to check that  
\begin{equation}
\label{eq:continuity-young-x}
  \norm{Q_Z-Q_{\wt Z}}_{(1+\delta)\gamma,I} \le K(C_{X,I}
  \epsilon^*_{W,I} + C_{Y,I} \rho_I)
\end{equation}
where we used the shorthands (defined in the proof of Prop.~\ref{eq:uniqueness_young}):
$$
\epsilon_{Z,I} = \norm{Z-\wt Z}_{\gamma,I}, \quad
\epsilon_{W,I}^* = \norm{W-\wt W}_{\delta\gamma,I}, \quad
\epsilon_{Y,I} = \norm{Y-\wt Y}_{\gamma,I}, \quad
\epsilon_{Y,I}^* = \norm{Y-\wt Y}_{\delta\gamma,I};
$$
$$
\rho_I = \norm{X-\wt X}_{\gamma,I} + |Y_0 - \wt Y_0|
$$
$$
C_{X,I} = \norm{X}_{\gamma,I} + \norm{\wt X}_{\gamma,I}
$$
$$
C_{Y,I} = \norm{Y}_{\gamma,I} + \norm{\wt Y}_{\gamma,I}
$$

Observe that
\begin{equation*}
  \begin{split}
\norm{\varphi(Y)-\varphi(\wt Y)}_{\infty,I} &\le
|\varphi(Y_0)-\varphi(\wt Y_0)| + T^{\delta\gamma}
\norm{\varphi(Y)-\varphi(\wt Y)}_{\delta\gamma,I}    
\\ & \le \norm{\varphi}_{1,\delta}\rho_I + T^{\delta\gamma} \epsilon^*_{W,I}
  \end{split}
\end{equation*}
\begin{equation*}
  \begin{split}
  \epsilon_{Z,I} & \le \norm{\varphi(Y) \der X-\varphi(\wt Y) \der \wt
    X}_{\gamma,I} + \norm{Q_Z-Q_{\wt Z}}_{\gamma,I}
\\ & \le \norm{\varphi(Y)-\varphi(\wt Y)}_{\infty,I}
  \norm{X}_{\gamma,I} + \norm{\varphi(\wt Y)}_{\infty,I} \norm{X-\wt
  X}_{\gamma,I} + T^{\delta\gamma} \norm{Q_Z-Q_{\wt
    Z}}_{(1+\delta)\gamma,I}    
\\ & \le \norm{\varphi}_{1,\delta} \rho_I C_{X,I}+T^{\delta\gamma}
    \epsilon^*_{W,I}+ K T^{\delta\gamma} (C_{X,I}
  \epsilon^*_{W,I} + C_{Y,I} \rho_I)
\\ & \le \norm{\varphi}_{1,\delta} \rho_I (C_{X,I}+1+KC_{Y,I}^\delta)
    + T^{\gamma\delta} \epsilon_{W,I}^* (C_{X,I}+KC_{Y,I})
  \end{split}
\end{equation*}
It remains to bound $\epsilon_{W,I}^*$:
Write
$$
\varphi(x)-\varphi(y) = \int_0^1 d\alpha \partial\varphi(\alpha x +
(1-\alpha)y) (x-y) = R\varphi(x,y)(x-y) 
$$
then
$$
\norm{R\varphi}_{\infty} = \sup_{x,y \in V} |R\varphi(x,y)| \le \norm{\varphi}_{1,\delta}
$$
and
\begin{equation*}
  \begin{split}
|R\varphi(x,y)-R\varphi(x',y')| & = \left|\int_0^1 (\partial \varphi(\alpha x +
(1-\alpha)y)-\partial \varphi(\alpha x' +
(1-\alpha)y') d\alpha \right|
    \\ & \le \norm{\varphi}_{1,\delta}
\int_0^1 |\alpha (x-x') +
(1-\alpha)(y-y')|^\delta d\alpha 
    \\ & \le \norm{\varphi}_{1,\delta}(|x-x'|^\delta+|y-y'|^\delta)
  \end{split}
\end{equation*}
so that
\begin{equation*}
  \begin{split}
    \epsilon^*_{W,I}  & = \norm{\varphi(Y)-\varphi(\wt
  Y)}_{\delta\gamma,I} 
=     \norm{R\varphi(Y,\wt Y)(Y-\wt Y)}_{\delta\gamma,I}
\\ & \le \norm{R\varphi(Y,\wt Y)}_{\infty,I} \norm{Y-\wt
  Y}_{\delta\gamma,I} + \norm{R\varphi(Y,\wt
  Y)}_{\delta\gamma,I}\norm{Y-\wt Y}_{\infty,I}
\\ & \le \norm{\varphi}_{1,\delta} \norm{Y-\wt
  Y}_{\delta\gamma,I} + \norm{Y-\wt Y}_{\infty,I}
  \norm{\varphi}_{1,\delta} (\norm{Y}^\delta_{\gamma,I} + \norm{\wt
  Y}^\delta_{\gamma,I})
\\ & \le K \norm{\varphi}_{1,\delta} C_{Y,I}^\delta \epsilon_{Y,I}^*
\\ & \le K \norm{\varphi}_{1,\delta} C_{Y,I}^\delta \epsilon_{Y,I}
  \end{split}
\end{equation*}
concluding:
\begin{equation}
\label{eq:epsilonZbound-young-00}
  \begin{split}
  \epsilon_{Z,I} &  \le K \norm{\varphi}_{1,\delta} C_{X,I}
  C_{Y,I}^\delta (\rho_I + T^{\gamma\delta} \epsilon_{Y,I})
  \end{split}
\end{equation}

The general case in which $\varphi \neq \wt\varphi$ can be easily
derived from Eq.~(\ref{eq:epsilonZbound-young-00}) and the continuity
of the integral, giving:
\begin{equation*}
  \begin{split}
  \epsilon_{Z,I} &  \le K  C_{X,I}
  C_{Y,I}^\delta [(1+\norm{\varphi}_{1,\delta})\rho_I + T^{\gamma\delta} \epsilon_{Y,I}].
  \end{split}
\end{equation*}
\qed

Using the notation in the proof of Prop.~\ref{eq:existence_rough} we have
\begin{lemma}
\label{lemma:zetabound-rough}
For any interval $I=[t_0,t_0+T] \subseteq J$ such that $T<1$
the following bound holds
\begin{equation}
\label{eq:zetabound-lemma}
   \epsilon_{Z,I} \le K  \norm{\varphi}_{2,\delta}(C_{X,I}^2 C_{Y,I}^3
   \rho_I + T^{\delta\gamma} C_{X,I}^3 C_{Y,I}^2 \epsilon_{Y,I}) + K 
\norm{\varphi-\wt \varphi}_{2,\delta} C_{X,I} C_{Y,I}^2
 \end{equation}  
\end{lemma}

\proof To begin assume that $\varphi = \wt \varphi$.
Let  $W = \varphi(Y)$, $\widetilde W = \varphi(\widetilde Y)$
and
write their decomposition as
$$
\der W^\mu = W^{\prime\,\mu}_{\nu} \der X^\nu + R^\mu_W, \qquad 
\der \widetilde W^\mu = \widetilde W^{\prime\,\mu}_{\nu} \der \widetilde X^\nu + R^\mu_{\widetilde W}, \qquad 
$$
with $W^{\prime\,\mu}_{\nu} = \partial_\kappa \varphi(Y)^{\mu}
Y^{\prime\,\kappa}_{\nu}$, $\widetilde W^{\prime\,\mu}_{\nu} = \partial_\kappa
\varphi(\widetilde Y)^{\mu} \widetilde Y^{\prime\,\kappa}_{\nu}$ 
Moreover let
\begin{equation*}
  \epsilon_{W,I}^* = \norm{W'-\wt W'}_{\infty,I} +  \norm{W'-\wt
  W'}_{\delta\gamma,I} + \norm{R_W+R_{\wt W}}_{(1+\delta)\gamma,I} +  \norm{W -\wt
  W}_{\gamma,I} 
\end{equation*}

Using the bound~(\ref{eq:betterbound-difference})
we have
\begin{equation}
  \label{eq:boundQ}
  \norm{Q- Q_{\wt Z}}_{(2+\delta)\gamma} \le  K (D_1+D_2)
\end{equation}
$$
D_1 = C_X \epsilon^*_{W,I}
$$
\begin{equation*}
  \begin{split}
D_2  & =
(\norm{\varphi(Y)}_{D(X,\gamma,(1+\delta)\gamma),I}+\norm{\varphi(\wt Y)}_{D(\widetilde X,\gamma,(1+\delta)\gamma),I})(\norm{X-\widetilde X}_{\gamma,I}  +
\norm{\mathbb{X}^2-\mathbb{\widetilde X}^2}_{2\gamma,I}) 
\\ & \le  K
\norm{\varphi}_{2,\delta} C_{Y,I}^2 \rho_I    
  \end{split}
\end{equation*}
where we used eq.~(\ref{eq:boundWxx}) to bound
$\norm{\varphi(Y)}_{D(X,\gamma,(1+\delta)\gamma),I}$ and
$\norm{\varphi(\wt Y)}_{D(\wt X,\gamma,(1+\delta)\gamma),I}$ in terms
of $C_{Y,I}$.

By Prop.~\ref{prop:functionD} we have
\begin{equation}
  \label{eq:function_difference-xx}
\epsilon^*_{W,I}  \le K
 \norm{\varphi}_{2,\delta} C_{X,I}
 C_{Y,I}^{1+\delta} 
(\|X-\widetilde X\|_{\gamma,I} + \epsilon^*_{Y,I}) \le K
 \norm{\varphi}_{2,\delta} C_{X,I}
 C_{Y,I}^{2} 
(\rho_I + \epsilon^*_{Y,I})   
\end{equation}
with 
\begin{equation*}
  \begin{split}
\epsilon^*_{Y,I} & = \|Y'-\widetilde Y'\|_\infty + \|Y'-\widetilde
Y'\|_{\delta\gamma} +
\|R_Y-R_{\widetilde Y}\|_{(1+\delta)\gamma} +\|Y-\widetilde Y\|_\gamma
  \end{split}
\end{equation*}
and
$$
C_{I} =  K
 \norm{\varphi}_{2,\delta} C_{X,I}
 C_{Y,I}^{1+\delta} 
$$

Taking $T < 1$ we can bound $\epsilon^*_{Y,I} \le \epsilon_{Y,I} +
\norm{Y-\wt Y}_{\gamma,I}$ and
\begin{equation}
\label{eq:bound-epsilon-Y-xx}
  \begin{split}
\epsilon^*_{Y,I} & \le \|Y'-\widetilde Y'\|_\infty + \|Y'-\widetilde
Y'\|_{\gamma} +
\|R_Y-R_{\widetilde Y}\|_{2\gamma} + C_{X,I} \epsilon_{Y,I} + C_{Y,I}
\norm{X-\wt X}_{\gamma,I}
\\ & \le 2 C_{X,I} \epsilon_{Y,I} + C_{Y,I}
\rho_I
  \end{split}
\end{equation}
where we used the following majorization for $\norm{Y-\wt Y}_{\gamma,I}$:
\begin{equation}
  \label{eq:boundYgamma}
  \begin{split}
\norm{Y-\wt Y}_{\gamma,I} & \le \norm{Y'\der X-\wt
  Y' \der \wt X}_{\gamma,I} +\norm{R_Y-R_{\wt Y}}_{\gamma,I} 
  \\ & \le 
\norm{Y'-\wt
  Y'}_{\infty,I} \norm{X}_{\gamma,I} + (\norm{Y'}_{\infty,I}+\norm{\wt
  Y'}_{\infty,I}) \norm{X-\wt X}_{\gamma,I} + \norm{R_Y-R_{\wt Y}}_{2\gamma,I} 
  \\ & \le C_{X,I} \epsilon_{Y,I} + C_{Y,I} \rho_I
  \end{split}
\end{equation}
Eq.~(\ref{eq:bound-epsilon-Y-xx}) together with
eq.~(\ref{eq:function_difference-xx}) imply
\begin{equation*}
  \epsilon^*_{W,I} \le  K
 \norm{\varphi}_{2,\delta} (C_{X,I}
 C_{Y,I}^{3}  \rho_I + 
  C^2_{X,I}
 C_{Y,I}^{2} \epsilon_{Y,I}) 
\end{equation*}
and so
\begin{equation}
  \label{eq:boundQ2}
  \begin{split}
\norm{Q_Z- Q_{\wt Z}}_{(2+\delta)\gamma} & \le  K C_X \epsilon_{W,I} + K
  \norm{\varphi}_{2,\delta} C_{Y,I}^2 \rho_{I}
   \\ & \le
 K (C_X C_I (1+2 C_Y)+\norm{\varphi}_{2,\delta} C_{Y,I}^2) \rho_I + 2K C_I C^2_X \epsilon_{Y,I}
   \\ & \le K \norm{\varphi}_{2,\delta}(C_X^2 C_Y^3 \rho_I + C_X^3
  C_Y^2 \epsilon_{Y,I})
  \end{split}
\end{equation}

\begin{equation*}
  \begin{split}
\epsilon_{Z,I} & = \norm{\varphi(Y)-\varphi(\wt Y)}_{\infty,I} +
\norm{\varphi(Y)-\varphi(\wt Y)}_{\gamma,I} + \norm{R_Z-R_{\wt
    Z}}_{2\gamma,I}    
\\ & \le   |\varphi(Y_0)-\varphi(\wt Y_0)|+
2 \norm{\varphi(Y)-\varphi(\wt Y)}_{\gamma,I} + \norm{R_Z-R_{\wt
    Z}}_{2\gamma,I}
  \end{split}
\end{equation*}

Proceed step by step:
\begin{equation*}
  \begin{split}
  \norm{\partial \varphi(Y)-\partial \varphi(\wt Y)}_{\infty,I} 
 & \le  |\partial \varphi(Y_{t_0})-\partial \varphi(\wt Y_{t_0})|+
 T^\gamma \norm{\partial \varphi(Y)-\partial \varphi(\wt Y)}_{\gamma,I}
 \\ & \le \norm{\varphi}_{2,\delta} |Y_{t_0}-\wt Y_{t_0}| +
 T^\gamma \norm{\varphi}_{2,\delta} \norm{Y-\wt Y}_{\gamma,I}
\\ & \le 
 T^\gamma \norm{\varphi}_{2,\delta} C_{X,I} \epsilon_{Y,I} +
  2 \norm{\varphi}_{2,\delta} C_{Y,I} \rho_I
  \end{split}
\end{equation*}
Next:
\begin{equation*}
  \begin{split}
 \norm{R_Z-R_{\wt Z}}_{2\gamma,I} & \le  \norm{\partial
 \varphi(Y)\mathbb{X}^2-\partial \varphi(\wt Y)\mathbb{\wt X}^2
 }_{2\gamma,I} 
+  \norm{Q_Z-Q_{\wt Z}}_{2\gamma,I}  
\\ & \le
 \norm{\partial
 \varphi(Y)-\partial \varphi(\wt Y)
 }_{\infty,I} (\norm{\mathbb{X}^2}_{2\gamma,I}+\norm{\mathbb{\wt
 X}^2}_{2\gamma,I} ) 
\\ &\qquad + (\norm{\partial
 \varphi(Y)}_{\infty,I}+\norm{\partial \varphi(\wt Y)}_{\infty,I})
\norm{\mathbb{X}^2-\mathbb{\wt
 X}^2}_{2\gamma,I}
\\ &\qquad+ T^{\delta\gamma} \norm{Q_Z-Q_{\wt Z}}_{(2+\delta)\gamma,I}  
\\ & \le K \norm{\varphi}_{2,\delta}(\rho_I C_X^2 C_Y^3 + \epsilon_{Y,I} T^{\delta \gamma}
 C_X^3 C_Y^2)
  \end{split}
\end{equation*}
and
\begin{equation*}
  \begin{split}
     \norm{\varphi(Y)&-\varphi(\wt Y)}_{\gamma,I}  \le
\norm{\partial \varphi(Y) \der X - \partial \varphi(\wt Y)\der \wt
  X}_{\gamma,I} + \norm{R_W-R_{\wt W}}_{\gamma,I}
\\ & \le  
\norm{\partial \varphi(Y)  - \partial \varphi(\wt
  Y)}_{\infty,I}(\norm{X}_{\gamma,I}+\norm{\wt X}_{\gamma,I})
\\ & \qquad +
(\norm{\partial \varphi(Y) }_{\infty,I}+\norm{\partial \varphi(\wt Y)
}_{\infty,I}) \norm{X-\wt X}_{\gamma,I}+ T^\gamma \norm{R_W-R_{\wt
    W}}_{2\gamma,I}
\\ & \le
 (\norm{\varphi}_{2,\delta} |Y_{t_0}-\wt Y_{t_0}| +
 T^\gamma \norm{\varphi}_{2,\delta} C_{X,I} \epsilon_{Y,I} +
 \norm{\varphi}_{2,\delta} C_{Y,I} \norm{X-\wt
   X}_{\gamma,I})(\norm{X}_{\gamma,I}+\norm{\wt X}_{\gamma,I})
\\ & \qquad +2
 \norm{\varphi}_{2,\delta}  \norm{X-\wt X}_{\gamma,I}+ T^\gamma
 \epsilon_{W,I}^*
\\ & \le K \norm{\varphi}_{2,\delta}(C_{X,I} C_{Y,I}^3 \rho_I + T^\gamma C_{X,I}^2 C_{Y,I}^2 \epsilon_{Y,I})
  \end{split}
\end{equation*}
finally we have
\begin{equation}
\label{eq:zetabound00}
  \epsilon_{Z,I} \le K  \norm{\varphi}_{2,\delta}(C_{X,I}^2 C_{Y,I}^3 \rho_I + T^{\delta\gamma} C_{X,I}^3 C_{Y,I}^2 \epsilon_{Y,I}).
 \end{equation}

When $\varphi \neq \wt \varphi$ rewrite the difference $Z-\wt Z$ as
\begin{equation*}
 Z_t-\wt Z_t = Y_{t_0} - \wt Y_{t_0} + \int_{t_0}^t
 [\varphi(Y)-\varphi(\wt Y)] dX + \int_{t_0}^t
 [\varphi(\wt Y)-\wt \varphi(\wt Y)] dX  
\end{equation*}
the contribution to $\epsilon_{Z,I}$ from the first integral is
bounded by  Eq.~(\ref{eq:zetabound00}) while the last integral can be
bounded by $K \norm{\varphi-\wt\varphi}_{2,\delta} C_{X,I} C_{Y,I}^2$
(cfr. Eq.~(\ref{eq:bound_on_gy})) giving the final result~(\ref{eq:epsilonZbound-young-lemma}).
\qed

\subsection{Proof of Lemma~\ref{lemma:besov}}
\label{sec:proof_besov}

\proof
 Let $B(u,r) = \{ w \in T : |w-u| \le r \}$.
Observe that by the monotonicity and convexity of $\psi$ for any
couple of measurable sets $A,B \subset T$ we have
\begin{equation}
\label{eq:mean-estimate}
  \begin{split}
  \left|\int_{A\times B} R_{st} \frac{dt ds}{|A| |B|}\right|    
&  \le p(d(A,B)/4) \psi^{-1}\left(\int_{A\times B}
  \psi\left(\frac{|A_{st}|}{p(d(t,s)/4)}\right) \frac{dt ds}{|A| |B|}\right)
\\ & \qquad \le p(d(A,B)/4) \psi^{-1}\left(\frac{U}{|A| |B|}\right)
  \end{split}
\end{equation}
where $d(A,B) = \sup_{t \in A, s \in B}  |t-s|$.
Let
\begin{equation*}
  \overline R(t,r_1,r_2) = \int_{B(t,r_1)}   \frac{du}{| B(t,r_1)|
  }\int_{ B(t,r_2)}   \frac{ dv}{| B(t,r_2)|} R_{uv}
\end{equation*}
Take $t,s \in T$, $a = |t-s|$, define the decreasing sequence of numbers
$\lambda_n \downarrow 0$ as $\lambda_0 = a$, $\lambda_{n+1}$ such that
$$
p( \lambda_{n}) = 2 p( \lambda_{n+1})
$$
then
\begin{equation*}
  \begin{split}
p((\lambda_n+\lambda_{n+1})/4) & \le p(\lambda_n) = 2 p(\lambda_{n+1}) 
\\ & = 4    p(\lambda_{n+1}) - 2 p(\lambda_{n+1})
\\ & = 4 [   p(\lambda_{n+1})  - p(\lambda_{n+2})].
  \end{split}
\end{equation*}
Using eq.~(\ref{eq:mean-estimate}) and the fact that $|B(t,\lambda_i)|
\ge \lambda_i$ for every $i \ge 0$ we have
\begin{equation*}
  \begin{split}
|\overline R(t,\lambda_{n+1},\lambda_{n})| &   
\le p((\lambda_{n} +\lambda_{n+1})/4) 
 \psi^{-1}\left(\frac{U}{ \lambda_{n} \lambda_{n+1}}\right)   
\\ & \le 4 [   p(\lambda_{n+1})  - p(\lambda_{n+2})]
 \psi^{-1}\left(\frac{U}{ \lambda_{n} \lambda_{n+1}}\right)   
\\ & \le 4 \int_{\lambda_{n+2}}^{\lambda_{n+1}} 
 \psi^{-1}\left(\frac{U}{ r^2}\right)   dp(r).
  \end{split}
\end{equation*}
Take a sequence $\{ t_i \}_{i=0}^\infty$ of variables in $T$ and note that, for
every $n \ge 0$, 
$$
R_{t\,t_n} = R_{t\, t_{n+1}} + R_{t_{n+1}\, t_{n}} + (\NN R)_{t t_{n+1} t_n}
$$
so that, by induction,
$$
R_{t\, t_0} = R_{t\, t_{n+1}}+ \sum_{i=0}^{n} [R_{t_{i+1} t_i} + (\NN R)_{t\, t_{i+1} t_i}  ].
$$

Average each $t_i$ over the ball $B(t,\lambda_i)$ and bound
as follows
\begin{equation}
\label{eq:finite-rep}
\overline R(t,0,\lambda_0) = \overline R(t,0,\lambda_{n+1}) 
+ \sum_{i=0}^{n} \overline R(t,\lambda_{i+1},\lambda_{i}) +
\sum_{i=0}^{n} \overline B(t,\lambda_{i+1},\lambda_{i}) 
\end{equation}
where
\begin{equation*}
\overline B(t,\lambda_{i+1},\lambda_i) = \int_{B(t,\lambda_{i+1})}
\frac{dv}{| B(t,\lambda_{i+1}) |}  \int_{B(t,\lambda_i)}
\frac{du}{| B(t,\lambda_i)|} \NN
R_{tvu}   
\end{equation*}
which, using~(\ref{eq:ext-bound-n}), can be majorized by
\begin{equation*}
  \begin{split}
 | \overline B(t,\lambda_{i+1},\lambda_i)| & \le \psi^{-1}\left(
 \frac{C}{  \lambda_i^2}\right) p( \lambda_i/2) \le  4 \psi^{-1}\left(
 \frac{C}{ \lambda_i^2}\right) [p(\lambda_{i+1})-p(\lambda_{i+2})]
\\ & \le  4 \int_{\lambda_{i+2}}^{\lambda_{i+1}} \psi^{-1}\left(
 \frac{C}{ r^2}\right) dp(r)    
  \end{split}
\end{equation*}

Then, taking the limit as $n \to \infty$ in Eq.~(\ref{eq:finite-rep}), using the continuity of
$R$ and that $R_{tt}=0$, we get
\begin{equation}
\label{eq:bound-from-t}
  \begin{split}
|\overline R(t,0,\lambda_0)| & \le 
 \sum_{i=0}^{\infty} 4 \int_{\lambda_{i+2}}^{\lambda_{i+1}} 
 \psi^{-1}\left(\frac{U}{ r^2}\right)   dp(r) 
+ \sum_{i=0}^{\infty} 4 \int_{\lambda_{i+2}}^{\lambda_{i+1}} \psi^{-1}\left(
 \frac{C}{ r^2}\right) dp(r)    
\\ &    \le 4 \int_0^{\lambda_1} \left[ \psi^{-1}\left(\frac{U}{ r^2}\right)+ \psi^{-1}\left(\frac{C}{ r^2}\right)\right]   dp(r) 
\\ &    \le 4 \int_0^{|t-s|} \left[ \psi^{-1}\left(\frac{U}{ r^2}\right)+ \psi^{-1}\left(\frac{C}{ r^2}\right)\right]   dp(r) 
  \end{split}
\end{equation}
and of course the analogous estimate
\begin{equation}
\label{eq:bound-from-s}
  |\overline R(s,0,\lambda_0)|  \le  4 \int_0^{|t-s|} \left[ \psi^{-1}\left(\frac{U}{ r^2}\right)+ \psi^{-1}\left(\frac{C}{ r^2}\right)\right]   dp(r) 
\end{equation}

Moreover
$$
R_{st} = R_{s u} + R_{uv} + R_{v t} + \NN R_{sut} + \NN R_{uvt} 
$$
so
$$
|R_{st}| \le |R_{s u}| + |R_{v t}| + |R_{uv}| + \sup_{r \in [s,t]}
|\NN R_{srt}| + \sup_{r \in [u,t]} |\NN R_{urt}|
$$
by averaging $u$ over the ball $B(s,a)$ and $v$ over the ball
$B(t,a)$ we get 
\begin{equation*}
\int_{B(s,a)} \frac{du}{|B(s,a)|} \int_{B(t,a)} \frac{dv}{|B(t,a)|}   |R_{uv}|
\le p(3 a/4)\psi^{-1}\left(\frac{U}{4 a^2}\right) \le \int_0^{|t-s|}
\psi^{-1}\left(\frac{U}{ r^2}\right) dp(r) 
\end{equation*}
and
\begin{equation*}
\int_{B(s,a)} \frac{du}{|B(s,a)|} \sup_{r \in [u,t]} |\NN R_{urt}| \le 
p(a/2)\psi^{-1}\left(\frac{C}{  a^2}\right)  \le
\int_0^{|t-s|} \psi^{-1}\left(\frac{C}{ r^2}\right) dp(r) 
\end{equation*}

Putting all toghether we end up with
\begin{equation*}
 |R_{st}| \le  
10 \int_0^{|t-s|} \left[ \psi^{-1}\left(\frac{U}{ r^2}\right)+
  \psi^{-1}\left(\frac{C}{ r^2}\right)\right]   dp(r)
\end{equation*}
\qed

\end{document}